\documentclass[12pt]{amsart}
\usepackage{amsmath,amsfonts,amssymb,amsthm}
\usepackage[dvips]{graphicx}
\usepackage{psfrag}
\usepackage{epsfig,subfigure}
\usepackage{color}
\usepackage[T1]{fontenc}
\usepackage{a4wide}

\def\struckint{\mathop{%
\def\mathpalette##1##2{\mathchoice{##1\displaystyle##2}%
  {##1\textstyle##2}{##1\scriptstyle##2}{##1\scriptscriptstyle##2}}%
\mathpalette
{\vbox\bgroup\baselineskip0pt\lineskiplimit-1000pt\lineskip-1000pt
\halign\bgroup\hfill$}
{##$\hfill\cr{\intop}\cr\diagup\cr\egroup\egroup}%
}\limits}

\newcommand{\gp}{\pi}

\newcommand\D{d\hspace{-0.5pt}}
\newcommand{\Dq}{q}
\newcommand{\HA}{\mathcal{H}}

\newcommand{\QQQ}{\mathcal{Q}}

\newcommand{\HHH}{\mathcal{H}}

\theoremstyle{definition} 
\newtheorem{defs}{Definition}[section]

\newtheorem{defslem}[defs]{Definition/Lemma}

\theoremstyle{plain} 
\newtheorem{prop}[defs]{Proposition}
\newtheorem{lem}[defs]{Lemma}
\newtheorem{ths}[defs]{Theorem} 
\newtheorem{conv}[defs]{Convention}

\newtheorem*{ThmA}{Theorem A}
\newtheorem*{ThmB}{Theorem B}
\newtheorem*{ThmC}{Theorem C}
\newtheorem*{ThmD}{Theorem D}

\newtheorem{cor}[defs]{Corollary}
\newtheorem*{NNths}{Theorem}
\newtheorem*{NNprop}{Proposition}

\theoremstyle{remark} 
\newtheorem{rem}[defs]{Remark}

\newtheorem{Example}[defs]{Example}

\newtheorem*{Acknowledgments}{Acknowledgments}
\newtheorem*{Reader's guide}{Reader's guide}

\begin{document}
\title[Rauzy-Veech induction for quadratic differentials]
{Dynamics and geometry of the Rauzy-Veech induction for quadratic
  differentials}
\author{Corentin Boissy, Erwan Lanneau}
\address{
IRMAR, Campus de Beaulieu, UMR CNRS 6625 \newline
Universit\'e de Rennes I \newline
35042 Rennes cedex, France
}
\email{corentin.boissy@univ-rennes1.fr}

\address{
Centre de Physique Th\'eorique (CPT), UMR CNRS 6207 \newline
Universit\'e du Sud Toulon-Var and \newline
F\'ed\'eration de Recherches des Unit\'es de 
Math\'ematiques de Marseille \newline
Luminy, Case 907, F-13288 Marseille Cedex 9, France
}

\email{lanneau@cpt.univ-mrs.fr}

\subjclass[2000]{Primary: 37E05. Secondary: 37D40}
\keywords{Interval exchange map, Linear involution, Rauzy-Veech induction, Quadratic
  differential, Moduli space, Teichm\"uller geodesic flow}

\date{January 28, 2008}

\begin{abstract}
Interval exchange maps are related to geodesic flows
on translation surfaces; they correspond to the first return maps of the
vertical flow on a transverse segment. The Rauzy-Veech induction
on the space of interval exchange maps provides a powerful tool
to analyze the Teichm\"uller geodesic flow on the moduli space of
Abelian differentials. Several major results have been proved using
this renormalization. 

Danthony and Nogueira introduced in $1988$ a natural generalization of 
interval exchange transformations, namely the linear
involutions. These maps are related to general measured foliations on
surfaces (orientable or not). In this paper we are interested by such
maps related to geodesic flow on (orientable) flat surfaces with
$\mathbb{Z}/2\mathbb{Z}$ linear holonomy. 
We relate geometry and dynamics of such maps to the combinatorics of
generalized permutations. We study an analogue of the Rauzy-Veech
induction and give an efficient combinatorial characterization
of its attractors. We establish a natural bijection between
the extended Rauzy classes of generalized permutations and
connected components of the strata of meromorphic quadratic
differentials with at most simple poles, which allows, in
particular, to classify the connected components of all
exceptional strata.
\end{abstract}
\setcounter{tocdepth}{1}

\maketitle
\tableofcontents
%

\section*{Introduction}

A geodesic flow in a given direction on a translation surface induces
on a transverse segment an interval exchange map. Dynamic of such
transformations has been extensively studied during these last thirty
years providing applications to billiards in rational polygons, to
measured foliations on surfaces, to Teichm\"uller geometry and
dynamics, etc. \medskip

Interval exchange transformations are closely related to Abelian
differentials on Riemann surfaces. It is well known that the continued
fractions encode cutting sequences of hyperbolic geodesics on the 
Poincar\'e upper half-plane. Similarly, the Rauzy-Veech induction
(analogous to Euclidean algorithm) provides a discrete model for the 
Teichm\"uller geodesics flow (\cite{Rauzy, Ve82, Arnoux}). \medskip

Using this relation H.~Masur in~\cite{Ma82} and W.~A.~Veech
in~\cite{Ve82} have independently proved the Keane's conjecture
(unique ergodicity of almost all interval exchange
transformations). Using combinatorics of 
Rauzy classes, Kontsevich and Zorich classified the connected
components of strata of the moduli spaces of Abelian
differentials~\cite{Kontsevich:Zorich}. More recently, Avila,
Gou\"ezel and Yoccoz proved the exponential decay of correlations for
the Teichm\"uller geodesic flow also using a renormalization of the
Rauzy-Veech induction (see~\cite{Zorich,Avila:Gouezel:Yoccoz}). Avila
and Viana used combinatorics of Rauzy-Veech induction to prove the
simplicity of the essential part of the Lyapunov spectrum of the
Teichm\"uller geodesic flow on the strata of Abelian differentials 
(see~\cite{Avila:Viana}). Recently Bufetov and Gurevich proved the 
existence and uniqueness of the measure of maximal entropy for the
Teichm\"uller geodesic flow on the moduli space of Abelian
differentials~\cite{Bu:Gu}. Avila and Forni proved the weak mixing for
almost all interval exchange transformations and translation
flows~\cite{Avila:Forni}.

These examples show that Rauzy-Veech induction which was initially
elaborated to prove ergodicity of interval exchange transformations
and ergodicity of the Teichm\"uller geodesic flow is, actually, very
efficient far beyond these initial problems.

However, all the aforementioned results concern only the moduli space
of Abelian differentials. The corresponding questions for strata of
strict quadratic differentials (i.e. of those, which are not global
squares of Abelian differentials) remain open. 

Note that the (co)tangent bundle to the moduli space of curves is
naturally identified with the moduli space of \emph{quadratic}
differentials. From this point of view, the strata of Abelian differentials
represent special orbifolds of high codimension in the total space of
the tangent bundle. Our interest in Teichm\"uller dynamics and
geometry of the strata of strict quadratic differentials was one of
the main motivations for developing Rauzy-Veech induction for
quadratic differentials. \medskip

Natural generalizations of interval exchange transformations were
introduced by Danthony and Nogueira in~\cite{DaNo:CRAS,DaNo} (see
also~\cite{No}) as cross sections of measured foliations on
surfaces. They introduced the notion of linear involutions, as well as the 
notion of Rauzy induction on these maps. \medskip

Studying Lyapunov spectrum of the Teichm\"uller geodesic flow
Kontsevich and Zorich have performed series of computer experiments
with linear involutions corresponding to
quadratic differentials~\cite{Kontsevich:Zorich:hodge}. These experiments indicated
appearance of attractors for the Rauzy-Veech induction, as well as
examples of generalized permutations such that the corresponding
linear involutions are minimal for a domain of parameters of
positive measure, and non minimal for a complementary domain of
parameters also of positive measure (examples of this type are
presented in Figure~\ref{rauzy:class:2ii} and Figure~\ref{fig:rauzy:2}
in Appendix $A$). But at this point, there was no combinatorial explanation.

Thus, in order to generalize technique of Rauzy-Veech induction to
quadratic differentials in a consistent way it was necessary to find
combinatorial criteria allowing to identify generalized permutations,
which belong to attractors and those ones, which represents cross
sections of vertical foliations of quadratic differentials. It was
also necessary to distinguish those generalized permutation which give
rise to minimal linear involution, and to
specify the domains of appropriate parameters. \medskip

In this paper we establish corresponding combinatorial criteria,
which enable us to develop technique of Rauzy-Veech induction for
quadratic differentials. Partial results in this direction were
obtained by the second author in~\cite{La1}. We also study relations
between combinatorics, geometry and dynamics of linear involutions.

To compare similarities and differences between linear involutions 
corresponding to Abelian and to quadratic
differentials let us first briefly review the situation in the
classical case. \medskip

An interval exchange transformation is encoded by a combinatorial data
(permutation $\pi$ on $d$ elements) and by a continuous data (lengths
$\lambda_1, \dots, \lambda_d$ of the intervals). Recall that the
Keane's property (see below) is a criterion of ``irrationality'' (which, in
particular, implies minimality) of an interval exchange
transformation. This property is satisfied for almost all parameters
$\lambda$ when the permutation $\pi$ is irreducible (i.e. $\pi(\{1,
\dots , k\}) \neq \{1, \dots , k\},\quad 1\le k \le d - 1)$, while
when $\pi$ is reducible, the corresponding interval exchange map is
\textit{never} minimal. On the other hand the irrational interval
exchange maps are precisely those that arise as cross sections of
minimal vertical flows on well chosen transverse intervals. \medskip 

The Rauzy-Veech induction consists in taking the first return map of
an interval exchange transformation to an appropriate smaller
interval. This induction can be viewed as a dynamical
system on a finite-dimensional space of interval exchange
maps. The behavior of an orbit of the induction provides
important information on dynamics of the interval exchange
transformation representing the starting point. This information is
especially useful when all iterates are well defined and when the
length of the underlying subintervals tends to zero. An interval
exchange transformation satisfying the latter conditions is said to
have \emph{Keane's property}. For a given irreducible permutation
$\pi$, the subset of parameters $\lambda$ which give rise to interval
exchange transformations satisfying Keane's property contains all
irrational parameters, and so it is a full Lebesgue measure
subset. Moreover, for the space of interval exchange transformations
with irreducible combinatorial data, the renormalized induction
process is recurrent with respect to the Lebesgue measure (and even
ergodic by a theorem of Veech). Note that the corresponding invariant
measure has infinite total mass. \medskip

In this paper we use the definition of {\it linear
  involution}~\footnote{Let $f$ be the involution of $X\times\{0,1\}$  
given by $f(x, \varepsilon)=(x,1-\varepsilon)$.
A linear involution is a map $T$, from $X\times\{0,1\}$ into itself, of the form 
$f\circ \tilde T$, where $\tilde T$ is an involution of $X\times\{0,1\}$ 
without fixed point, continuous except in finitely many points, 
and which preserves the Lebesgue measure. In this paper we will only
consider linear involutions with the additional condition. The
derivative of $\tilde T$ is $-1$ at $(x,\varepsilon)$ if
$(x,\varepsilon)$ and $T(x,\varepsilon)$ belong to the same connected
component, and $-1$ otherwise; see also
Convention~\ref{convention:oriented}.} proposed by Danthony and
Nogueira (see~\cite{DaNo:CRAS,DaNo}).

As above, a \emph{linear involution} is encoded by a
combinatorial data (``generalized permutation'') and by continuous
data. A generalized permutation of type $(l,m)$ (with $l + m = 2d$) is
a two-to-one map $\pi : \{1,\dots,2d\} \to \mathcal{A}$ to an
alphabet $\mathcal{A}$. \medskip

A generalized permutation is called \emph{irreducible} if there exists
a linear involution associated to this generalized
permutation, which represents an appropriate cross section of the
vertical foliation of some quadratic differential. A generalized
permutation is called \emph{dynamically irreducible} if there exists a
minimal linear involution associated to this
generalized permutation. It is easy to show that any irreducible
generalized permutation is dynamically irreducible; the converse is
not true in general as we will see.

\begin{ThmA}
Irreducible and dynamically irreducible generalized permutations can
be characterized by natural criteria expressed in elementary
combinatorial terms.
\end{ThmA}

The corresponding criteria are stated as Definitions~\ref{def:irred}
and Definition~\ref{def:dyn:irr} respectively.

Consider a dynamically irreducible generalized permutation
$\pi$. The parameter space of normalized linear involutions
associated to $\pi$ is represented by a hyperplane section of a
simplex. We describe an explicit procedure which associates to each 
generalized permutation $\pi$ an open subset in the parameter space
defined by a system of linear inequalities determined by $\pi$. This
subset is called the set of \textit{admissible parameters}. When $\pi$
is irreducible, the set of admissible parameters coincides with entire
parameter space; in general it is smaller. The next result gives a
more precise statement than Theorem~$A$ in the dynamically irreducible
case.

\begin{ThmB}
$\ $
\begin{enumerate}
\item If $\pi$ is not dynamically irreducible, or if $\pi$ is dynamically
irreducible, but $\lambda$ does not belong to the set of admissible
parameters, the linear involution $T=(\pi,\lambda)$ is not minimal.
   
\item
If $\pi$ is dynamically irreducible, then for almost all admissible
parameters $\lambda$  the linear involution 
$T=(\pi,\lambda)$ satisfies the Keane's property, and hence is
minimal.
\end{enumerate}
\end{ThmB}

Since the Rauzy-Veech induction commutes with dilatations, it
projectivizes to a map $\mathcal{R}_r$ on the space of normalized
linear involutions; we shall call this map the 
\textit{renormalized Rauzy-Veech induction}.

\begin{ThmC}
Let T be a linear involution on the unit interval and
let us consider a sequence 
$\big(\mathcal{R}_r^{(n)}(T) = (\pi^{(n)},
\lambda^{(n)})\big)_{n\in\mathbb{N}}$ of iterates by the
renormalized Rauzy-Veech induction $\mathcal{R}_r$.
\begin{enumerate}
\item If $T$ has the Keane's property, then there exists $n_0$ such that
$\pi^{(n)}$ is irreducible for all $n \ge n_0$.

\item The renormalized Rauzy-Veech induction, defined on
the set $\{(\pi, \lambda)\,|\ \pi\ \text{irreducible}\}$, is
  recurrent.
\end{enumerate}
\end{ThmC}

Having a generalized permutation $\pi$ we can define one or two other
generalized permutations $\mathcal{R}_0(\pi)$ and $\mathcal{R}_1(\pi)$
reflecting the possibilities for the image of the Rauzy-Veech
induction $\mathcal{R}(T)$. These combinatorial Rauzy operations
define a partial order in the set of irreducible permutations
represented by an oriented graph. A \textit{Rauzy class} is a
connected component of this graph. 

Note that geometry of the Rauzy graphs is very different and more complicated
than in the case of ``true'' permutations since for some irreducible
generalized permutations one of the Rauzy operations might be not
defined. From Theorem~$C$ we will deduce that a Rauzy class is an
equivalence class for the equivalence relation given by these
combinatorial operations (see Proposition~\ref{prop:classes}).

In analogy with the case of the ``true'' permutations, we introduce one
more combinatorial operation on generalized permutations and define
\textit{extended Rauzy classes} as minimal subsets of irreducible
generalized permutations invariant under these corresponding three
operations. 

The moduli spaces of Abelian differentials and of quadratic
differentials are stratified by multiplicities of the zeroes of
the corresponding differentials. We denote a stratum of the moduli space
of strict quadratic differentials (with at most simple poles) 
by $\mathcal{Q}(k_1,\dots, k_n)$, where $k_i \ge -1$ are the
multiplicities of the zeroes ($k_i=-1$ corresponds to a pole).

\begin{ThmD}
Extended Rauzy classes of irreducible generalized permutations are in
one-to-one correspondence with connected components of strata in the
moduli spaces of quadratic differentials. 
\end{ThmD}

Historically, extended Rauzy classes where used to prove the
non-connectedness of some strata of Abelian differentials. For
permutations of a small number of elements, it is easy to construct
explicitly the subset of irreducible permutations and then using the
Rauzy operations to decompose it into a disjoint union of extended
Rauzy classes. Using this approach Veech proved that the minimal
stratum in genus $3$ has two connected components and Arnoux proved
that the minimal stratum in genus $4$ has three connected
components (for Abelian differentials). \medskip

Having established an explicit combinatorial criterion of
irreducibility of a generalized permutation (namely Theorem~$A$) one
can apply Theorem~$D$ to classify the connected components of all
strata of quadratic differentials of sufficiently small
dimension. This justifies, in particular, the following experimental
result of Zorich.

\begin{NNths}[Zorich]
Each of the following four exceptional strata of quadratic
differentials $\mathcal{Q}(-1, 9), \mathcal{Q}(-1, 3, 6),
\mathcal{Q}(-1, 3, 3, 3)$ and $\mathcal{Q}(12)$ contains exactly two
connected components.
\end{NNths}

Note that a theorem of the second author~\cite{La1} classifies all
connected components of all other strata of meromorphic quadratic
differentials with at most simple poles. These strata are either
connected, or contain exactly two connected components one of which being
hyperelliptic. The same theorem~\cite{La1} proves that each of the
remaining four exceptional strata might have at most two connected
components. However, the only currently available proof of the fact
they are disconnected is the one based on explicit calculation of the
extended Rauzy classes and corresponds to the theorem of Zorich. It
would be interesting to have an algebraic-geometrico proof of the last
theorem; namely a topological invariant as in the Kontsevich-Zorich's
classification~\cite{Kontsevich:Zorich}.
Note also that a paper of Zorich~\cite{Zo:representants} gives 
explicits representatives elements for each extended Rauzy class. 
See also~\cite{Zorich:experiments} for programs concerning 
calculations of these Rauzy classes.

\begin{Reader's guide}
In Section~\ref{sec:iem} we recall basic properties of flat surfaces,
moduli spaces and interval exchange maps. In particular we recall the
Rauzy-Veech induction and its dynamical properties. We relate these
properties to irreducibility. \\
In section~\ref{sec:giem} we recall the definition of a linear involution and give 
basics properties. Then in section~\ref{combinatoric}
we define a combinatorial notion of irreducibility, and prove the first part of 
Theorem~$A$. The main tool we use to prove this theorem is
the presentation proposed by Marmi, Moussa and Yoccoz which appears
in~\cite{Marmi:Moussa:Yoccoz} \\
In section~\ref{sec:dyn} we introduce the Keane's property for the 
linear involutions and prove the second  part of
Theorem~$A$, that is Theorem~$B$. For that we prove that $T$ satisfies
the Keane's property if and only if the Rauzy-Veech 
induction is always well defined and the length parameters tends to
zero. Then if $T$ does not satisfy the Keane's property we show that
there exists $n_0$ such that $\mathcal R^{n_0}(T)$ is dynamically
reducible which then implies that $T$ is also dynamically
reducible.\\
In section~\ref{dyn:rauzy:veech}, we study the dynamics of the renormalized Rauzy-Veech 
map on the space of the linear involutions, and prove
Theorem~$C$. For that we use the Teichm\"uller geometry
and the finiteness of the volume of the strata proved by Masur and
Veech (see~\cite{Ma82,Ve3}). \\
Section~\ref{rauzy:classes:giem} is devoted to a proof of
Theorem~$D$ on extended Rauzy classes; we present a result of Zorich
based on an explicit calculation of these classes in low genera.\\
In the Appendix we present some explicit Rauzy classes as illustration
of the problems which appear in the general case. We also give a
property concerning the extended Rauzy classes.

\end{Reader's guide}

\begin{Acknowledgments}
We thank Anton Zorich for useful discussions. We thank Arnaldo Nogueira for 
comments and remarks on a preliminary version of this text. \\
This work was partially supported by the ANR ``Teichm\"uller projet
blanc'' ANR-06-BLAN-0038.
\end{Acknowledgments}

\section{Background}
\label{sec:iem}

In this section we review basic notions concerning flat surfaces, moduli spaces 
and interval exchange maps. For general references see 
say~\cite{Ma82, Rauzy,Ve78,Ve82,Zorich} and~\cite{Masur:Tabachnikov}.
In this paper we will mostly follows notations presented in the
paper~\cite{Marmi:Moussa:Yoccoz}, or equivalently~\cite{Yoccoz}.

\subsection{Flat surfaces}

A {\it flat surface} is a (real, compact, connected) genus
$g$ surface equipped with a flat metric (with isolated conical
singularities) such that the holonomy group belongs to $\{\pm
\textrm{Id}\}$. Here holonomy means that the parallel transport of a
vector along a long loop  brings the vector back to itself or to its opposite.
This implies that all cone angles are integer multiples of $\pi$. Equivalently 
a flat surface is a triple $(S,\mathcal U,\Sigma)$ 
such that $S$ is a topological compact connected surface,
$\Sigma$ is a finite subset of $S$ (whose elements are called
{\em singularities}) and $\mathcal U = \{(U_i,z_i)\}$ is an
atlas of $S \setminus \Sigma$ such that the transition maps
$z_j \circ z_i^{-1} : z_i(U_i\cap U_j) \rightarrow
z_j(U_i\cap U_j)$ are translations or half-turns: $z_i = \pm
z_j + const$, and for each $s\in \Sigma$, there is a neighborhood of
$s$ isometric to a Euclidean cone.  Therefore we get a {\it quadratic
  differential} defined locally in the coordinates $z_i$ by
the formula $\Dq=\D  z_i^2$. This form extends to the points of $\Sigma$ to zeroes, 
simple poles or marked points (see~\cite{Masur:Tabachnikov}). We will
sometimes use the notation $(S, \Dq)$ or simply $S$. \smallskip

Observe that the holonomy is trivial if and only if there exists a sub-atlas such that all transition functions are
translations or equivalently if the quadratic differentials $\Dq$ is the
global square of an Abelian differential. We will then say that $S$ is a
translation surface.

\subsection{Moduli spaces}
\label{sec:moduli:spaces}

For $g \geq 1$, we define the moduli space of Abelian differentials
$\HA_g$ as the set of pairs $(S,\omega)$ 
modulo the equivalence relation generated by:
$(S,\omega) \sim (S',\omega')$ if there
exists an analytic isomorphism $f:S \rightarrow S'$ such that 
$f^\ast \omega'=\omega$. \\
For $g \geq 0$, we also define the moduli space of quadratic
differentials $\QQQ_g$ as the moduli space of pairs $(S,\Dq)$
(where $\Dq$ is not the global square of any Abelian differential) modulo 
the equivalence relation generated by:
$(S,\Dq) \sim (S',\Dq')$ if there
exists an analytic isomorphism $f:S \rightarrow S'$ such that 
$f^\ast \Dq'=\Dq$. \medskip

The moduli space of Abelian differentials (respectively quadratic differentials) is stratified by the 
multiplicities of the zeroes. We will denote by $\HHH(k_1, \dots, k_n)$ (respectively $\QQQ(k_1, \dots, k_n)$)
the stratum consisting of holomorphic one-forms (respectively quadratic differentials) with $n$ 
zeroes (or poles) of multiplicities $(k_1,\dots,k_n)$. These strata are 
non-connected in general (for a complete classification see~\cite{Kontsevich:Zorich}
in the Abelian differentials case and~\cite{La1} in the quadratic differentials case). \medskip

The linear action of the $1$-parameter subgroup of diagonal matrices
$g_t:=\textrm{diag}(e^{t/2},e^{-t/2})$ on the flat surfaces presents a particular
interest. It gives a measure-preserving flow with respect to a natural measure 
$\mu^{(1)}$, preserving each stratum of area one flat surfaces. This
flow is known as the {\it Teichm\"uller geodesic flow}. Masur and
Veech proved the following theorem. 

\begin{NNths}[Masur;~Veech]
The Teichm\"uller geodesic flow acts ergodically on each connected
component of each stratum of the moduli spaces of area one Abelian and quadratic differentials 
(with respect to a finite measure in the Lebesgue class).
\end{NNths}

This theorem was proved by Masur~\cite{Ma82} and 
Veech~\cite{Ve82} for the $\HHH(k_1, \dots, k_n)$ case and for the
$\QQQ(4g-4)$ case. \\ 
The ergodicity of the Teichm\"uller geodesic flow is proved in full 
generality in~\cite{Ve86}, Theorem~$0.2$. The finiteness of the measure appears
in two $1984$ preprints of Veech: Dynamical systems on analytic manifolds of quadratic 
differentials I,II (see also~\cite{Ve86} p.445). These preprints have been published in 
$1990$~\cite{Ve3}.

\subsection{Interval exchange maps}

In this section we recall briefly the theory of interval exchange
maps. We will show that, under simple combinatorial conditions, such
transformations arise naturally as Poincar\'e return maps of measured
foliations and geodesic flows on translation surfaces. Moreover we
will present the Rauzy-Veech induction and its geometric and
dynamical properties (see~\cite{Ve82} for more details). \medskip

Let $I \subset \mathbb R$ be an open interval and let us choose 
a finite subset $\{sing\}$ of $I$. Its complement is a union of $d\geq 2$ open  subintervals. An interval exchange map is a
one-to-one map $T$ from $I\backslash \{sing\}$ to a co-finite subset of $I$
that is a translation on each subinterval of its definition domain. It is easy to see
that $T$ is precisely determined by the following data: a permutation $\overline{\pi}$
that encodes how the intervals are exchanged (expressing that the k-th
interval, when numerated from the left to the right, is sent by $T$ to
the place $\overline{\pi}(k)$), and a vector with positive entries
that encodes the lengths of the intervals. \medskip

Following Marmi, Moussa, Yoccoz~\cite{Marmi:Moussa:Yoccoz}, we denote these intervals by 
$\{I_\alpha, \ \alpha\in \mathcal A\}$, with $\mathcal{A}$ a finite
alphabet. The length of the intervals is a vector
$\lambda=(\lambda_\alpha)_{\alpha\in \mathcal{A}}$, and the
combinatorial data is a pair 
$\pi=(\pi_0,\pi_1)$ of one-to-one maps $\pi_\epsilon: \mathcal A
\rightarrow \{1,\dots,d\}$. Then $\overline{\pi}$ is a one-to-one map from
$\{1,\dots,d\}$ into itself given by $\overline{\pi}=\pi_1 \circ \pi_0^{-1}$.
We will usually represent such a permutation by a table:
\begin{multline*}
\overline{\pi}=\left(\begin{array}{ccccc}1&2&\ldots&n \\\overline{\pi}^{-1}(1) 
&\overline{\pi}^{-1}(2)  &
  \ldots & \overline{\pi}^{-1}(n)& \end{array}\right) = \\ 
= \left(\begin{array}{ccccc}\pi_0^{-1}(1)&\pi_0^{-1}(2)&\ldots&\pi_0^{-1}(n)
  \\
\pi_1^{-1}(1)&\pi_1^{-1}(2)&\ldots&\pi_1^{-1}(n) \end{array}\right).\\
\end{multline*}

\begin{Example}
\label{example:permutation}
Let us consider the following alphabet $\mathcal{A}=\{A,B,C,D\}$ with
$d=4$. Then we define a permutation $\pi$ as follows.
$$
\overline{\pi}=\left(\begin{array}{cccc}
A & B & C & D \\
D & C & B & A
\end{array}\right). 
$$
\begin{figure}[htbp]
   \begin{center}
\psfrag{i}{$I$}  \psfrag{T}{$T$}  
\psfrag{a}{$I_A$}  \psfrag{b}{$I_B$}  
\psfrag{c}{$I_C$}  \psfrag{d}{$I_D$} 
\psfrag{ta}{$T(I_A)$}  \psfrag{tb}{$T(I_B)$}  
\psfrag{tc}{$T(I_C)$}  \psfrag{td}{$T(I_D)$} 

   \includegraphics{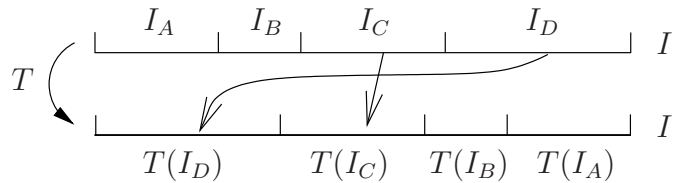}
     \caption{An interval exchange map.}
   \end{center}
\end{figure}
\end{Example}

\subsubsection{Rauzy-Veech induction}

In this section we introduce the notion of winner and loser, following
the terminology of the paper of Avila, Gou\"ezel and
Yoccoz~\cite{Avila:Gouezel:Yoccoz}. For $T=(\pi,\lambda)$ we define
the {\it type} $\varepsilon$ of $T$ by
$\lambda_{\pi_{\varepsilon}^{-1}(d)}>\lambda_{\pi_{1-\varepsilon}^{-1}(d)}$.
We will then say that $I_{\pi_{\varepsilon}^{-1}(d)}$ is the
winner and $I_{\pi_{1-\varepsilon}^{-1}(d)}$ is the loser. 
Then we define a subinterval $J$ of $I$ by removing the loser
of $I$ as follows.
$$
\left\{ \begin{array}{ll}
J=I \backslash T(I_{\pi_{1}^{-1}(d)}) & \textrm{if $T$ is of type 0}\\
J=I \backslash I_{\pi_0^{-1}(d)} & \textrm{if $T$ is of type 1.} 
\end{array} \right.
$$
The Rauzy-Veech induction $\mathcal R(T)$ of $T$ is defined as the first
return map of $T$ to the subinterval $J$. It is easy to see that this
is again an interval exchange transformation, defined on $d$
letters (see e.g.~\cite{Rauzy}). We now see how to compute the data of
the new map. \medskip

There are two cases to distinguish depending on whether $T$ is of type
$0$ or $1$; the combinatorial data of $\mathcal{R}(T)$ only
depends on $\pi$ and on the type of $T$. This defines two maps
$\mathcal{R}_0$ and $\mathcal{R}_1$ by
$\mathcal{R}(T)=(\mathcal{R}_\varepsilon(\pi), \lambda^{\prime})$,
with $\varepsilon$ the type of $T$.

\begin{enumerate}
\item $T$ has type $0$; equivalently the winner is $I_{\pi_0^{-1}(d)}$. \\
In that case, we define  $k$ by $\pi_1^{-1}(k)=\pi_0^{-1}(d)$ where
$k\leq d-1$. In an equivalent way $k=\pi_1 \circ
\pi_0^{-1}(d)=\overline{\pi}(d)$. 
Then $\mathcal R_0(\pi_0,\pi_1)=(\pi_0',\pi_1')$ where $\pi_0=\pi_0'$ and
$$
\pi_1'^{-1}(j) = \left\{
\begin{array}{ll}
\pi_1^{-1}(j) & \textrm{if $j\leq k$}\\
\pi_1^{-1}(d) & \textrm{if $j = k+1$}\\
\pi_1^{-1}(j-1) & \textrm{otherwise.}
\end{array} \right.
$$
We have $\lambda_\alpha'=\lambda_\alpha$ if $\alpha \not =
\pi_0^{-1}(d)$ and
$\lambda_{\pi_0^{-1}(d)}'=\lambda_{\pi_0^{-1}(d)}-\lambda_{\pi_1^{-1}(d)}$.

\item $T$ has type $1$; equivalently the winner is $I_{\pi_1^{-1}(d)}$. \\
In that case, we define $k$ by $\pi_0^{-1}(k)=\pi_1^{-1}(d)$ where $k \leq
d-1$. In an equivalent way $k=\pi_0 \circ \pi_1^{-1}(d)=\overline{\pi}^{-1}(d)$.
Then $\mathcal R_1(\pi_0,\pi_1)=(\pi_0',\pi_1')$ where $\pi_1=\pi_1'$ and
$$
\pi_0'^{-1}(j) = \left\{
\begin{array}{ll}
\pi_0^{-1}(j) & \textrm{if $j\leq k$}\\
\pi_0^{-1}(d) & \textrm{if $j = k+1$}\\
\pi_0^{-1}(j-1) & \textrm{otherwise.}
\end{array} \right.
$$
We have $\lambda_\alpha'=\lambda_\alpha$ if $\alpha \not =
\pi_1^{-1}(d)$ and
$\lambda_{\pi_1^{-1}(d)}'=\lambda_{\pi_1^{-1}(d)}-\lambda_{\pi_0^{-1}(d)}$.

\end{enumerate}

\begin{Example}
Let $\mathcal{A}=\{A,B,C,D\}$ be an alphabet. Let us consider the
permutation $\pi$ of Example~\ref{example:permutation}. Then
$$
\mathcal R_0\pi=\left(\begin{array}{cccc}A & B & C & D 
\\D & A & C & B \end{array}\right) \qquad \textrm{and} \qquad 
\mathcal R_1\pi=\left(\begin{array}{cccc}A & D & B & C  
\\D & C & B & A \end{array}\right).
$$
\end{Example}

We stress that the Rauzy-Veech induction is well defined if and only if the two
rightmost intervals do not have  the same length
i.e. $\lambda_{\pi_0^{-1}(d)} \not = \lambda_{\pi_1^{-1}(d)}$. In the
next, we want to study the Rauzy-Veech induction as a dynamical system defined on the space of interval exchange transformations. 
Thus we want the iterates of the Rauzy-Veech induction on $T$ to be
always well defined. We also want this induction to be a good
renormalization process, in the sense that the iterates correspond to
inductions on  subintervals that tend to zero. 
This leads to the definition of reducibility and to the
Keane's property.

\subsubsection{Rauzy-Veech induction and Keane's property} 

We will say that $\pi=(\pi_0,\pi_1)$ is reducible if there exists $1
\leq k \leq d-1$ such that $\{1,\dots,k\}$ is invariant under
$\overline{\pi}=\pi_1 \circ \pi_0^{-1}$. This means exactly that $T$ splits into
two interval exchange transformations.  \medskip

We will say that $T$ satisfies the Keane's property (also called the infinite
distinct orbit condition or i.d.o.c. property), if the orbits of the
singularities of $T^{-1}$ by $T$ are infinite. This ensures that $\pi$
is irreducible and the iterates of the Rauzy-Veech induction are
always well defined.

If the $\lambda_\alpha$ are rationally independent vectors, that is
$\sum r_\alpha \lambda_\alpha \not = 0$ for all nonzero integer vectors
$(r_\alpha)$, then $T$ satisfies the Keane's property (see \cite{Keane}). However the
converse is not true. Note that if $T$ satisfies the Keane's property
then $T$ is minimal. \medskip

Let $T=(\pi,\lambda)$ be an interval exchange map. Let us denote by
$\lambda_\alpha^{(n)}$ the length of the interval associated to the
symbol $\alpha\in \mathcal A$ for the $n$-th iterate of $T$ by
$\mathcal R$; we denote $\mathcal R^n(T)=:(\pi^{(n)},\lambda^{(n)})$
if it is well defined.

\begin{NNprop}
The following are equivalent.
\begin{enumerate}
\item $T$ satisfies the Keane's property.

\item The Rauzy-Veech induction $\mathcal R$ is always well-defined and 
 for any $\alpha\in \mathcal A$, the length of the intervals 
  $\lambda^{(n)}_\alpha$ goes to zero as $n$ tends to infinity. 
\end{enumerate}
\end{NNprop}

As we will see this situation is very similar in the case of
linear involutions. \medskip

If we want to study the Rauzy-Veech induction as a dynamical system on
the space of interval exchange maps, it is useful to consider the
Rauzy-Veech renormalisation on the projective space of 
lengths parameters space. The natural associated object is the 
renormalized Rauzy-Veech induction defined on length one intervals:
$$
\textrm{if } \mathcal R(\pi,\lambda) = (\pi',\lambda') \textrm{ then } 
\mathcal R_r(\pi,\lambda):= (\pi',\lambda'/|\lambda'|).
$$

\subsubsection{Rauzy classes}\label{rauzy:classes:iem}

Given a permutation $\pi$, we can define two other permutations
$\mathcal R_\varepsilon(\pi)$ with $\varepsilon=0,1$. Conversely, any
permutation $\pi'$ has exactly two predecessors: there exist exactly
two permutations~$\pi^0$ and~$\pi^1$ such that $\mathcal
R_\varepsilon(\pi^\varepsilon)=\pi'$. Note that $\pi$ is irreducible
if and only if $\mathcal R_\varepsilon(\pi)$ is irreducible. Thus the
relation generated by $\pi \sim \mathcal{R}_\varepsilon(\pi)$ is a
partial order on the set of irreducible permutations; we represent it
as a directed graph $G$. We call Rauzy classes the connected
components of this graph.

\begin{NNprop}[Rauzy]
The above relation is an equivalence relation on the set of
permutations. In particular, the equivalent class of a permutation is
the Rauzy class.
\end{NNprop}

\begin{proof}
The key remark is the following: if $\pi^{\prime}=\mathcal{R}_\varepsilon(\pi)$ 
then there exists $n>0$ such that $\pi=\mathcal{R}_\varepsilon^n(\pi^{\prime})$.
Now assume there exists an oriented path in $G$ joining  $\pi$
and $\pi^{\prime}$, \emph{i.e.} there exist
$\varepsilon_1,\dots,\varepsilon_r$ such that
$\pi^{\prime}=\mathcal{R}_{\varepsilon_1}\circ\dots\circ\mathcal{R}_{\varepsilon_r}(\pi)$.
Then there exists $n_1$ such that
$\mathcal{R}_{\varepsilon_1}^{n_1}(\pi')=\mathcal{R}_{\varepsilon_2}\circ\dots\circ\mathcal{R}_{\varepsilon_r}(\pi)$.
Iterating this argument, there exist $n_1,\dots,n_r$ such that
$\pi=\mathcal{R}_{\varepsilon_r}^{n_r}\circ\dots\circ\mathcal{R}_{\varepsilon_1}^{n_1}(\pi^{\prime})$.
Thus there is an oriented path in $G$ that joins $\pi^{\prime}$ and
$\pi$.
\end{proof}
We will see that there is an analogous proposition in the case of generalized permutations although the situation is much more complicated.

\subsubsection{Suspension data over an interval exchange
  transformation}
\label{susp:data:iem}

Here we describe the construction of a \emph{suspension} over an
interval exchange map $T$, that is a flat surface for which $T$ is the
first return map of the vertical flow on a well chosen segment.  

Let $T=(\pi,\lambda)$ be an interval exchange transformation.  A
\emph{suspension data} for $T$ is a collection of vectors
$(\zeta_\alpha)_{\alpha\in \mathcal{A}}$ such that: 
\begin{enumerate}
\item $\forall \alpha \in \mathcal A,\  
Re(\zeta_\alpha)=\lambda_\alpha$.

\item $\forall 1 \leq k \leq d-1,\ Im(\sum_{\pi_0(\alpha)\leq k} \zeta_\alpha)>0$. 

\item $\forall 1 \leq k \leq d-1,\ Im(\sum_{\pi_1(\alpha)\leq k} \zeta_\alpha)<0$. 

\end{enumerate}

\begin{figure}[htbp]
   \begin{center}
    \psfrag{1}[][]{\small $\zeta_A$}
    \psfrag{2}[][]{\small $\zeta_B$}
    \psfrag{3}[][]{\small $\zeta_C$}
    \psfrag{4}[][]{\small $\zeta_D$}
   \includegraphics{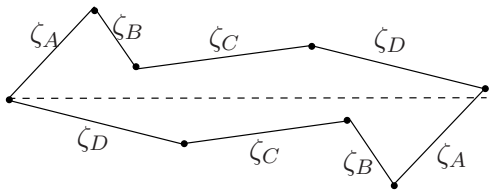}
     \caption{Suspension over an interval exchange transformation.}
     \label{fig:suspension:iet}
   \end{center}
\end{figure}

Given a suspension datum $\zeta$, we consider the broken line $L_0$ on
$\mathbb{C}=\mathbb R^2$ defined by concatenation of the vectors
$\zeta_{\pi_0^{-1}(j)}$ (in this order) for $j=1,\dots,d$ with
starting point at the origin (see
Figure~\ref{fig:suspension:iet}). Similarly, we consider the broken
line $L_1$  defined by concatenation of the vectors
$\zeta_{\pi_1^{-1}(j)}$ (in this 
order) for $j=1,\dots,d$ with starting point at the origin.
If the lines $L_0$ and $L_1$ have no
intersections other than the endpoints, then we can construct a
translation surface $S$ as follows: we can identify each side
$\zeta_\alpha$ on $L_0$ with the side $\zeta_\alpha$ on $L_1$
by a translation (in the general case, we must use the Veech zippered
rectangle construction, see section \ref{subsec:zip}). Let $I \subset
S$ be the horizontal interval defined by $I = (0,\sum_{\alpha \in \mathcal
    A}\lambda_\alpha) \times \{0\}$.
 Then the interval exchange map $T$ is precisely the one defined by the
first return map to $I$ of the vertical flow on $S$.

We have not yet discussed the existence of such a suspension datum for
a general interval exchange map. A necessary condition for $T$ to have
suspension data is that $\pi$ is irreducible. Indeed, if we have
$1\leq k \leq d-1$ such that $\pi_1\circ
\pi_0^{-1}(\{1,\dots,k\})=\{1,\dots,k\}$, and let
$\zeta=(\zeta_\alpha)_\alpha$ be a collection of complex numbers,
then:
$$
\sum_{\pi_0(\alpha)\leq k} \zeta_\alpha=\sum_{\pi_1(\alpha)\leq k}
\zeta_\alpha.
$$
So the imaginary part of this number cannot be both positive and
negative, and $\zeta$ is not a suspension data for $T$. If $\pi$ is
irreducible, the existence of a suspension data is given by Masur and
Veech independently (see~\cite{Ma82} page 174 and~\cite{Ve82}
formula $3.7$ page 207). We explain the construction here.

First, let us remark that $\pi=(\pi_0,\pi_1)$ is irreducible if
and only if
\begin{equation}
\label{eq:irr:masurveech1}
\sum_{i=1}^k \pi_1 \circ \pi_0^{-1}(i) - i > 0 \qquad \textrm{for any } 1 \leq
k \leq d-1.
\end{equation}
Of course if $\pi$ is irreducible, then so is $\pi^{-1}$, therefore 
\begin{equation}
\label{eq:irr:masurveech2}
\sum_{i=1}^k \pi_0 \circ \pi_1^{-1}(i) - i > 0 \qquad \textrm{for any } 1 \leq
k \leq d-1.
\end{equation}
Let us define a collection of complex number $\zeta=(\zeta_\alpha)_\alpha$ as follows:
$$
\zeta_\alpha = \lambda_\alpha + i (\pi_1(\alpha) - \pi_0(\alpha))
\qquad \textrm{for any } \alpha \in \mathcal A.
$$
Then following~\eqref{eq:irr:masurveech1} and
\eqref{eq:irr:masurveech2}, the collection $(\zeta_\alpha)_{\alpha \in
  \mathcal A}$ is a suspension data over $T$.

\subsubsection{Zippered rectangles}
\label{subsec:zip}
\begin{figure}[htbp]
    \psfrag{A}[][]{\small $A$}
    \psfrag{B}[][]{\small $B$}
    \psfrag{C}[][]{\small $C$}
    \psfrag{D}[][]{\small $D$}
   \begin{center}
\includegraphics{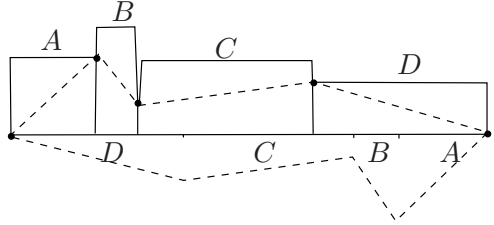}
    \caption{Zippered rectangles construction.}
     \label{fig:zip:iet}
   \end{center}
\end{figure}

Here we describe an alternative construction of the
suspension over an interval exchange transformation that works for
\emph{any} suspension data, namely the so 
called zippered rectangles construction due to Veech~\cite{Ve82}. Let
$T=(\pi,\lambda)$ be an interval exchange map, and let us assume that
$\pi$ is irreducible. Let $\zeta$ be any suspension over $T$. Then we
define $h=(h_\alpha)_{\alpha\in \mathcal A}$ by
$$
h_\alpha = \sum_{\pi_0(\beta)<\pi_0(\alpha)} Im(\zeta_\beta) -
\sum_{\pi_1(\beta)<\pi_1(\alpha)} Im(\zeta_\beta) > 0.
$$
For each $\alpha \in \mathcal A$ let us consider a rectangle
$R_\alpha$ of width $Re(\zeta_\alpha)$ and of height
$h_\alpha$ based on $I_{\pi_0(\alpha)}\subset I$.
The zippered rectangle construction is the translation surface
$\bigcup_{\alpha \in \mathcal A} R_\alpha / \sim$
where $\sim$ is the following equivalence relation: we identify the
top and the bottom of these rectangles by $(x,h_\alpha) \sim
(T(x),0)$ for $x\in I_{\pi_0(\alpha)}$. Then we ``zip'' the vertical
boundaries of these rectangles that are adjacent (see
figure~\ref{fig:zip:iet}; see also~\cite{Ve82} for a more precise description).

\subsubsection{Rauzy-Veech induction on suspensions}
We can define the Rauzy-Veech induction on the space of suspensions, as
well as on the space of zippered rectangles. Let $T=(\pi,\lambda)$ be
an interval exchange map and let $\zeta$ be a suspension over $T$. Then we
define $\mathcal R(\pi,\zeta)=(\pi',\zeta')$ as follows.

We define $(\pi',Re(\zeta'))=\mathcal R(\pi,Re(\zeta))$ (the
standard Rauzy-Veech induction). If $I_{\pi_\varepsilon^{-1}(d)}$ is
the winner for $T=(\pi,Re(\zeta))$ then
$$
\left\{ \begin{array}{l}
\zeta'_{\pi_\varepsilon^{-1}(d)}=\zeta_{\pi_\varepsilon^{-1}(d)}
- \zeta_{\pi_{1-\varepsilon}^{-1}(d)} \\
\zeta'_\alpha=\zeta_\alpha \qquad \textrm{otherwise.}
\end{array} \right.
$$

\begin{figure}[htbp]
   \begin{center}
    \psfrag{1}[][]{\small $\zeta_A$}
    \psfrag{2}[][]{\small $\zeta_B$}
    \psfrag{3}[][]{\small $\zeta_C$}
    \psfrag{4}[][]{\small $\zeta_D$}

    \psfrag{1p}[][]{\small $\zeta'_A$}
    \psfrag{2p}[][]{\small $\zeta'_B$}
    \psfrag{3p}[][]{\small $\zeta'_C$}
    \psfrag{4p}[][]{\small $\zeta'_D$}

   \includegraphics{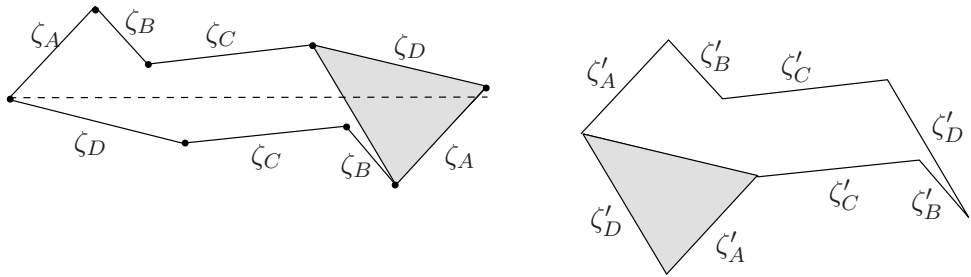}
     \caption{Rauzy-Veech induction on a suspension over an interval exchange transformation. 
The corresponding map is of type $0$ hence the new suspension data are $\zeta'_A=\zeta_A$, 
$\zeta'_B=\zeta_B$, $\zeta'_C=\zeta_C$ and $\zeta'_D=\zeta_D-\zeta_A$.
}
     \label{fig:rauzy:suspension}
   \end{center}
\end{figure}

\begin{rem}
\label{rk:moduli}
Since $(\pi^{\prime},\zeta^{\prime})$ is obtained from $(\pi,\zeta)$
by ``cutting'' and ``gluing'', these two surfaces differ by an element
of the mapping class group, hence they define the same point in the moduli space (see Figure~\ref{fig:rauzy:suspension} for an example).
\end{rem}

If $C$ is a Rauzy class, we define
$$
\mathcal{T}_C = \{(\pi,\zeta),\ \pi\in C,\zeta \textrm{ is a
  suspension data for  $\pi$} \}.
$$ 
We have thus defined the Rauzy-Veech map on the space
$\mathcal{T}_C$. It is easy to check that it defines an almost
everywhere invertible map: If $\sum Im(\zeta'_\alpha) \not = 0$ then
every $(\pi',\zeta')$ has exactly one preimage for $\mathcal R$.

\subsubsection{Moduli spaces and Rauzy-Veech induction}

We define the quotient $\mathcal{H}_C=\mathcal{T}_C / \sim$ of $\mathcal{T}_C$ by 
the equivalence relation generated by $(\pi,\zeta)\sim \mathcal{R}(\pi,\zeta)$. 
The zippered rectangle construction,  provides a mapping $p$ from 
$\mathcal{H}_C$ to a  stratum $\HHH(k_1,\dots,k_n)$ of the moduli space of Abelian 
differentials (see Remark~\ref{rk:moduli}). Observe that
$(k_1,\dots,k_n)$ can be calculated in terms of $C \ni \pi$. One can
also show that $\mathcal{H}_C$ is connected and so the image belongs
to a connected component of a stratum.

We will denote by $m$ the natural Lebesgue measure  on $\mathcal{T}_C$
i.e. $m=d\pi d\zeta$, were $d\pi$ is the counting measure on $C$ and
$d\zeta$ is the Lebesgue measure. The mapping $\mathcal{R}$ preserves
$m$, so it induces a measure, denoted again by  $m$ on
$\mathcal{H}_C$.

There is  natural action of the matrix
\begin{eqnarray*}
g_t=
\left(
\begin{array}{cc}
 e^{\frac{t}{2}} &  0    \\
 0  &    e^{-\frac{t}{2}}
\end{array}
\right)
\end{eqnarray*}
on $\mathcal{T}_C$ by
$g_t(\pi,\zeta)=(\pi,(g_t(\zeta_\alpha))_\alpha)$, where $g_t$ acts on
$\zeta_\alpha\in \mathbb{C}=\mathbb{R}^2$ linearly. This action
preserves the measure $m$ on $\mathcal{T}_C$ and commutes
with~$\mathcal{R}$, so it descends to a 1-parameter action on
$\mathcal{H}_C $ called the {\it Teichm\"uller flow}. 
Since the action of $g_t$ on $\mathcal{H}_C$ preserves the area of the
corresponding  flat surface, the Teichm\"uller flow also acts on the
subset $\mathcal{H}_C^1$ corresponding to area one surfaces, and
preserves the measure $m^{(1)}$ induced by the measure $m$ on that
subset. Note also that
$$
\left\{(\pi,\zeta)\in \mathcal{T}_C;\ 1 \leq |Re(\zeta)| \leq
1+\min\bigl(Re(\zeta_{\pi_{0}^{-1}(d)}),Re(\zeta_{\pi_1^{-1}(d)}\bigr)
\right\} 
$$
is a fundamental domain of $\mathcal{T}_C$ for the relation $\sim$
and the Poincar\'e map of the Teichm\"uller flow on  
$$
\mathcal{S}=\{(\pi,\zeta);\ \pi \textrm{ irreducible},\
|Re(\zeta)|=1\}/ \sim 
$$ 
is precisely the renormalized Rauzy-Veech induction on suspensions.

One can show (see \cite{Ve82}) that the mapping $p$ is a finite
covering from $\mathcal H_C^{1}$ onto a subset of full measure in a
connected component of a stratum and the measure $m$ 
projects to the measure $\mu^{(1)}$ defined in
section~\ref{sec:moduli:spaces}. Moreover the action of $g_t$ is
equivariant with respect to $p$, that is $p\circ g_t(\pi,\zeta) = g_t
\circ p(\pi,\zeta)$. Hence if we restrict to area one surfaces, the
result of Masur and Veech (finiteness of the measure) implies that the
measure $m^{(1)}$ is finite on $\mathcal{H}_C^1$.

\begin{cor}
The renormalized Rauzy-Veech induction is recurrent on $\mathcal S$.
\end{cor}

\begin{rem}
Veech proved a stronger result, that is the ergodicity of $g_t$ (on
the level of $\mathcal H_C$ for any Rauzy class $C$), which
implies the ergodicity of the Teichm\"uller flow for Abelian
differentials (see~\cite{Ve82}). He also proved that the induced
measure on $\mathcal{S}$ is always infinite. 
\end{rem}

\section{Linear involutions}
\label{sec:giem}

\subsection{linear involutions and generalized  permutations}
\label{section:def:giem}

Let $S$ be a (compact, connected, oriented) flat surface with
$\mathbb{Z}/2\mathbb{Z}$ linear holonomy and let 
$X$ be a horizontal segment with a choice of a positive vertical
direction (or equivalently, a choice of left and right ends). We
consider the first return map $T_0:X\rightarrow X$ of vertical geodesics starting from $X$ 
 in the positive direction. Any vertical
geodesic which start from $X$ and doesn't hit a singularity  will intersect $X$ again. Therefore, the
map $T_0$ is well defined outside a finite number of points $\{sing\}$ (called
singular points) that correspond to vertical geodesics that stop at a
singularity before intersecting again the interval $X$. 
The set $X\backslash \{sing\}$ is a finite union of open intervals $(X_i)$
and the restriction of $T_0$ on 
each of these intervals is of the kind $x\mapsto \pm x+c_i$. 

The map
$T_0$ alone does not properly correspond to the dynamics of  vertical geodesics
since when $T_0(x)=-x+c_i$ on the interval $X_i$, then $T_0^2(x)=x$,
and $(x,T_0(x), T_0^2(x))$ does not correspond to the successive
intersections of a vertical geodesic with $X$ starting from $x$. To fix
this problem, we have to consider $T_1$ the first return map of the
vertical geodesics starting from $X$ 
 in the negative direction. Now if $T_0(x)=-x+c_i$ then the
successive intersections with $X$ of the vertical geodesic starting
from $x$ will be $x,T_0(x), T_1(T_0(x))$\ldots

We get a dynamical system on $X\times \{0,1\}$. Following Danthony and
Nogueira (see~\cite{No,DaNo:CRAS,DaNo}) we will call such a dynamical
system a linear involution. We recall here the definition that we have 
restricted to our purpose.

\begin{defs}\label{def:giem}
Let $X$ be an open interval and let $\widehat X=X\times \{0,1\}$ be
two disjoint copies of $X$. A linear involution on $X$
is a  map $T:=f\circ \tilde{T}$, where: 
\begin{itemize}
\item $\tilde{T}$ is a smooth involution without fixed point defined
  on  $\widehat{X}\backslash \{sing\}$, where $\{sing\}$ is a finite
  subset of $\widehat{X}$. 
\item If $p=(x,\varepsilon)$ and $T(p)$ belong to the same connected component of
  $\widehat X$ then the derivative of $\tilde{T}$ at  
$p$ is $-1$ otherwise the derivative of $\tilde{T}$ at $p$ is $1$.
\item $f$ is the involution $(x,\varepsilon)\mapsto (x,
  1-\varepsilon)$.
\end{itemize}
\end{defs}

\begin{conv}
\label{convention:oriented}
In this paper, we are interested with non oriented measured foliations
defined on {\it oriented} surfaces. Observe that the orientability of
the surface $S$ forces the second condition on the derivative of $T$
in Definition~\ref{def:giem}.
\end{conv}

\begin{figure}[htb]
\begin{center}
\psfrag{11}{$1$}
\psfrag{12}{$2$}
\psfrag{13}{$3$}
\psfrag{14}{$4$}
\psfrag{15}{$5$}
\psfrag{x}{\small $(x,\varepsilon)$}
\psfrag{X}{\small $X$}
\psfrag{ttx}{\small $\tilde{T}(x,\varepsilon)$}
\psfrag{txi}{\small $T(x,\varepsilon)$?}
\psfrag{tx}{\small $T(x,\varepsilon)$}
\psfrag{t2x}{\tiny $T^2(x,\varepsilon)$}
\includegraphics[scale=0.52]{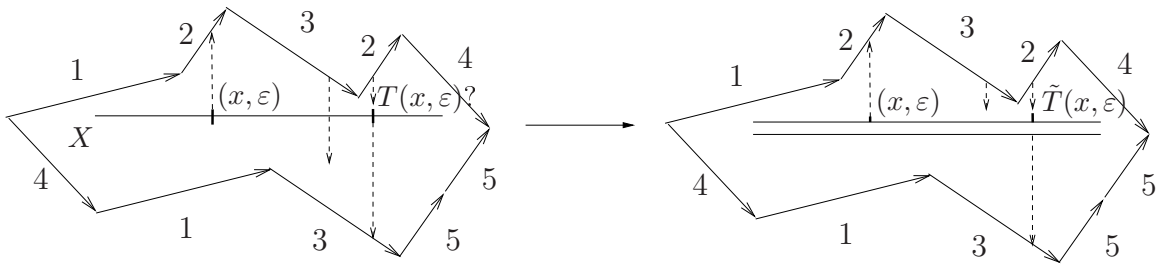}
\caption{Linear involution defined by the vertical foliation of a flat surface.}
\label{ttilde}
\end{center}
\end{figure}

The previous definition is motivated by the following remark.
 
\begin{rem}
The first return map of the vertical geodesic foliation on a horizontal
segment $X$ in a flat surface $S$ defines a linear involution
in the following way. 
Choose a positive vertical direction in a neighborhood of $X$, and
replace $X$ be two copies of $X$ as in Figure~\ref{ttilde}. We denote
by $X\times\{0\}$ the one on the top and by $X\times \{1\}$ the one on
the bottom. Then we consider the first return map on $X\times\{0,1\}$
of vertical geodesics, where a geodesic starting from $X\times\{0\}$
is taken in the positive vertical direction, and a geodesic starting
from $X\times\{1\}$ is taken in the negative direction. We obtain a
map $\tilde{T}$ and it is easy to check that $\tilde{T}$ satisfies the
condition of Definition~\ref{def:giem}. Then it is clear that the map
$T=f\circ \tilde{T}$ encodes the successive intersections of a vertical
geodesic with $X$.
\end{rem}

\begin{figure}[htbp]

\begin{center}

\psfrag{a}{$\scriptstyle X_A$} \psfrag{c}{$\scriptstyle X_C$}
\psfrag{b}{$\scriptstyle X_B$} \psfrag{d}{$\scriptstyle X_D$} 
\psfrag{e}{$\scriptstyle X_E$}
\psfrag{ta}{$\scriptstyle T(X_A)$} \psfrag{tc}{$\scriptstyle T(X_C)$}
\psfrag{tb}{$\scriptstyle T(X_B)$} \psfrag{td}{$\scriptstyle T(X_D)$}
\psfrag{te}{$\scriptstyle T(X_E)$}
\psfrag{x0}{$\scriptstyle X\times\{0\}$}\psfrag{x1}{$\scriptstyle X\times\{1\}$}
\psfrag{T}{$T$}
\psfrag{x}{$\scriptstyle (x,0)$}\psfrag{y}{$\scriptstyle (y,1)$}
\psfrag{tx}{$\scriptstyle T(x,0)$}\psfrag{ty}{$\scriptstyle T(y,1)$}
\psfrag{11}{$1$}
\psfrag{12}{$2$}
\psfrag{13}{$3$}
\psfrag{14}{$4$}
\psfrag{15}{$5$}
\psfrag{aa}{$\scriptstyle A$}  \psfrag{bb}{$\scriptstyle B$}
\psfrag{cc}{$\scriptstyle C$}  \psfrag{dd}{$\scriptstyle D$}
\psfrag{ee}{$\scriptstyle E$} 

\subfigure[]{\epsfig{figure=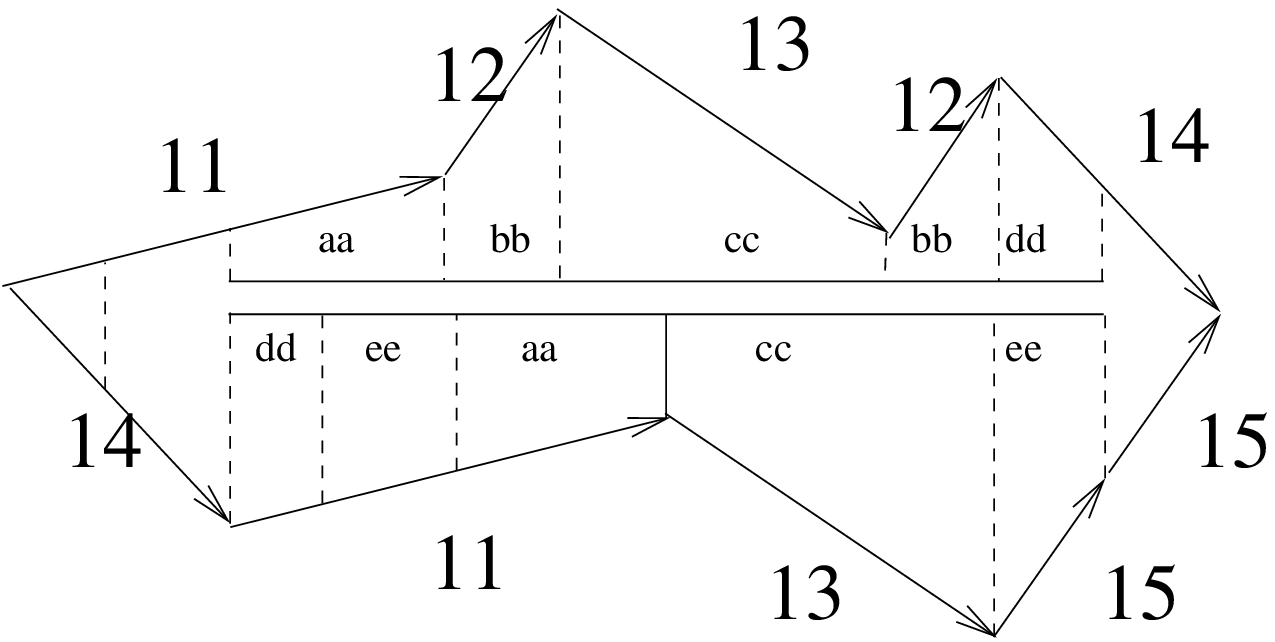,width=8cm}} \vskip 3mm
\subfigure[]{\epsfig{figure=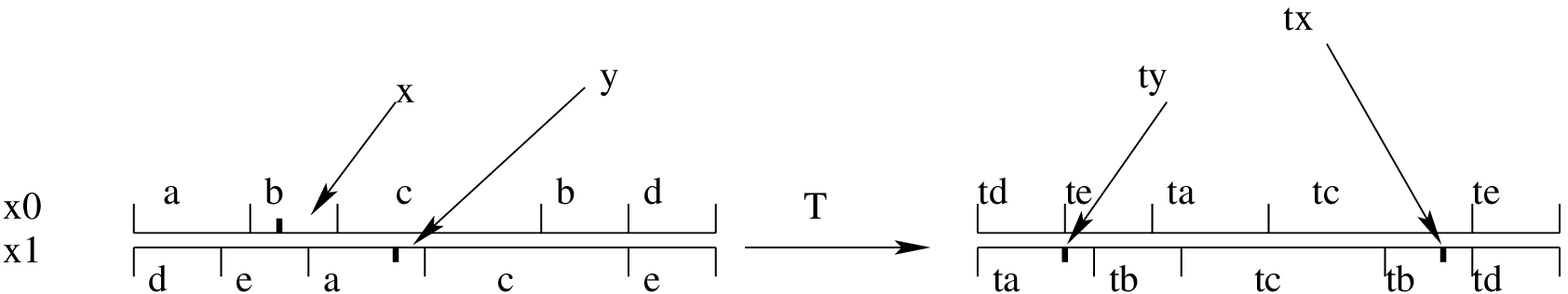,width=16.5cm}}

    \caption{A linear involution associated to a measured foliation on a 
flat surface.}
     \label{ex:giem}
   \end{center}
\end{figure}

Recall that interval exchange maps are encoded by combinatorial and
metric data: these are a permutation and a vector with positive entries.
We define an analogous object for linear involutions.

\begin{defs}\label{def2}
Let $\mathcal A$ be an alphabet of $d$ letters. A \emph{generalized
permutation} of type $(l,m)$, with $l+m=2d$, is a two-to-one map
$\gp:\{1,\ldots,2d\}\rightarrow \mathcal{A}$. We will usually
represent such generalized permutation by the table:
$$
\left(\begin{array}{ccc}
\pi(1) & \dots& \pi(l) \\
\pi(l+1) & \dots & \pi(l+m)
\end{array}\right). 
$$
A generalized permutation $\gp$ defines an involution $\sigma$
without fixed points by the following way
$$
\gp^{-1}(\{\gp(i)\}) = \{i,\sigma(i) \}.
$$
\end{defs}
     
Note that a permutation defines in a natural way a generalized permutation.

We now describe how a linear involution naturally defines
a generalized permutation. Let $T$ be a linear involution 
and let $\tilde{T}$ be the corresponding involution as in
Definition~\ref{def:giem}. The domain of definition of $\tilde{T}$ is
a finite union $X_1,\ldots,X_{l+m}$ of open intervals, where
$X_1,\ldots,X_l$ are subintervals of $X\times\{0\}$ and
$X_{l+1},\dots,X_{l+m}$ are subintervals of $X\times\{1\}$. Since
$\tilde{T}$ is an isometric involution without fixed point, each $X_i$
is mapped isometrically to a $X_j$, with $j\neq i$ , hence $\tilde{T}$
induces an involution without fixed point $\sigma_T$ on
$\{1,\dots,l+m\}$. As in section \ref{sec:iem}, we choose a name
$\alpha_i\in \mathcal{A}$ to each pair $\{i,\sigma_T(i)\}$ and we get
a generalized permutation in the sense of the above definition which is
defined up to a one-to-one map of $\mathcal{A}$.

\begin{Example}
\label{Example1}
In view of Figure~\ref{ex:giem}, 
let us consider the following alphabet $\mathcal{A}=\{A,B,C,D\}$ with
$d=5$. Then we define a generalized permutation $\pi$ as
follows.
\begin{eqnarray*}
l=m=5,\quad \pi(1)=\pi(8)=A,\ \pi(2)=\pi(4)=B,\\
\pi(3)=\pi(9)=C,\ \pi(5)=\pi(6)=D,\ \pi(7)=\pi(10)=E.
\end{eqnarray*}
In an equivalent way, we can define an involution without fixed point
in order to define $\pi$.
$$
\sigma(1)=8,\ \sigma(2)=4,\ \sigma(3)=9,\ \sigma(5)=6,\ \sigma(7)=10.
$$
We  represent $\gp$ by the following table
$$
\pi=\left(\begin{array}{ccccc}
A & B & C & B & D \\
D & E & A & C & E
\end{array}\right). 
$$
One can check that the discrete datum associated to the linear
involution described in Figure~\ref{ex:giem} is the generalized
permutation $\pi$.
     
\end{Example}

\begin{Example}
Note that $\pi$ is a ``true'' permutation on $d$ letters if and only
if $l=m=d$ and for any $i \leq l$, $\sigma(i) > l$. In this case (if
$\mathcal A=\{1,\dots,d\}$):
$$
\pi=\left(\begin{array}{ccccc}1&2&\ldots&d \\ \sigma(1)-d & \sigma(2)-d  &
  \ldots & \sigma(d)-d & \end{array}\right).
$$
\end{Example}

Conversely, let $\gp$ be a generalized permutation of type $(l,m)$
and let $\sigma$ be the associated involution. If $\gp$ is not
a ``true'' permutation, then an obvious necessary and sufficient condition for
$\gp$ to come from a linear involution is  
that there exist at least two indices $i\leq l$ and $j>l$ such that
$\sigma(i)\leq l$ and $\sigma(j)>l$.

\begin{conv}
\label{conv:normal}
From now, unless explicitly stated (in particular in
section~\ref{suff:cond}), we will always assume that generalized
permutations will satisfy the following convention. There exist at
least two indices $i\leq l$ and $j>l$ such that $\sigma(i)\leq l$ and
$\sigma(j)>l$.
\end{conv}

Let $(\lambda_{\alpha})_{\alpha \in
  \mathcal{A}}$ be a collection of positive real numbers such that
\begin{eqnarray}\label{condition:lambda:giem}
L:=\sum_{i=1}^l \lambda_{\gp(i)}=\sum_{i=l+1}^{l+m} \lambda_{\gp(i)}.
\end{eqnarray} 
It is easy to construct a linear involution on
the interval $X=(0,L)$ with combinatorial data $(\gp,\lambda)$. As 
in section \ref{sec:iem}, we will denote by $T=(\gp,\lambda)$ a
linear involution.

\subsection{Rauzy-Veech induction on linear involutions}

We recall the \emph{Rauzy-Veech induction} on linear 
involutions introduced by Danthony and Nogueira (see~\cite{DaNo} p.~473).

Let $T=(\gp,\lambda)$ be a linear involution on 
$X=(0,L)$, with $\gp$ of type $(l,m)$. If $\lambda_{\gp(l)}\neq
\lambda_{\gp(l+m)}$, then the Rauzy-Veech induction 
$\mathcal R(T)$ of $T$ is the linear involution
obtained by the first return map of $T$ to
$$
\bigl(0,\max(L-\lambda_{\gp(l)}, L-\lambda_{\gp(l+m)})\bigr) \times
\{0,1\}.
$$
As in the case of interval exchange maps, the combinatorial data of
the new linear involution depends only on the
combinatorial data of $T$ and whether
$\lambda_{\gp(l)}>\lambda_{\gp(l+m)}$ or $\lambda_{\gp(l)}<
\lambda_{\gp(l+m)}$. As before, we say that $T$ has type $0$ or type
$1$ respectively.
The corresponding combinatorial operations are
denoted by $\mathcal R_\varepsilon$ for $\varepsilon=0,1$
respectively. 
Note that if $\gp$ is a given generalized permutation, the subsets 
$\{T=(\gp,\lambda) ,\
\lambda_{\gp(l)}>\lambda_{\gp(l+m)}\}$ and  $\{T=(\gp,\lambda),
\lambda_{\gp(l)}< \lambda_{\gp(l+m)}\}$ can be empty because
$\gp(l)=\gp(l+m)$ or because of the linear relation on the ~$\lambda_i$ that must be satisfied. 

We first describe the combinatorial Rauzy operations $\mathcal
R_\varepsilon$. Let $\sigma$ be the associated involution to $\gp$.

\begin{enumerate}
\item map $\mathcal R_0$. \\
$\bullet$ If $\sigma(l)>l$ and if $\gp(l)\neq \gp(l+m)$ then we define $\mathcal
R_0\gp$ to be of type $(l,m)$ and such that: 
\[ \mathcal R_0\gp(i)= \left\{
\begin{array}{ll}
\gp(i) & \textrm{if $i\leq \sigma(l)$}\\
\gp(l+m) & \textrm{if $i=\sigma(l)+1$}\\
\gp(i-1) & \textrm{otherwise.}
\end{array}
\right. \]
$\bullet$ If $\sigma(l)\leq l$, and if there exists a pair
$\{x,\sigma(x)\}$ included in $\{l+1,\ldots,l+m-1\}$ then we define $\mathcal R_0\gp$ to
be of type $(l+1,m-1)$ and such that: 
$$
\mathcal R_0\gp(i)= \left\{
\begin{array}{ll}
\gp(i) & \textrm{if $i < \sigma(l)$}\\
\gp(l+m) & \textrm{if $i=\sigma(l)$}\\
\gp(i-1) & \textrm{otherwise.}
\end{array}
\right.
$$
$\bullet$ Otherwise $R_0 \gp$ is not defined.

\item map $\mathcal R_1$. \\
$\bullet$ If $\sigma(l+m)\leq l$ and if $\gp(l)\neq \gp(l+m)$ then we define
$\mathcal R_1\gp$ to be of type $(l,m)$ such that: 
\[ \mathcal R_1\gp(i)= \left\{
\begin{array}{ll}
\gp(l) & \textrm{if $i=\sigma(l+m)+1$}\\
\gp(i-1) & \textrm{If $\sigma(l+m)+1<i\leq l$}\\
\gp(i) & \textrm{otherwise.}
\end{array}
\right. \]
$\bullet$ If $\sigma(l+m)>l$ and if there exists a pair $\{x,\sigma(x)\}$ included in
$\{1,\ldots,l-1\}$ then $\mathcal R_1\gp$ is of type $(l-1,m+1)$ and: 
\[ \mathcal R_1\gp(i)= \left\{
\begin{array}{ll}
\gp(i+1) & \textrm{if $l\leq i <\sigma(l+m)-1$}\\
\gp(l) & \textrm{If $i=\sigma(l+m)-1$}\\
\gp(i) & \textrm{otherwise.}
\end{array}
\right. \]
$\bullet$ Otherwise $R_1 \gp$ is not defined.
\end{enumerate}

We now describe the Rauzy-Veech induction $\mathcal{R}(T)$ of $T$:
\begin{itemize}
\item 
If $T=(\gp,\lambda)$ has type $0$, then
$\mathcal{R}(T)=(\mathcal{R}_0\pi,\lambda^{\prime})$, with
$\lambda^{\prime}_\alpha=\lambda_\alpha$ if $\alpha\neq \pi(l)$ and
$\lambda_{\pi(l)}^{\prime}=\lambda_{\pi(l)}-\lambda_{\pi(l+m)}$.
\item
If  $T=(\gp,\lambda)$ has type $1$, then
$\mathcal{R}(T)=(\mathcal{R}_1\pi,\lambda^{\prime})$, with
$\lambda^{\prime}_\alpha=\lambda_\alpha$ if $\alpha\neq \pi(l+m)$ and
$\lambda_{\pi(l+m)}^{\prime}=\lambda_{\pi(l+m)}-\lambda_{\pi(l)}$.
\end{itemize}

\begin{Example}
Let us consider the permutation of Example~\ref{Example1}, namely 
$\gp=\left(\begin{smallmatrix} A & B & C & B & D \\D & E & A & C & E
\end{smallmatrix}\right)$. Then
$$
\mathcal R_0(\gp)=\left(\begin{array}{ccccc} A & B & C & B & D \\ D &
  E & E & A & C \end{array}\right) \textrm{ and }  
\mathcal R_1(\gp)=\left(\begin{array}{cccccc} A & B & C & B \\D & D & E & A & C & E
 \end{array}\right).
$$
\end{Example}

\begin{Example}\label{ex2}
Let us consider the permutation $\gp$ defined on the alphabet
$\mathcal A=\{A,B,C,D\}$ by $\gp=\left(\begin{smallmatrix} A & B & A &&\\ B & D & C & C & D
\end{smallmatrix}\right)$. Then 
$$
\mathcal R_0(\gp)=\left(\begin{array}{cccc} D & A & B & A \\ B & D & C & C
\end{array}\right)
$$
and $\mathcal R_1(\gp)$ is not defined. Indeed, consider any
linear involution with $\gp$ as combinatorial
data. Then we must have 
$$
2\lambda_A+\lambda_B= \lambda_B+2\lambda_C+2\lambda_D.
$$
Therefore we necessarily have $\lambda_D<\lambda_A$ and
$\lambda_D>\lambda_A$ never happens.
\end{Example}

\begin{Example}\label{ex3}
Consider the permutation $\gp$ defined on the alphabet
$\mathcal A=\{A,B,C\}$ by $\gp=\left(\begin{smallmatrix} A & B & A \\
  B & C & C \end{smallmatrix}\right)$. Then $\mathcal
R_\varepsilon(\gp)$ is not defined for any $\varepsilon$. 
Indeed, consider any
linear involution with $\gp$ as combinatorial
data. Then we must have  $\lambda_A=\lambda_C$, hence the Rauzy-Veech
induction of $T$ is not defined for any parameters.
\end{Example}

In the case of interval exchange maps, one usually define the Rauzy-Veech induction 
 only for irreducible combinatorial data. Here we have not yet
defined irreducibility. However, it will appear in
section~\ref{combinatoric} that some interesting phenomena with respect
to Rauzy-Veech induction appear also in the reducible case.

In the next section we will define a notion of irreducibility
which is equivalent to have a suspension data. It is easy to see that
a generalized permutation $\gp$ such that neither $\mathcal R_0(\gp)$ nor
$\mathcal R_1(\gp)$ is defined is necessarily reducible. However, the
permutation $\gp$ of Example~\ref{ex2} is irreducible (see
Definition~\ref{def:irred} and Theorem~\ref{CNS}) while $\mathcal R_1(\gp)$ is
not defined.

\subsection{Suspension data and zippered rectangles construction}
\label{susp:giem}
Starting from a linear involution $T$, we want to
construct a flat surface and a horizontal segment whose corresponding
first return maps $(T_0,T_1)$ of the vertical foliation give $T$. Such
surface will be called a \emph{suspension} over $T$, and the
parameters encoding this construction will be called \emph{suspension
data}.

\begin{defs}
Let $T$ be a linear involution and let
$(\lambda_{\alpha})_{\alpha \in \mathcal{A}}$ be the lengths of the
corresponding intervals. Let $\{\zeta_{\alpha}\}_{\alpha \in
\mathcal{A}}$ be a collection of complex numbers such that: 
\begin{enumerate} 
\item $\forall \alpha \in \mathcal{A} \quad Re(\zeta_{\alpha})=\lambda_{\alpha}$.  
\item $\forall 1\leq i \leq l-1 \quad Im(\sum_{j\leq i} \zeta_{\gp(j)})>0$ 
\item $\forall 1\leq i \leq m-1 \quad Im(\sum_{1\leq j\leq i} \zeta_{\gp(l+j)})<0$
\item $\sum_{1\leq i\leq l} \zeta_{\gp(i)} = \sum_{1\leq j\leq m}\zeta_{\gp(l+j)}$.
\end{enumerate} 
The collection $\zeta=\{\zeta_{\alpha}\}_{\alpha\in \mathcal{A}}$ is
called a \emph{suspension data} over $T$.
\end{defs}

We will also speak in an obvious manner of a suspension data for a
generalized permutation.

\begin{figure}[htbp]
\psfrag{a}{$\scriptstyle \zeta_A$}  \psfrag{b}{$\scriptstyle \zeta_B$}
\psfrag{c}{$\scriptstyle \zeta_C$}  \psfrag{d}{$\scriptstyle \zeta_D$}
\psfrag{e}{$\scriptstyle \zeta_E$}

   \begin{center}
     \includegraphics{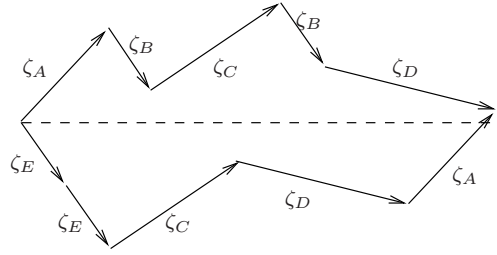}
    \caption{A suspension over a linear involution.}
     \label{figure:suspension:data}
   \end{center}
\end{figure}

Let $L_0$ be  a broken line (with a finite number of edges) on the
plane such that the edge number $i$ is 
represented by the complex number $\zeta_{\gp(i)}$, for $1\leq i\leq
l$, and $L_1$ be a broken line that starts on the same point as $L_0$,
and whose edge number~$j$ is represented by the complex number
$\zeta_{\gp(l+j)}$ for $1\leq j\leq m$ (Figure~\ref{figure:suspension:data}).

If $L_0$ and $L_1$ only intersect on their endpoints, then $L_0$ and
$L_1$ define a polygon whose sides comes by pairs and for each pair
the corresponding sides are parallel and have the same length.  Then
identifying these sides together, one gets a flat surface. It is
easy to check that the first return map of the vertical foliation on
the segment corresponding 
to $X$ in $S$ defines the same linear involution as $T$,
so we have constructed a suspension over $T$. We will say in this case
that $\zeta$ defines a \emph{suitable polygon}.

\begin{figure}[htbp]
\psfrag{a}{$\scriptstyle \zeta_A$}  \psfrag{b}{$\scriptstyle \zeta_B$}
\psfrag{c}{$\scriptstyle \zeta_C$}  \psfrag{d}{$\scriptstyle \zeta_D$}
\psfrag{e}{$\scriptstyle \zeta_E$}

   \begin{center}
     \includegraphics{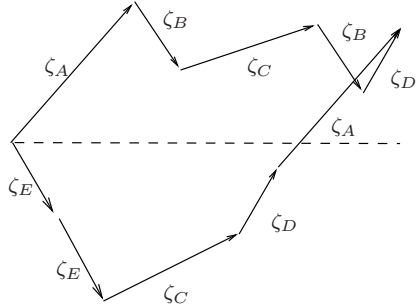}
    \caption{Suspension data that does not give a suitable polygon.}
     \label{figure:wrong:polygon}
   \end{center}
\end{figure}

The broken lines $L_0$ and $L_1$ might intersect at other points (see
Figure~\ref{figure:wrong:polygon}).
However, we can still define a flat surface by using an analogous
construction as the zippered rectangles construction.
We now give a sketch of this construction (see e.g.~\cite{Ve82,Yoccoz}
for the case of interval exchange maps, or
section~\ref{subsec:zip}). This construction is very similar 
to the usual one, although its precise description is very
technical. Still, for completeness, we give an equivalent but rather
implicit formulation. \medskip

We first consider the previous case when $L_0$ and $L_1$ define a
suitable polygon. For each pair of interval $X_i, X_{\sigma(i)}$ on
$X$, the return time $h_{\gp(i)}$ of the vertical foliation starting from
$x\in X_i$ and returning in $y\in X_{\sigma(i)}$ is constant. This
value depends only on the generalized permutation and on the
imaginary part of the suspension data $\zeta$.  
There is a natural embedding of the open rectangle
$R_{\gp(i)}=(0,\lambda_i)\times(0,h_{\gp(i)})$ into the flat surface
$S$ and this surface is obtained from $\sqcup_\alpha R_\alpha$ by
identifications on the boundaries of the $R_{\alpha}$. Identifications
for the horizontal sides 
$[0,\lambda_{\alpha}]$ are given by the linear involution 
and identifications for the vertical sides only depend on the
generalized permutation and of $\{Im(\zeta_{\alpha})\}_{\alpha\in
\mathcal{A}}$.

\begin{figure}[htbp]
\psfrag{a}{$\scriptstyle R_{\{1,10\}}$}  \psfrag{b}{$\scriptstyle R_{\{2,4\}}$}
\psfrag{1}{$\scriptstyle 1$}  \psfrag{2}{$\scriptstyle 2$}
   \begin{center}
     \includegraphics{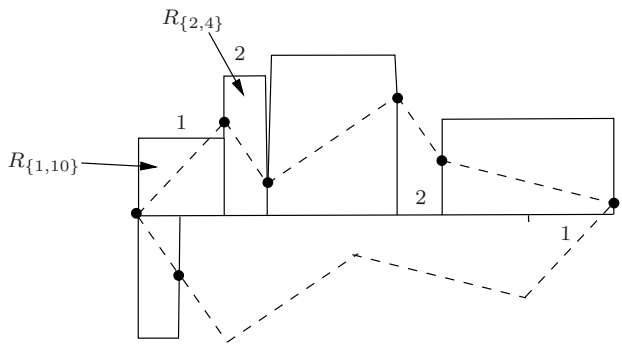}
    \caption{Zippered rectangle construction of the flat surface of
     Figure~\ref{figure:suspension:data}.}
     \label{polygtozip}
   \end{center}
\end{figure}

For the general case, we construct the rectangles $R_{\alpha}$ using
the same formulas. Identifications for the horizontal sides are 
straightforward. Identifications for the vertical sides, that do not
depends on the horizontal parameters, will be well defined after the
following lemma.

\begin{lem}
Let $\zeta$ be a suspension data for a linear involution 
$T$, and let $\gp$ be the corresponding generalized
permutation. There exists a linear involution $T'$ and a suspension data 
$\zeta^{\prime}$ for $T^{\prime}$ such that:
\begin{itemize}
\item The generalized permutation associated to $T^{\prime}$ is $\gp$.
\item For any $\alpha$ the complex numbers $\zeta_{\alpha}$ and $\zeta'_{\alpha}$ 
have the same imaginary part.
\item The suspension data $\zeta^{\prime}$ defines a suitable polygon.
\end{itemize}
\end{lem}

\begin{proof}
We can assume that $\sum_{k=1}^{l} Im(\zeta_{\pi(k)})>0$ (the negative
case is analogous and there is nothing to prove when the sum is
zero). It is clear that $\sigma(l+m)\neq l$ otherwise there would be
no possible suspension data. If $\sigma(l+m)<l$, then we can shorten
the real part of  $\zeta_{\gp(l+m)}$, keeping conditions (1)--(4)
satisfied, and get a suspension data $\zeta^{\prime}$ with the same
imaginary part as $\zeta$, and such that
$Re(\zeta^{\prime}_{\gp(l+m)})<Re(\zeta^{\prime}_{\gp(l)})$. This last
condition implies that $\zeta^{\prime}$ defines a suitable polygon.

If $\sigma(l+m)>l$, then condition $(4)$ implies that
$Re(\zeta_{\gp(l+m)})$ is necessary bigger than
$Re(\zeta_{\gp(l)})$. However, we can still change $\zeta$ into a
suspension data $\zeta^{\prime}$, with same imaginary part, and such
that $Re(\zeta^{\prime}_{\gp(l+m)})$ is very close to
$Re(\zeta^{\prime}_ {\gp(l)})$. In that case, $\zeta^{\prime}$ also
defines a suitable polygon. See \cite{B}, Lemma~2.1 for more details.
\end{proof}

We have therefore defined the zippered rectangle construction for any
suspension data. Note that we have not yet discussed  the existence
of a suspension data. This will be done in the upcoming section. This notion 
is natural. See~\cite{Ve82} and the following Proposition.

\begin{prop}
\label{prop:irr:natural}
Let $S$ be a flat surface with no vertical saddle connections and let $X$ be a
horizontal interval attached to a singularity on the left. Let $\gamma$ be the vertical 
leaf passing through the right endpoint of $X$, we assume that
$\gamma$ meets a singularity before returning to $X$, in positive or
negative direction. Let $T=(\pi,\lambda)$ be  the 
linear involution given by the cross section on $X$ of the vertical flow. There 
exists a suspension data $\zeta$ such that $(\pi,\zeta)$ defines a surface 
isometric to $S$.
\end{prop}

\begin{proof}
See the construction given in the proof of Proposition $2.2$ in~\cite{B}.
\end{proof}

We can define the Rauzy-Veech induction on the space of suspensions, as
well as on the space of zippered rectangles. Let $T=(\gp,\lambda)$ be
a linear involution and let $\zeta$ be a suspension over $T$. Then we
define $\mathcal R(\gp,\zeta)=(\gp',\zeta')$ as follows.

\begin{itemize}
\item 
If $T=(\gp,\lambda)$ has type $0$, then
$\mathcal{R}(\gp,\zeta)=(\mathcal{R}_0\pi,\zeta^{\prime})$, with
$\zeta^{\prime}_\alpha=\zeta_\alpha$ if $\alpha\neq \pi(l)$ and
$\zeta_{\pi(l)}^{\prime}=\zeta_{\pi(l)}-\zeta_{\pi(l+m)}$. 
\item
If  $T=(\gp,\lambda)$ has type $1$, then
$\mathcal{R}(\gp,\zeta)=(\mathcal{R}_1\pi,\zeta^{\prime})$, with
$\zeta^{\prime}_\alpha=\zeta_\alpha$ if $\alpha\neq \pi(l+m)$ and
$\zeta_{\pi(l+m)}^{\prime}=\zeta_{\pi(l+m)}-\zeta_{\pi(l)}$. 
\end{itemize}

We can show that $(\gp^{\prime},\zeta^{\prime})$ is a suspension over $\mathcal{R}(T)$ 
and defines a surface isometric to the one corresponding to $(\gp,\zeta)$. 

As in the case of interval exchange maps we consider the 
renormalized Rauzy-Veech induction defined on lengths one intervals:
$$
\textrm{if } \mathcal R(\pi,\lambda) = (\pi',\lambda') \textrm{ then } 
\mathcal R_r(\pi,\lambda):= (\pi',\lambda'/|\lambda'|).
$$
One can define obviously the corresponding renormalized Rauzy-Veech induction on 
the suspensions data by contracting the imaginary parts by a factor $|\lambda'|$ 
which preserves the area of the corresponding flat surface.

\section{Geometry of generalized permutations}
\label{combinatoric}

In this section we give a necessary and sufficient condition for a 
generalized permutation to admit a suspension; this will prove the
first part of Theorem~$A$. Let us first introduce some notations to
make clear the definition.

\emph{Notation:}
If $\mathcal{A}=\{\alpha_1,\dots,\alpha_d\}$ is an alphabet, we will denote by 
$\mathcal A \sqcup \mathcal A$ the set with multiplicities
$\{\alpha_1,\alpha_1,\dots,\alpha_d,\alpha_d\}$ of cardinal $2d$, and
we will use analogous notations for subsets of $\mathcal{A}$. \medskip

We will also call \emph{top} (respectively \emph{bottom}) the
restriction of a generalized permutation $\gp$ to $\{1,\ldots,l\}$
(respectively $\{l+1,\ldots,l+m\}$) where $(l,m)$ is the type of
$\gp$. \medskip

\emph{Notation:}
Let $F_1,F_2,F_3,F_4$ be  (possibly empty) unordered subsets of
$\mathcal{A}$ or $\mathcal A \sqcup \mathcal A$.  We say that a
generalized permutation $\gp$ of type $(l,m)$ is decomposed if
\[
\gp=
\left(\begin{array}{c|c|c}F_1 & *** & F_2 \\\hline F_3& *** &F_4
  \end{array}\right),
\]
and there exist $0 \leq i_1 \leq i_2 \leq l$  and $l \leq i_3\leq i_4
\leq l+m=2d$ such that
\begin{itemize}
\item $\{\gp(1),\dots,\gp(i_1)\}=F_1$
\item $\{\gp(i_2),\dots,\gp(l)\}=F_2$
\item $\{\gp(l+1),\dots,\gp(i_3)\}=F_3$
\item $\{\gp(i_4),\dots,\gp(2d)\}=F_4$. 
\end{itemize}
The sets $F_1, F_2, F_3$, and $F_4$ will be referred as top-left, top-right,
bottom-left and bottom-right corners respectively.

We do not assume that $card(F_1)=card(F_3)$, or $card(F_2)=card(F_4)$.
\begin{defs}
\label{def:irred}
We will say that $\gp$ is {\it reducible} if $\gp$ admits a decomposition
\begin{equation*}(*)\quad
\left(\begin{array}{c|c|c}A\cup B & *** & D\cup B \\\hline A\cup C &
  *** & D \cup C  \end{array}\right), \ A,B,C,D \textrm{
  disjoint subsets of } \mathcal{A},
\end{equation*}

where the subsets $A,B,C,D$ are not all empty and 
one of the following statements holds
\begin{enumerate}
\item[i-] No corner is empty
\item[ii-] Exactly one corner is empty and it is on the left.
\item[iii-] Exactly two corners are empty and they are either both on the left, either both on the right.
\end{enumerate}
A permutation that is not reducible is {\it irreducible}.
\end{defs}

The main result of this section is the next theorem which, being
combined with Proposition~\ref{prop:irr:natural} implies
first part of Theorem~$A$. We make clear that in this section, we only speak of
suspensions given by the construction of section \ref{susp:giem}.

\begin{ths}
\label{CNS}
Let $T=(\pi,\lambda)$ be a linear involution. 
Then $T$ admits a suspension $\zeta$ if and only if the underlying generalized 
permutation $\pi$ is irreducible.
\end{ths}

\begin{rem}
Note that the existence or not of a suspension is independent of the length data $\lambda$.
\end{rem}

\begin{rem}
 One can see that this reducibility notion is not symmetric with
 respect to the left/right, contrary to the case of interval exchange
 maps. Therefore, the choice of attaching a singularity on the left
 end of the segment in the construction of section \ref{susp:giem} is
 a real choice. This will have an important consequence in terms of
 extended Rauzy classes.
\end{rem}

\begin{rem}
\label{rem:strongly:irreducible}
In the usual case of interval exchange maps, one can always
choose $\zeta$ in such a way that $Im(\sum_{i=1}^l \zeta_{\gp(i)})=0$
(\emph{i.e.} there is a singularity on the left and on the
right of the interval). Here it is not always possible. More precisely one can show
that $T$ admits such a suspension with this extra condition if and
only if for any decomposition of $\gp$ as in equation~$(*)$ above, all
the corners are empty.
\end{rem}

\subsection{Necessary condition}

\begin{prop} 
\label{CN}
A reducible generalized permutation does not admit any suspension data. 
\end{prop}

\begin{proof}[Proof of the Proposition]
Let us consider $\pi$ a reducible generalized permutation. It is convenient to
introduce some notations. Let us assume that there exists a suspension
$\zeta$ over $\pi$. Then we define $a$ the real number
$a=\sum_{j\in A} Im(\zeta_{\gp(j)})$; we define $a=0$ if the set $A$
is empty. Finally we define $b$, $c$ and $d$ in an analogous manner for
$B,C$ and $D$. We also define  $t=\sum_{i=1}^l Im(\zeta_{\gp(i)})$. We
distinguish three cases following Definition~\ref{def:irred}. \medskip

\noindent {\it i- No corner is empty.} \\
Then the following inequalities hold

$$
\left\{
\begin{array}{lll}
a+b &>&0 \\
a+c &<&0 \\
t-d-b &>&0 \\
t-d -c &<&0
\end{array}
\right.
$$

Subtracting the second one from the first one, and the fourth one from the third one, we get:
$$
\left\{
\begin{array}{lll}
b-c &>&0 \\
c-b &>&0
\end{array}
\right.
$$
which is a contradiction. \medskip

\noindent {\it ii- Exactly one corner is empty, and it is on the left.} \\
We can assume without loss of generality  that it is  the
top-left one. That means that $A$, $B$ are empty, and $C$, $D$ are
nonempty. Therefore the following inequalities holds: 
$$
\left\{
\begin{array}{lll}
c &<&0 \\
t-d &>&0\\
t-d-c &<&0
\end{array}
\right.
$$
Subtracting the third  inequality from the second one, we get $c>0$, which contradicts
the first one.\medskip 

\noindent {\it iii- Exactly two corners are empty.} \\
If they are both on the left side, then we have $B$ and $C$ empty and
$D$ non empty. This implies that $t-d$ is both positive and negative,
which is impossible. \\
If they are both on the right side, it is similar. 
If the two corners forming  a diagonal were empty, then it is easy to see
that all the corners would be empty, hence this case doesn't occur by
assumption. The proposition is proven.
\end{proof}

\subsection{Sufficient condition}
\label{suff:cond}

In this section, we will not necessarily assume that generalized permutations satisfy 
Convention~\ref{conv:normal}, since
for technical reasons, some intermediary results of this section must be stated 
for an arbitrary generalized permutation.

We will have to work only on the imaginary part of the $\zeta_i$ in
order to built a suspension. Hence, in order to simplify the notations
we will use the following ones. We will use this vocabulary only in this section.

\begin{defs} \label{def:pseudo}
A \emph{pseudo-suspension} is a collection of real numbers
$\{\tau_i\}_{i\in \mathcal{A}}$ such that: 
\begin{itemize}
\item For all $k\in\{1,\ldots,l\}\quad\sum_{i\leq k} \tau_{\gp(i)}\geq 0$ .
\item For all $k\in\{1,\ldots,m\}\quad \sum_{l<i\leq l+k} \tau_{\gp(i)}\leq 0$ .
\item $\sum_{i\leq l} \tau_{\gp(i)}=\sum_{l<i\leq l+m} \tau_{\gp(i)}=0$
\end{itemize}

A pseudo-suspension is \emph{strict} if all the previous inequalities
are strict except for the extremal ones.

A \emph{vanishing index} on the top (respectively bottom) of a
pseudo-suspension is an integer $k_0 < l$ (respectively $k_0 < m$) such
that $\sum_{i\leq k_0} \tau_{\gp(i)}=0$ (respectively $ \sum_{l<i\leq l+k_0} \tau_{\gp(i)}=0$).

A pseudo-suspension $\tau^{\prime}$ is \emph{better} than $\tau$ if the set
of vanishing indices of $\tau'$ is strictly included into the set of
vanishing indices of $\tau$. 

We will say that $\gp$ is {\it strongly irreducible} if for any
decomposition of $\gp$ as in $(*)$ of Definition~\ref{def:irred}, all the
corners are empty. Of course strong irreducibility implies
irreducibility. 
\end{defs}

The following lemma is obvious and left to the reader.

\begin{lem}\label{lm:strongly:irred}    
Let $\gp$ be generalized permutation satisfying Convention~\ref{conv:normal} 
that admits a strict pseudo-suspension. Then $\gp$ admits a suspension $\zeta$
with $Im(\sum_{1\leq i \leq l} \zeta_{\pi(i)})=0$. 
\end{lem}

Let us assume that $\pi$ is any irreducible permutation. One has to
find a suspension 
$\zeta$ over $\pi$. We will first assume that $\pi$ is strongly
irreducible and we will show that $\pi$ admits such a suspension with
the extra equality $Im(\sum_{1\leq i \leq l} \zeta_{\pi(i)})=0$. This
corresponds to a special case of Proposition~\ref{max:donc:reducible:star}.
We will then relax the condition on the irreducibility of $\pi$ and
prove our main result. Note that one can extend
the proof of Proposition~\ref{CN} to show that if $\zeta$ is a
suspension data such that $Im(\sum_{1\leq i \leq l}
\zeta_{\pi(i)})=0$, then $\pi$ is strongly irreducible. 
\medskip

From Lemma~\ref{lm:strongly:irred} we have reduced the problem to
the construction of a strict pseudo-suspension. 
As we have seen in section \ref{sec:iem}, in the case of true
permutations, there is an explicit formula, due to Masur and Veech,
that gives a 
suspension when the permutation is irreducible. We will first build a
pseudo-suspension $\tau_{MV}$ by extending this formula to generalized
permutations. This will not give in general a strict pseudo-suspension. \medskip

Let $\gp: \{1,\ldots l+m\} \rightarrow \mathcal{A}$ be a generalized
permutation. We can decompose $\mathcal{A}$ into three disjoint subsets
\begin{itemize}
\item The subset $\mathcal{A}_{01}$ of elements $\alpha \in \mathcal{A}$
  such that $\gp^{-1}(\{\alpha \})$ contains exactly one element in
  $\{1,\ldots,l\}$ and one element in $\{l+1,\ldots,l+m\}$. 
The restriction of $\gp$ on $\gp^{-1}(\mathcal{A}_{01})$ defines a true permutation.
\item The subset $\mathcal{A}_0$ of elements $\alpha \in \mathcal{A}$ such
  that $\gp^{-1}(\{\alpha\})$ contains exactly two elements in
  $\{1,\ldots,l\}$ (and hence no elements in $\{l+1,\ldots,l+m\}$). 
\item The subset $\mathcal{A}_1$ of elements $\alpha \in \mathcal{A}$ such
  that $\gp^{-1}(\{\alpha\})$ contains exactly two elements in $\{l+1,\ldots,l+m\}$ (and hence no elements in   $\{1,\ldots,l\}$). 
\end{itemize}

The next lemma is just a reformulation of the construction of a suspension data in section~\ref{susp:data:iem}
\begin{lem}[Masur; Veech]
\label{solMV}
Let $\pi$ be a true permutation defined on $\{1,\ldots,d\}$, then the
integers $\tau_i = \pi(i)-i$ for $1\leq i \leq d$ define a pseudo-suspension
over $\pi$. Furthermore, we have:
$$
\sum_{i\leq i_0}\tau_i=0 \Leftrightarrow \sum_{i\leq i_0}
\tau_{\pi^{-1}(i)}=0 \Leftrightarrow
\pi(\{1,\ldots,i_0\})=\{1,\ldots,i_0\}.
$$
\end{lem}

Recall that we do not assume any more that a generalized permutation
satisfies Convention~\ref{conv:normal}.

\begin{lem}
\label{solMV2}
Let $\gp$ be a generalized permutation of type $(l,m)=(2d,0)$ and $\sigma$ the 
associated involution. There exists a collection of real numbers $(\tau_1,\ldots,\tau_{2d})$
with $\sum_{i\leq i_0} \tau_i \geq 0$ for all $i_0$ and such that
$$\sum_{i\leq i_0} \tau_i=0 \Leftrightarrow
\sigma(\{1,\ldots,i_0\})=\{2d,\ldots, 2d-i_0+1\}).
$$
\end{lem}

\begin{proof}
We will construct from $\gp_0:=\gp$ a new permutation $\tilde \pi$ on 
$d$ symbols. Let us consider the ``mirror symmetry'' $\gp_1$ of $\gp_0$ as follows. In 
tabular representation $\gp_0$ is $(\tau(1),\dots,\tau(2d))$; $\gp_1$ is of type $(0,2d)$ and 
its tabular representation is $(\tau(2d),\dots,\tau(1))$.

Then $\tilde \pi$ is in tabular representation $\left( \begin{smallmatrix} L_0 \\ L_1 
\end{smallmatrix} \right)$ with $L_i$ is obtained from $\gp_i$ by removing the second occurrence 
of each letter. For instance, if $\gp_0=(A\ B\ C\ C\ D\ D\ A\ B)$ then 
$\gp_1=(B\ A\ D\ D\ C\ C\ B\ A)$ and 
$\tilde \pi = \left( \begin{smallmatrix} A & B & C & D \\ B & A & D & C 
\end{smallmatrix} \right)$.
It is easy to check that $\tilde \pi$ is reducible 
if and only if there exists $i_0$ such that
$\sigma(\{1,\ldots,i_0\})=\{2d-i_0+1,\dots,2d\}$. Moreover 
the solution of Lemma~\ref{solMV} gives the
desired collection of numbers $\tau_i$.
\end{proof}

\begin{defslem}
We define the pseudo-suspension $\tau_{MV}$ over $\pi$ by the
collection of real numbers given by
\begin{itemize}
\item The solutions given by Lemma~\ref{solMV} and Lemma~\ref{solMV2}
  for the restrictions of $\gp$ on $\gp^{-1}(\mathcal{A}_{01})$ and on $\gp^{-1}(\mathcal{A}_0)$. 
\item The solution of Lemma~\ref{solMV2} for the
  restriction of $\gp$ on  $\gp^{-1}(\mathcal{A}_1)$, taken with opposite sign.
\end{itemize}
\end{defslem}

\begin{lem}
\label{lm:dec}
Let $k\in \{1,\ldots,l\}$ be any vanishing index on the top of $\tau_{MV}$. 
Setting $A=\gp(\{1,\ldots, k\}) \cap \mathcal{A}_{01}$ and $B=\gp(\{1,\ldots,k\})\cap
(\mathcal{A}_0 \sqcup \mathcal{A}_0)$, there exists $C\subset
\mathcal{A}_1 \sqcup \mathcal{A}_1$ and $D\subset \mathcal{A}_{01}$
such that the generalized permutation $\gp$ decomposes as
$$
\left(\begin{array}{c|c|c}A\cup B & *** & D\cup B^{\prime} \\\hline A\cup C & *** & *** 
\end{array}\right)
$$
with $A\cup B\neq \emptyset$ and 
with one of the following properties: 
either $B= B^{\prime}\subset \mathcal{A}$ or there exist 
$i_1,i_2\leq k$ such that $\gp(i_1)=\gp(i_2)\in B$ and $B'\subset B$.\medskip

There is an analogous decomposition for vanishing indices in
$\{l+1,\ldots,l+m\}$ but with different subsets $A^{\prime}$
$B^{\prime}$, $C^{\prime}$ and $D^{\prime}$ \emph{a priory}.
\end{lem}

\begin{proof}

It follows from Lemmas~\ref{solMV} and~\ref{solMV2}.

\end{proof}

\begin{rem}
\label{flip}
If $\tau$ is a pseudo-suspension of $\gp =
\left(\begin{smallmatrix}\alpha_1 & \alpha_2 &*** & \alpha_l \\
  \alpha_{l+1} & \alpha_{l+2} & *** & \alpha_{l+m}
\end{smallmatrix}\right)$ then $\tau^{\prime}=-\tau$ is a pseudo-suspension
of $\gp^{\prime} =  \left(\begin{smallmatrix} \alpha_{l+1} &
  \alpha_{l+2} & *** & \alpha_{l+m}\\ \alpha_1 & \alpha_2 &*** &
  \alpha_l \end{smallmatrix}\right),$ and $\tau$ is a
  pseudo-suspension of the generalized permutation 
  $\gp^{\prime\prime}=\left(\begin{smallmatrix}\alpha_l & \alpha_{l-1}
  &*** & \alpha_1 \\ \alpha_{l+m} & \alpha_{l+m-1} & *** &
  \alpha_{l+1} \end{smallmatrix}\right)$. \\
Hence we can ``flip'' the generalized permutation $\gp$ by top/bottom
or left/right without loss of generality.
\end{rem}
 
In the next two lemmas, we denote by $\tau$ a pseudo-suspension
that is better than $\tau_{MV}$ and maximal (i.e. there is no better
pseudo-suspensions).

\begin{lem} 
\label{lemmeA}
Let $i_1$ and $i_2$ be the two first top and bottom vanishing
indices for $\tau$ (possibly $i_1=l,i_2=m$). Let $A=\gp(\{1,\ldots,i_1\})\cap \mathcal{A}_{01}$ and
$A^{\prime}=\gp(\{l+1,\ldots,l+i_2\})\cap \mathcal{A}_{01}$. Then either
$A=A^{\prime}$ or $A=\emptyset$ or $A^{\prime}=\emptyset$.
\end{lem}

\begin{figure}[htb]
\begin{center}
 \psfrag{x}[][]{$A=\{1,2\}$}
  \psfrag{y}[][]{$A'=\{1,2,3,4\}$}
   \psfrag{1ps}[][]{\tiny $\gp(j_1)=1$}
    \psfrag{4ps}[][]{\tiny $4=\gp(j_2)$}
\includegraphics[scale=0.6]{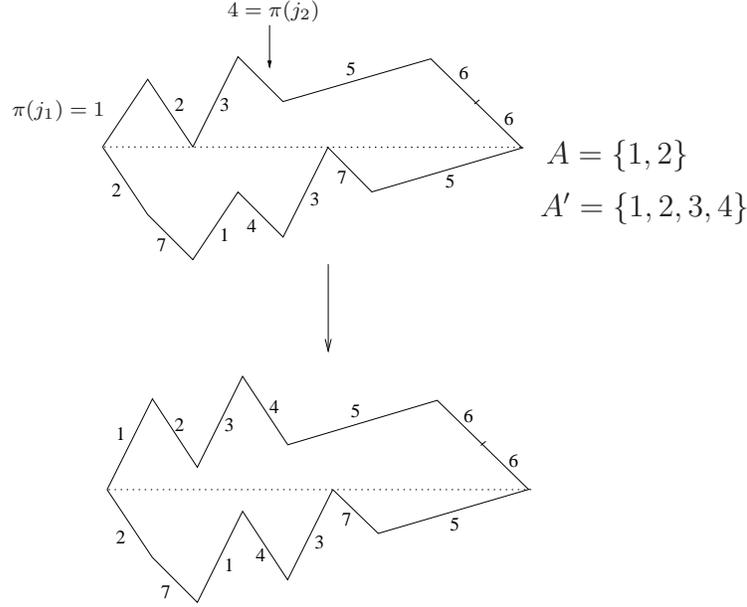}
\caption{Construction of a pseudo-suspension $\tau^{\prime}$ better than $\tau$.}
\label{exemple1}
\end{center}
\end{figure}

\begin{proof}
We assume that neither $A$ nor $A^{\prime}$ is
empty. Lemma~\ref{lm:dec} implies that one of this set is a subset of
the other one. 

Without loss of generality, we can assume that $A\subseteq
A^{\prime}$. Let us assume $A\not = A'$; we will get a contradiction.
So there exist $j_1,j_2$ in $\gp^{-1}(\mathcal{A}_{01})$ such that $1\leq j_1\leq i_1 <j_2\leq l$.
But by definition of $A$ and $A^{\prime}$, we also have $\sigma(j_1) < \sigma(j_2)$.

The definition of $i_2$ implies that there exists $c<0$ such that, for
$l+1\leq k<i_2$, the following inequality holds:
$$
\sum_{1+l\leq i \leq l+k} \tau_{\gp(i)} \leq c<0.
$$
Now we replace $\tau_{\gp(j_1)}$  (respectively $\tau_{\gp(j_2)}$) by
$\tau_{\gp(j_1)}- \frac{c}{2}$ (respectively  $\tau_{\gp(j_2)}+ \frac{c}{2}$) and
get a vector $\tau^{\prime}$, see Figure~\ref{exemple1}. We have
\begin{itemize}
\item $\sum_{1\leq i \leq k} \tau^{\prime}_{\gp(i)}>0$ for $k<j_2$. 
\item $\sum_{1\leq i \leq k} \tau^{\prime}_{\gp(i)}=\sum_{1\leq i \leq k} \tau_{\gp(i)}$ for $k\geq j_2$.
\item $\sum_{1+l\leq i \leq k} \tau^{\prime}_{\gp(i)}\leq c/2<0$ for $l+1\leq k<\sigma(j_2)$.
\item $\sum_{1+l\leq i \leq k} \tau^{\prime}_{\gp(i)}=\sum_{1+l\leq i
  \leq k} \tau_{\gp(i)}$ for $k\geq \sigma(j_2)$ (since
  $\sigma(j_1)<\sigma(j_2)$). 
\end{itemize}
Hence, $\tau^{\prime}$ is a pseudo-suspension better than $\tau$,
contradicting its maximality. Therefore $A=A'$ and the lemma is proven.
\end{proof}

\begin{lem} \label{lemmeB}
Let $i_1$ and $i_2$ be the first and last  top vanishing indices of $\gp$.
Let $B=\gp(\{1,\ldots,i_1\})\cap (\mathcal{A}_0 \sqcup \mathcal{A}_0)$ and $B^{\prime}
=\gp(\{i_2+1,\ldots,l\})\cap (\mathcal{A}_0 \sqcup \mathcal{A}_0)$. Then either $B^{\prime}=B\subset \mathcal{A}_0$ or
$B=\emptyset$ or $B'=\emptyset$. Moreover if there
exist $i_{b_1}\neq i_{b_2}$ in $\{1,\ldots,i_1\}$ such that
$\gp(i_{b_1})=\gp(i_{b_2})$ then $B=\mathcal{A}_0 \sqcup \mathcal{A}_0$.
\end{lem}

\begin{proof}
We sketch the proof here. We assume that there exist $i_{b_1}$ and $i_{b_2}$ in
$\{1,\ldots,i_1\}$ such that $\gp(i_{b_1})=\gp(i_{b_2})$. If there
exists $i_{b_3}>i_1$ such that $\gp(i_{b_3})\in B$, then we set: 
\begin{eqnarray*}
\tau^{\prime}_{\pi (i_{b_1})}=\tau_{\pi (i_{b_1})} +\varepsilon\\
\tau^{\prime}_{\pi (i_{b_3})}=\tau_{\pi (i_{b_3})} -\varepsilon.
\end{eqnarray*}
Then is is easy to see that, for $\varepsilon$ small enough,
$\tau^{\prime}$ is a pseudo-suspension and is better than $\tau$,
contradicting its maximality. Remark~\ref{flip} implies that the same
statement is true for $B'$; hence, we can assume that $B,B' \subset
\mathcal{A}_0$. We conclude using the same argument as the one of the
proof of the previous Lemma~\ref{lemmeA}.
\end{proof}

\begin{prop}
\label{max:donc:reducible:star}
Let  $\gp$ be a  strongly irreducible generalized permutation. Let $\tau$ be any 
pseudo-suspension which is better than $\tau_{MV}$ and maximal. Then 
$\tau$ is a strict pseudo-suspension.
\end{prop}

\begin{proof}[Proof of Proposition~\ref{max:donc:reducible:star}]
Let us assume that $\tau$ is not strict. From Lemmas~\ref{lemmeA}
and~\ref{lemmeB} and Remark~\ref{flip} we have the following
decomposition of $\gp$.
$$
\left(\begin{array}{c|c|c}A\cup B & *** & D\cup B^{\prime} \\\hline
  A^{\prime}\cup C & *** & D^{\prime} \cup C^{\prime}  
\end{array}\right)
$$
with $A,A',D,D'\subset \mathcal{A}_{01}$, $B,B'\subset \mathcal{A}_0$ and $C,C'\subset \mathcal{A}_1$ by assumption, and 
with the condition that either $A,A^{\prime}$ are  equal, or at least
one of them is empty (and similar statement for the pair
$(D,D^{\prime})$); and the condition that if $B,B' \subset
\mathcal{A}_0$ then they are either equal, or at least one of them is
empty, otherwise one of them is  $\mathcal{A}_0 \sqcup \mathcal{A}_0$
(and similar statements for $C,C'$). By
convention from now on, we will keep the notation $B$ or $C$ only when
they are not equal to $\mathcal{A}_0 \sqcup \mathcal{A}_0$ or 
$\mathcal{A}_1 \sqcup \mathcal{A}_1$, and therefore subsets of $\mathcal{A}_0$ or $\mathcal{A}_1$. \medskip

Let us note that if there is no vanishing index in $\{1,\dots,l-1\}$ or in
$\{l+1,\dots,l+m-1\}$, the corresponding right corner is just
empty. But if $\tau$ is not strict, then there exists at least a pair of
nonempty corners in the top or in the bottom. \medskip

If there is a vanishing index on the top, then the two corresponding
corners are non-empty. Then it is easy to see that either there is a
corner with only $A, B$ or $D$, or the corners are respectively $A\cup
B$ or $D\cup B$, with $A,B,D$ nonempty. In this case
Lemma~\ref{lemmeA} implies that there must be a vanishing index in
$\{l+1,\dots,l+m\}$.  

Since there must be a vanishing index in the top, or in the bottom,
the previous argument implies that either~$\gp$ is not strongly
irreducible, or there is one corner that only consists of one set
$A,B,C$ or $D$. Thanks to Remark~\ref{flip}, we assume that this is
the top-left corner; this leads to the two next cases. \medskip 

The general idea of the next part of the proof is first to remove the
cases that correspond to not strongly irreducible permutations, and
then show that the other cases correspond to a non-maximal
pseudo-solution. \medskip 

\noindent {\it First case: The top-left corner is $B$.} \\
There is  necessary a vanishing index in $\{1,\dots,l-1\}$, and hence 
the top-right corner is not empty. It also does not contains all
$\mathcal{A}_0 \sqcup \mathcal{A}_0$, hence it is necessary $B$, $D$
or $D\cup B$. Recall that $\pi$ is assumed to be strongly irreducible,
so the top-right corner is not $B$.  If the bottom-right corner were
$D$, the generalized permutation $\gp$ would decompose as  
$$
\left(\begin{array}{c|c|c}B & *** & D\cup B \\\hline  & *** & D 
\end{array}\right),
$$
or 
$$
\left(\begin{array}{c|c|c}B & *** & D \\ \hline & *** & D 
\end{array}\right)
$$
which are not strongly irreducible. Hence the bottom-right corner is not $D$. This also implies that $\mathcal{A}_1$ cannot be empty.
\medskip

Let us assume that there are no vanishing indices in the bottom line.
 We choose any
element $b\in B$, $c\in \mathcal{A}_1$, and $d\in D$ and change $\tau_b$ by
$\tau_b+\varepsilon$, $\tau_c$ by $\tau_c+\varepsilon$ and $\tau_d$ by
$\tau_d-2\varepsilon$. If $\varepsilon$ is small enough, then the new
vector $\tau^{\prime}$ is better than $\tau$, which contradicts its
maximality. \medskip

So, the bottom admits vanishing indices; then the
bottom-left corner can be $C,\mathcal{A}_1 \sqcup \mathcal{A}_1,\mathcal{A}_1 \sqcup \mathcal{A}_1 \cup A,A$ or $A\cup
C$. Let us discuss these cases in details.

\begin{itemize}
\item $C$: the bottom-right corner is $C$, $D$ or $C\cup
D$. In the first and second cases, $\gp$ is not strongly irreducible. If for
instance, the top-right is  $D\cup B$, then $\gp$ decomposes as   
$$
\left(\begin{array}{c|c|c}B & *** & D\cup B \\\hline C& *** & D \cup C
\end{array}\right),
$$
and therefore $\gp$ is not strongly irreducible. The other case is
similar.

\item $\mathcal{A}_1 \sqcup \mathcal{A}_1$ or $\mathcal{A}_1 \sqcup \mathcal{A}_1 \cup A$: in that case, the bottom-right
  corner is necessary $D$ and we have already proved that $\pi$ is not strongly irreducible in this situation.

\item $A$ or $A\cup C$: 
We construct a better pseudo-suspension $\tau'$. \\
Let $j_1\leq l$ be the smallest index such that $\sigma(j_1)>l$ and
let $j_2\leq l$ be the largest one. Let $i_1$ be the first vanishing
index. There exists $j_b\in \{1,\dots,i_1\}$ such that $\sigma(j_{b})<
j_2$ otherwise the top-line would have a decomposition as $B|***|B$,
and $\gp$ would be not strongly irreducible.  Let $j_c$ be the first
index in $\gp^{-1}(\mathcal{A}_1)$ (see Figure~\ref{exemple2}).

Now we define $\tau^{\prime}$ in the following way:
\begin{eqnarray*}
\tau^{\prime}_{\gp(j_1)}=\tau_{\gp(j_1)}-\varepsilon\\
\tau^{\prime}_{\gp(j_2)}=\tau_{\gp(j_2)}-\varepsilon\\
\tau^{\prime}_{\gp(j_b)}=\tau_{\gp(j_b)}+\varepsilon\\
\tau^{\prime}_{\gp(j_c)}=\tau_{\gp(j_c)}+\varepsilon\\
\forall \alpha \notin \gp(\{j_1,j_2,j_b,j_c\}) \quad
\tau^{\prime}_{\alpha}=\tau_{\alpha}.
\end{eqnarray*}
In the extremal case $j_1=j_2$, the following arguments will work
similarly if we define $\tau^{\prime}_{\gp(j_1)}$ by
$\tau_{\gp(j_1)}-2\varepsilon$. We have 
\begin{eqnarray*}
\forall k\in \{1,\ldots,l\}\quad \sum_{i=1}^{k} \tau_{\gp(i)}^{\prime}=\sum_{i=1}^{k} \tau_{\gp(i)}+n_k \varepsilon\\
\forall k\in \{1,\ldots,m\}\quad \sum_{i=l+1}^{l+k}
\tau_{\gp(i)}^{\prime}=\sum_{i=l+1}^{l+k} \tau_{\gp(i)}+m_k
\varepsilon 
\end{eqnarray*}
Here $n_k$ is the difference between the number of indices in
$\{j_b,\sigma(j_b)\}$ smaller than or equal to $k$, and number of
indices in $\{j_1,j_2\}$ smaller than or equal to $k$. This value is
\emph{always} greater than or equal to zero for $k\in \{1,\dots,l\}$,
and is strictly greater than zero when $k$ is the first vanishing
index.

Similarly $m_k$ is the difference between the number of indices in
$\{j_c,\sigma(j_c)\}$ that are in $\{l+1,\ldots,k\}$, and number of
indices in  $\{\sigma(j_1),\sigma(j_2)\}$ that are in
$\{l+1,\ldots,k\}$. This value might be positive.  
Let $i_3\leq i_4<l+m$ be respectively the first and last bottom
vanishing indices. We have the following facts:
\begin{itemize}
\item $\sigma(j_1)\leq i_3$ otherwise the bottom-left corner is $C$.
\item $\sigma(j_c)>i_4$ otherwise the bottom-right corner is $D$.
\end{itemize}
Hence it is easy to check that $m_k$  can be strictly positive only
for $l<k<i_3$ or $i_4<k<l+m$. 

Then if $\varepsilon$ is small enough, $\tau^{\prime}$ is a
pseudo-suspension, and is better than $\tau$ (see Figure~\ref{exemple2}),
which contradicts the maximality of $\tau$.
\end{itemize}

\begin{figure}[htbp]
\psfrag{1}{$\scriptstyle 1$} \psfrag{2}{$\scriptstyle 2$}
\psfrag{3}{$\scriptstyle 3$} \psfrag{4}{$\scriptstyle 4$}
\psfrag{5}{$\scriptstyle 5$} \psfrag{6}{$\scriptstyle 6$}
\psfrag{7}{$\scriptstyle 7$} 
\psfrag{a}{$\scriptstyle \pi(\sigma(j_b))=1$} 
\psfrag{b}{$\scriptstyle \pi(j_1)=2$} 
\psfrag{c}{$\scriptstyle \pi(j_b)=1$} 
\psfrag{d}{$\scriptstyle \pi(j_2)=5$} 
\psfrag{e}{$\scriptstyle \pi(j_c)=6$} 
\psfrag{B}{$B$} \psfrag{BuD}{$ B\cup D$} 
\psfrag{C}{$C$} \psfrag{AuC}{$ A\cup C$} 

   \begin{center}
     \includegraphics{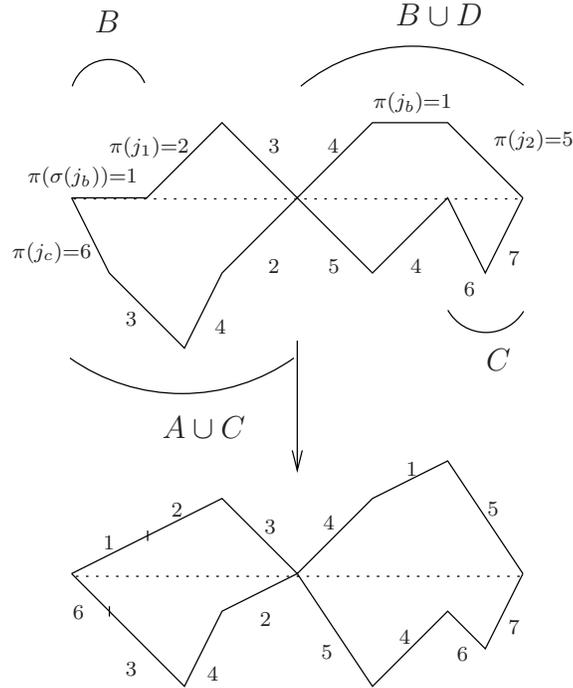}
    \caption{Construction a pseudo-solution $\tau^{\prime}$ better than $\tau$.}
     \label{exemple2}
   \end{center}
\end{figure}

\noindent {\it Second case: The top-left corner is $A$.} \\
We assume that there are no corners $B$ or $C$, since this case has already been discussed. \medskip

Let us assume that there is no vanishing index in the bottom
line. Then, according to Lemma~\ref{lemmeA}, $A=\mathcal{A}_{01}$; therefore the
top-right corner is $\mathcal{A}_0 \sqcup \mathcal{A}_0$ or $B$. If $\mathcal{A}_1$ is empty, then
$\gp$ decomposes as 
\[
\left(\begin{array}{c|c}A & \mathcal{A}_0 \sqcup \mathcal{A}_0\\ \hline & A  
\end{array}\right)
\]
so $\gp$ is not strongly irreducible. If $\mathcal{A}_1$ is not empty,
we choose any element $a\in A$, $b\in \mathcal{A}_0$, $c\in
\mathcal{A}_1$, and replace $\tau_a$ by $\tau_a+2\varepsilon$,
$\tau_b$ by $\tau_b-\varepsilon$, and $\tau_c$ by
$\tau_c-\varepsilon$. This new pseudo-suspension we have constructed is better
that the old one for $\varepsilon$ small enough. 
\medskip

If there are vanishing indices in the bottom, then the bottom-left
corner belongs to the list: $A,A\cup C,A\cup \mathcal{A}_1 \sqcup \mathcal{A}_1$ or
$\mathcal{A}_1 \sqcup \mathcal{A}_1$.

\begin{itemize}
\item $A$: the permutation $\gp$ is then obviously not strongly irreducible.
\item $A\cup C$: the bottom-right corner is necessary $D$ or $D\cup
C$. If the top-right corner where $D$, then $\gp$ would be
not strongly irreducible. In particular, that means $\mathcal{A}_0$ is not empty. Hence
there exists $j_d < j_1\leq l$ such that $d=\gp(j_d)\in D$ and
$b=\gp(j_1)\in \mathcal{A}_0$. Then we choose any index $a\in A$ and any
index $c\in C$, and set: 
\begin{eqnarray*}
\tau^{\prime}_a=\tau_a+\varepsilon\\
\tau^{\prime}_d=\tau_d+\varepsilon\\
\tau^{\prime}_b=\tau_b-\varepsilon\\
\tau^{\prime}_c=\tau_c-\varepsilon\\
\forall i\notin\{a,b,c,d\}\quad \tau^{\prime}_i=\tau_i.
\end{eqnarray*}
Then $\tau^{\prime}$ is better than $\tau$ for $\varepsilon>0$ small enough,
which contradicts its maximality. 

\item $A\cup \mathcal{A}_1 \sqcup \mathcal{A}_1$, or $\mathcal{A}_1 \sqcup \mathcal{A}_1$. The bottom-right corner is
  necessary  $D$. If $\mathcal{A}_0$ is empty, then the top-right corner is
  also  $D$, and therefore $\gp$ is not strongly irreducible. If $\mathcal{A}_0$ is not
  empty, then we choose $a\in A$, $b\in \mathcal{A}_0$ and $c\in \mathcal{A}_1$,
  and set: 
\begin{eqnarray*}
\tau^{\prime}_a=\tau_a+2\varepsilon\\
\tau^{\prime}_b=\tau_b-\varepsilon\\
\tau^{\prime}_c=\tau_c-\varepsilon\\
\forall i\notin\{a,b,c\} \quad \tau^{\prime}_i=\tau_i
\end{eqnarray*}
And $\tau^{\prime}$ is better than $\tau$.
\end{itemize}

The proposition is proved.
\end{proof}

We now have all necessary tools for proving our main result.

\begin{proof}[Proof of Theorem~\ref{CNS}] We only have to prove the sufficient condition.
We consider a pseudo-suspension $\tau$ better than $\tau_{MV}$ and
maximal for this property. We can assume that $\gp:\{1,\ldots,l+m\}
\rightarrow \mathcal{A}$ is not strongly irreducible (i.e. at least one corner
is non empty in the decomposition) otherwise the theorem
follows from Lemma~\ref{lm:strongly:irred} and Proposition~\ref{max:donc:reducible:star}. Let us consider
a decomposition of $\gp$ as
$$
\left(\begin{array}{c|c|c}A\cup B & U & D\cup B \\\hline A\cup C& V & D \cup C 
\end{array}\right).
$$
where $A\cup B\cup C \cup D$ is maximal. 
Note that  $\gp^{\prime}= \bigr( \begin{smallmatrix} U \\ V
\end{smallmatrix} \bigl)$ defines a generalized permutation which is
not strongly irreducible by assumption. Note also that $\gp'$ does not
necessary satisfy Convention~\ref{conv:normal}, even if $\gp$ satisfies that
convention. \medskip

We define $\mathcal A^{\prime}=\mathcal A \backslash \bigl( A\cup B \cup
C\cup D\bigr)$; from Proposition~\ref{max:donc:reducible:star}, the restriction of $\tau$ to $\mathcal A^{\prime}$ is strict for $\gp'$.

Since $\gp$ is irreducible, there is one or two empty corners in
the decomposition. 
\begin{itemize}
\item If only one corner is empty, then it is on the right. So we can
  assume that $\gp$ decomposes as: 
\[
\left(\begin{array}{c|c|c}A\cup B & U & B \\\hline A & V & \\ 
\end{array}\right)\]
with $\gp^{\prime}= \bigr( \begin{smallmatrix}U\\V\end{smallmatrix} \bigl) $ strongly irreducible.

Let $i_1$ be the first vanishing index in the top line of $\gp$ and $i_2$
be the first vanishing index of the second line. Consider $i_b$ the
first index such that $b=\gp(i_b)\in B$. Then $i_b\leq i_1$ otherwise
there would be a subdecomposition of $\gp$ as  
\[ 
\left(\begin{array}{c|c}A^{\prime} & ***\\\hline A^{\prime} & ***
\end{array}\right)\]
and $\gp$ would be reducible. Now let $a=\gp(l+1)\in A$ and let $c\in \mathcal{A}_1$. We set:
\begin{eqnarray*}
\tau_b^{\prime}=\tau_b+2\varepsilon\\
\tau_c^{\prime}=\tau_c+2\varepsilon\\
\tau_a^{\prime}=\tau_a-\varepsilon\\
\forall j\notin \{a,b,c\} \quad \tau_j^{\prime}=\tau_j
\end{eqnarray*}
If $\varepsilon$ is small enough, then $\tau^{\prime}$ satisfies:
\begin{itemize}
\item For all $k\in\{1,\ldots,l\}\quad\sum_{i\leq k} \tau^{\prime}_{\gp(i)}> 0$ .
\item For all $k\in\{1,\ldots,m-1\}\quad \sum_{l<i\leq l+k} \tau^{\prime}_{\gp(i)}< 0$ .
\end{itemize}
And then, we can deduce from $\tau^{\prime}$ a suspension over $\gp$.

\item If two corner are empty, then we can  assume that $\gp$ decomposes as:
\[
\left(\begin{array}{c|c|c} B & U & B \\\hline  & V & \\ 
\end{array}\right)\]
\end{itemize}
with $\gp^{\prime}= \bigr( \begin{smallmatrix} U\\V\end{smallmatrix} \bigl) $ irreducible.
Now we choose $b\in B$ and $c\in \mathcal{A}_1$, and then set:
\begin{eqnarray*}
\tau_b^{\prime}=\tau_b+2\varepsilon\\
\tau_c^{\prime}=\tau_c+2\varepsilon\\
\forall j\notin \{b,c\} \quad \tau_j^{\prime}=\tau_j.
\end{eqnarray*}
Then $\tau^{\prime}$ defines a suspension over $\gp$ for $\varepsilon$
small enough. The theorem is proven.
\end{proof}

\section{Irrationality of linear involutions}\label{sec:dyn}

For an interval exchange map $T=(\pi,\lambda)$ either the underlying
permutation is reducible and then the transformation is never minimal
or $\pi$ is irreducible and $T$ has the Keane's property (and hence is minimal) for
almost every $\lambda$ (see section~\ref{sec:iem}). 
Furthermore $T$ admits a suspension if and only if $\pi$ is irreducible. \\
Hence the combinatorial set for which the dynamics of $T$ is good 
coincides with the one for which the geometry is good. As we will see,
the situation is more complicated in the general case. In this section
we prove Theorem~$B$ and the second half of Theorem~$A$.

\subsection{Keane's property}

\begin{defs}
A linear involution
has a connection (of length $r$) if there exist $(x,\varepsilon)\in
X\times \{0,1\}$ and $r\geq 0$ such that
\begin{itemize}
\item $(x,\varepsilon)$ is a singularity for $T^{-1}$.
\item $T^r(x,\varepsilon)$ is a singularity for $T$.
\end{itemize}
A linear involution with no connection is said to have
the Keane's property.
\end{defs}

Note that, by definition of a singularity, if we have a connection of 
length $r$ starting from $(x,\varepsilon)$, then $\forall r' < r,\ T^{r'}(x,\varepsilon)$ 
is not a singularity for $T$.

We first prove the following proposition:

\begin{prop} 
\label{no:connection}
Let $T$ be a linear involution. The
following statements are equivalent. 
\begin{enumerate}
\item $T$ satisfies the Keane's property.
\item $\mathcal R^n(T)$ is well defined for any $n\geq 0$ and the lengths 
of the intervals $\lambda^{(n)}$ tends to $0$ as $n$   tends to infinity.
\end{enumerate}
Moreover in the above situation the transformation $T$ is minimal.
\end{prop}

\begin{proof}[Proof of Proposition~\ref{no:connection}]
We denote by $\lambda^{(n)}$ the length parameters of the map
$\mathcal{R}^{(n)}(T)$,  by $\gp^{(n)}$, $\sigma^{(n)}$, $(l^{(n)},
m^{(n)})$ the combinatorial data, and by $X^{(n)}$ the subinterval of
$X$ corresponding to $\mathcal{R}^{(n)}(T)$. Let us assume that $T$ has no connection.
Then all the iterates of $T$ by the Rauzy-Veech induction are well defined.  Indeed 
it is easy to see that $T$ has the Keane's property if and only if its image
$\mathcal{R}(T)$ by the Rauzy-Veech induction  is well defined and has
the Keane's property. Hence if $T$ has the Keane's property, then by
induction,  all its iterates by $\mathcal{R}$ are well defined and
have the Keane's property. \medskip

Now we have to prove that $\lambda^{(n)}$ 
goes to zero as $n$ tends to infinity. Let $\mathcal{A}^{\prime}$ be the subset of elements $\alpha\in
\mathcal{A}$ such that $(\lambda_{\alpha}^{(n)})_n$ decreases an
infinite number of time in the sequence $\{\mathcal{R}^n(T))\}_{n}$, 
and let $\mathcal{A}^{\prime\prime}$ be its complement.

Repeating the arguments for the Proposition and Corollary 1 and 2 of
section 4.3 in~\cite{Yoccoz}, we have that:

\begin{itemize}
\item For $n$ large enough, the permutation $\gp^{(n)}$ can be written as:
$$
\left(\begin{array}{ccc|c}\alpha_1&\ldots&\alpha_{i_0}& *** \\\hline
  \beta_1&\ldots &\beta_{j_0} & *** \\  
\end{array}\right),
$$
with $\{\alpha_1,\ldots,\beta_{j_0}\}=\mathcal A'' \sqcup \mathcal A''$
\item For all $\alpha\in\mathcal{A}^{\prime}$,
  $\lambda^{(n)}_{\alpha}$ tends to zero.  
\end{itemize}

If $\mathcal{A}^{\prime}=\mathcal{A}$, then the proposition is proven. So we
can assume that $\mathcal{A}^{\prime}$ is a strict subset of
$\mathcal{A}$.  
Note that $\mathcal{A}^{\prime}$ cannot be empty. Therefore, we must have
\[\sum_{i=1}^{i_0}\lambda_{\alpha_i}=\sum_{j=1}^{j_0}\lambda_{\beta_j},\]
for some $1\leq i_0 \leq l^{(n)}-1$ and $1\leq j_0 \leq
m^{(n)}-1$. This means that $\mathcal{R}^n(T)$ has a connection of
length zero, hence $T$ has a connection. This contradicts the
hypothesis. So we have proven that if $T$ has no connections, then the
sequence $\{\mathcal{R}^n(T))\}_{n}$ of iterates of $T$ by the Rauzy-Veech
induction is infinite and  all length parameters of
$\mathcal{R}^{n}(T)$ tend to zero when $n$ tends to infinity. \medskip

Now we assume that $T$ has a connection. So, there exists  $u_0=(x,\varepsilon)$ in 
$X\times\{0,1\}$ which is  a singularity of $T^{-1}$, and such that its sequence 
$u_1,\dots,u_m$ of iterates  by $T$ is finite, with $u_m$ a
singularity of $T$. We denote by
$\overline{u_1},\dots,\overline{u_m}$ the projections of
$u_0,\dots,u_m$ on $X$. 
Let $u_{min}$ be the element of $\{u_0,\dots,u_m\}$ whose
corresponding projection to $X$ is minimal. We have $u_{min}>0$. If for all $n\geq 0$, 
the map $\mathcal{R}^n(T)$ is well defined and $u_{min}\in X^{(n)}$, then $X^{(n)}$ does 
not tend to zero, and hence there exists $\alpha\in \mathcal{A}$ such that $\lambda_\alpha^{(n)}$ 
does not tend to zero.
Hence we can assume that there exists a maximal $n_0$ such that
$\mathcal{R}^{n_0}(T)$ is well defined, and $X^{(n_0)}$ contains
$\overline{u_{min}}$. We want to show that $\mathcal{R}^{n_0+1}(T)$ is
not defined.  \medskip

Assume that $\mathcal{R}^{n_0+1}(T)$ is defined, then
$\overline{u_{min}}\notin X^{(n_0+1)}$.  
Since $\mathcal{R}^{n_0}(T)$ is an acceleration of $T$, there must
exists an iterate of $u_{min}$ by $T$, say $u_{k}$ which is a
singularity for  $\mathcal{R}^{n_0}(T)$. Either  $\overline{u_{k}}$ is
in $X^{(n_0+1)}$, either it is its right end. However,  $X^{(n_0+1)}$
does not contain $\overline{u_{min}}$, and $\overline{u_{min}}\leq
\overline{u_{k}}$. Therefore, we must have $u_{min}=u_k$, so
$u_{min}$ is a singularity for $\mathcal{R}^{n_0}(T)$. \smallskip

We prove in the same way that $u_{min}$ is also a singularity for
$\mathcal{R}^{n_0}(T)^{-1}$. This implies that we are precisely in the
case when the Rauzy-Veech induction is not defined. Hence we have proven
that if $T$ has a connection, then either the sequence
$(\mathcal{R}^n(T))_n$ is finite, either the length parameters do not
all tend to zero. \medskip

This proves the first part of the proposition. Now let $T$ be a
linear involution on $X$ satisfying the Keane's
property. Recall that $T$ is defined on the set $X\times \{0,1\}$. Let
us consider the first return map $T_0$ on $X\times \{0\}$. By
definition of $T_0$, one has for each $(x,\varepsilon)$ some return
time $k=k(x,\varepsilon)>0$ such that
$T_0(x,\varepsilon)=T^k(x,\varepsilon)$. But $T$ 
is piecewise linear thus for any $(y,\varepsilon)$ in a small
neighborhood of $(x,\varepsilon)$, the return time
$k(y,\varepsilon)=k(x,\varepsilon)$. Since the derivative of $T$ is
$1$ if $(x,\varepsilon)$ and $T(x,\varepsilon)$ belong to the same
connected component and $-1$ otherwise, the derivative of
$T_0$ is necessary $1$. Hence $T_0$ is an \emph{interval exchange}
map. Obviously $T_0$ has no connexion since it is an acceleration of
$T$, hence $T_0$ is minimal. Similarly, the first return map $T_1$ of
$T$ on $X\times \{1\}$ is also minimal. Since $T$ satisfies 
Convention~\ref{conv:normal}, any orbit of $T$ is
dense on $X\times \{0\}$ and $X\times \{1\}$ therefore $T$ is
minimal. The proof is complete.
\end{proof}

\subsection{Dynamical irreducibility}
\label{sec:red12}

\begin{defs}
\label{def:dyn:irr}
Let $T=(\gp,\lambda)$ be a linear involution. We will
say that $\lambda$ is \emph{admissible} for $\pi$ (or $T$ has
admissible parameters) if none of the following assertions holds:
\begin{enumerate}

\item
\label{def:1}
$\gp$ decomposes as $\left( \begin{smallmatrix} A|& *** \\\hline A| & *** \\ 
\end{smallmatrix}\right),\left(\begin{smallmatrix} ***&| D \\\hline  *** & |D \\ 
\end{smallmatrix}\right) \textrm{ or } 
\left(\begin{smallmatrix} A\cup B|& D\cup B \\\hline A \cup C|& D\cup C \\ 
\end{smallmatrix}\right)$ \\
with $A,D \subset \mathcal{A}_{01} \textrm{ and } B=\mathcal{A}_0,\ C=\mathcal{A}_1$
and $A,D$ non empty in the two first cases.

\item
\label{def:2}
There is a decomposition of $\gp$ as $\left(\begin{smallmatrix}A\cup B|& &*** & &|B \cup D\\\hline 
A\cup C| & \alpha_0& *** &\alpha_0 & |C\cup D \\ \end{smallmatrix}\right)$, 
with (up to switching the top and the bottom of $\gp$) $A,D \subset
\mathcal{A}_{01} \textrm{ and } \emptyset \not = B\subset \mathcal{A}_0,\
C\subset \mathcal{A}_1$  
and the length parameters $\lambda$ satisfy the following inequality 
$$
\sum_{\alpha \in C} \lambda_\alpha \leq \sum_{\alpha \in B} \lambda_\alpha
\leq \lambda_{\alpha_0} + \sum_{\alpha \in C} \lambda_\alpha.
$$
\end{enumerate}
A generalized permutation $\gp$ will be called \emph{dynamically
  irreducible} if the corresponding  set of admissible parameters is
nonempty. 
\end{defs}

The set of  admissible parameters of a generalized permutation is
always open. 

\begin{rem}
\label{rm:th:un}
These two combinatorial notions of reducibility were introduced by the second author 
(see~\cite{La1}). Observe that if $\lambda$ is not admissible for
$\gp$, then  $T=(\pi,\lambda)$ have a connection of length $0$ or $1$
depending on cases~\eqref{def:1} or~\eqref{def:2} of
Definition~\ref{def:dyn:irr}, and is never minimal. More precisely
there exist two invariant sets of positive measure.\\
One can also note that if $\pi$ is irreducible then $\pi$ is
dynamically irreducible (the set of admissible parameters being the
entire parameters space).
\end{rem}

The length parameters for $T$ cannot be linearly
independent over $\mathbb{Q}$ since they must satisfy a nontrivial relation with
integer coefficients. A linear involution 
$T=(\lambda,\gp)$ is said to have \emph{irrational parameters}
if  $\{\lambda_{\alpha}\}$ generates a
$\mathbb{Q}$-vector space of dimension $\#\mathcal{A}-1$. Almost all
linear involutions have irrational parameters, and this property 
is preserved by the Rauzy-Veech induction.

\begin{proof}[Proof of Theorem~$B$]
If $\pi$ is dynamical reducible, the non minimality comes from
Remark~\ref{rm:th:un}. Conversely let us assume that $\pi$ is a
dynamical irreducible permutation and let $T=(\pi,\lambda)$ be a
linear involution with irrational parameters and
$\lambda$ admissible for $\pi$.

We still denote by $\lambda^{(n)}$ the length parameters of
$\mathcal{R}^{(n)}$ and by $\gp^{(n)}$, $\sigma^{(n)}$, $(l^{(n)},
m^{(n)})$ the combinatorial data. 

The proof has two steps: first we show using Proposition~\ref{no:connection} 
that if $T$ does not have the Keane's property, then
there exists $n_0$ such that $\mathcal{R}^{n_0}(T)$ does not have admissible parameter (case~\eqref{def:1} of Definition~\ref{def:dyn:irr}). Then we show that in this case $\lambda$ is not admissible for $\gp$. This will imply the theorem.

\emph{First step:}
We assume that the sequence is finite. Then there exists
$\mathcal{R}^{n_0}(T)$ that admits no Rauzy-Veech induction.  
Since $\lambda^{(n_0)}$ is irrational then either
$\sigma^{(n_0)}(l^{(n_0)})=l^{(n_0)}+m^{(n_0)}$, or $l^{(n_0)}$
belongs to the only pair $\{i,\sigma^{(n_0)}(i)\}$ on the top of
the permutation and $l^{(n_0)}+m^{(n_0)}$ belongs to the only pair
$\{j,\sigma^{(n_0)}(j)\}$ on the bottom of the permutation. In
each case, $\mathcal{R}^{n_0}(T)$ does not have admissible parameter (case~\eqref{def:1}). 

Now we assume that the lengths parameters do not all tend to zero. 
As in the proof of Proposition~\ref{no:connection},  for $n$ large
enough, the generalized permutation $\gp^{(n)}$ decomposes as: 
\[
\left(\begin{array}{ccc|c}a_1&\ldots&a_{i_0}& *** \\\hline b_1&\ldots &b_{j_0} & *** \\ 
\end{array}\right), \]
with $\{a_1,\ldots,b_{j_0}\}=\mathcal A'' \sqcup \mathcal A''$, for some  $\emptyset \neq 
\mathcal{A}^{\prime\prime}\subset \mathcal{A}$ and some  $1\leq i_0<
l^{(n)}$ and  $1\leq j_0< m^{(n)}$. Recall that
$$
\sum_{i=1}^{i_0}\lambda_{\gp^{(n)}(i)}=\sum_{j=1}^{j_0}\lambda_{\gp^{(n)}(j)}.
$$
The map $\mathcal{R}^n(T)$ has irrational parameters, therefore
$\gp^{(n)}$ must decompose as: 
$$
\left(\begin{array}{c|c}A& *** \\\hline A & *** \\ 
\end{array}\right), \textrm{ or }
\left(\begin{array}{c|c}*** &D\\\hline *** &D \\ 
\end{array}\right),
$$
so $\mathcal{R}^{n}(T)$ does not have admissible parameter
(case~\eqref{def:1}). 
\medskip

\emph{Second step:}
It is enough to prove that if $T^{\prime}=\mathcal{R}(T)$ does not
have admissible parameter, then so is $T$.  
We can assume without loss of generality that the combinatorial
Rauzy-Veech transformation is  $\mathcal{R}_0$. 
We denote by $\gp,\sigma,\lambda$ the data of $T$ and by
$\gp^{\prime},\sigma^{\prime}, \lambda^{\prime}$ the data of
$T^{\prime}$. If $\gp^{\prime}$ decomposes as: 
\[
\left(\begin{array}{c|c} ***& D \\\hline *** & D \\ 
\end{array}\right),
\]
let us consider $l^{\prime}$ the last element of the top line. Its twin
$\sigma^{\prime}(l^{\prime})$  is on the bottom-right corner, but is
not $l^{\prime}+m^{\prime}$. 
We denote by
$\beta=\gp^{\prime}(\sigma^{\prime}(l^{\prime})+1)$. Then it is clear
that we obtain $\gp$ by removing $\beta$ from that place and
putting it at the right-end of the bottom line. Then $T$ does not have admissible parameter (case~\eqref{def:1}). 

Now we assume that $\gp^{\prime}$ decomposes as:
\[
\left(\begin{array}{c|c} A& *** \\\hline A & *** \\ 
\end{array}\right).
\]
If $\sigma^{\prime}(l^{\prime})$ is on the bottom line, the situation
is analogous to the previous case. If not, then we denote by
$\beta=\gp^{\prime}(\sigma^{\prime}(l^{\prime})-1)$ and
$\alpha=\gp^{\prime}(l^{\prime})$, and we get $\gp$ by removing
$\beta$ from $\sigma^{\prime}(l^{\prime})-1$ and putting it on the
right-end of the bottom line. If this place is in the top-right
corner, then clearly,  $T$ does not have admissible parameter (case~\eqref{def:1}). However, it might be the
last element of the top-left corner. In that case, setting
$A=A^{\prime}\cup \{\beta\}$, the generalized permutation $\gp$
decomposes as: 
\[
\left(\begin{array}{cc|c|cc} 
A^{\prime} & \alpha & *** & \alpha & \\ \hline 
A^{\prime} \cup \{\beta\}& & ***& &\beta \\ 
\end{array}\right),
\]
with $\lambda_\beta=\lambda_\beta'>0$ and
$\lambda_\alpha=\lambda_\alpha'+\lambda_\beta'>\lambda_\beta$, 
hence  $T$ does not have admissible parameter (case~\eqref{def:2}).

Now we assume that $\gp^{\prime}$  decomposes as
\[
\left(\begin{array}{c|c} A\cup B& B\cup D \\\hline A\cup C & C\cup D \\ 
\end{array}\right).
\]
Then we obtain $\gp$ from $\gp^{\prime}$ by removing an element on the
top-left corner or on the bottom-right corner, and putting it at the
right-end of the bottom line. Then  $T$ does not have admissible
parameter (case~\eqref{def:1}). The other cases are similar.
\end{proof}

\section{Dynamics of the renormalized Rauzy-Veech induction}
\label{dyn:rauzy:veech}

As we have seen previously, there are two notions of irreducibility
for a linear involution.
\begin{itemize}
\item ``Geometrical irreducibility'' as stated in section
  \ref{combinatoric}, that we just called irreducibility. 
\item Dynamical irreducibility as stated in section~\ref{sec:dyn}.
\end{itemize}

In this section, we first prove that the set of irreducible
linear involutions in an attractor for the
renormalized Rauzy-Veech induction. Then we show that (analogously to the case
of interval exchange transformations) the renormalized Rauzy-Veech
induction is recurrent for almost all irreducible linear involutions.

\subsection{An attraction domain}\label{attract_dom}

\begin{proof}[Proof of the first part of Theorem~$C$]
We can find a non-zero pseudo-suspension $(\tau_\alpha)_{\alpha\in
\mathcal{A}}$ (see Definition \ref{def:pseudo}) otherwise it is easy
to show that $T$ does not have admissible parameter (case (1)). For
all $\alpha$, we denote by $\zeta_\alpha$ the complex number  
$\zeta_\alpha=\lambda_\alpha+i\tau_\alpha$. Then, as in section
\ref{sec:red12}, we consider a broken line $L_0$ which starts at $0$, and
whose edge number~$i$ is represented by the complex number
$\zeta_{\gp(i)}$, for $1\leq i\leq l$. Then we consider  a broken line
$L_1$, which starts on the same point as $L_0$, and whose edge
number~$j$ is represented by the complex number $\zeta_{\gp(l+j)}$ for
$1\leq j\leq m$.

\begin{figure}[htbp]
\psfrag{1}{$\scriptstyle 1$} \psfrag{2}{$\scriptstyle 2$}
\psfrag{3}{$\scriptstyle 3$} \psfrag{4}{$\scriptstyle 4$}
\psfrag{5}{$\scriptstyle 5$} \psfrag{6}{$\scriptstyle 6$}
\psfrag{7}{$\scriptstyle 7$} 
\psfrag{a}{$\scriptstyle X$} \psfrag{b}{$\scriptstyle X^{(n)}$}

   \begin{center}
     \includegraphics{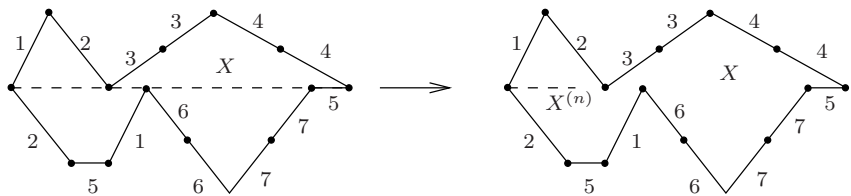}
    \caption{The transformation $T$ is the first return map of 
the vertical foliation on a union of saddle connections.}
     \label{special:case}
   \end{center}
\end{figure}

\emph{Special case:} We assume that  $L_0$ and $L_1$ only intersect on
their endpoints. Then they define a flat surface $S$, and $T$ appears
as a first return map of the vertical foliation on a segment $X$ which is a
union of horizontal saddle connections  (see Figure
\ref{special:case}).  
After $n$ steps of the Rauzy-Veech induction, the resulting
linear involution $\mathcal{R}^n(T)$ is the
first return map of the vertical flow of $S$ on a shorter segment
$X^{(n)}$, which is adjacent to the same singularity as $X$. Since~$T$
has no connection, then the length of $X^{(n)}$ tends to zero when $n$
tends to infinity by the first part of
Proposition~\ref{no:connection}. Hence for $n$ large enough,
$\mathcal{R}^n(T)$ is the first return map of the vertical flow of $S$ on a segment,
adjacent to a singularity, and with no singularities in its interior. 
With our construction of $S$, it is clear that any vertical saddle connection would
intersect $X$ and would give a connexion on $S$. Since $T$ has no
connection, the surface $S$ has no vertical saddle connection (note
that this is not true in general for a first return map on a
transverse segment). According to Proposition~\ref{prop:irr:natural},
$(\pi^{(n)},\lambda^{(n)})$ admits a suspension and hence
Theorem~\ref{CNS} implies that $\pi^{(n)}$ is irreducible. 
The theorem is proven for that case.

\begin{figure}[htbp]
\psfrag{1}{$\scriptstyle 1$} \psfrag{2}{$\scriptstyle 2$}
\psfrag{3}{$\scriptstyle 3$} \psfrag{4}{$\scriptstyle 4$}
\psfrag{5}{$\scriptstyle 5$} \psfrag{6}{$\scriptstyle 6$}
\psfrag{7}{$\scriptstyle 7$} 
\psfrag{a}{$\scriptstyle a$} \psfrag{b}{$\scriptstyle b$}
\psfrag{xe}{$X_\varepsilon$} \psfrag{se}{$S_\varepsilon$}

   \begin{center}
     \includegraphics{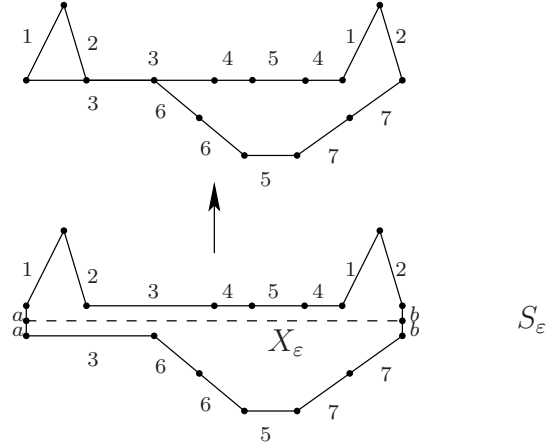}
    \caption{Constructing $T$ as a first return map on a regular segment
  of a surface $S_\varepsilon$.}
     \label{general:case}
   \end{center}
\end{figure}

\emph{General case:}
The two broken lines $L_0$ and $L_1$ might have other intersection
points. We first show this still defines a flat surface. We
consider the line $L^{\varepsilon}_0$ that starts at the complex number
$2 i \varepsilon$. Then we join the first points of
$L_0^\varepsilon$ and $L_1$ by a vertical segment, and do the same for
their ends points (see Figure \ref{general:case}). This defines a
polygon and the non vertical sides come by pairs, so we can glue them
as previously. There  are two vertical segments left.  
We decompose each vertical segment into a pair of vertical segments of
the same length and glue them together (see the figure). This creates a
pole for each initial segment. We denote by $S_\varepsilon$ the
resulting flat surface.  
The first return map of the vertical flow on the horizontal segment
$X_{\varepsilon}$ joining the two poles is $T$. The surface
$S_\varepsilon$ has two vertical saddle connections of length
$\varepsilon$ starting from the poles, but there is no other vertical
saddle connection on $S_\varepsilon$ since $T$ has no connections. 
When $\varepsilon$ tends to zero, the two vertical saddle connections
are the only ones that shrink to zero. Hence there is no loop that
shrink to zero. Furthermore, the initial pseudo-suspension is nonzero, so the area of
$S_\varepsilon$ is bounded from below. Hence, the surface $S_\varepsilon$ does not
degenerate when $\varepsilon$ tends to zero and so there exists a
sequence $(\varepsilon_k)$ that tends to zero such that
$(S_{\varepsilon_k})$ tends to a surface $S$.
 
The segment $X\subset S$ corresponding to the limit of $X_{\varepsilon_k}$, as
$k$ tends to infinity, might be very complicated and the first
return map on $X$ is not well defined. 

The transformation  $\mathcal{R}^n(T)$ is the first return map of the
vertical flow of $S_{\varepsilon_k}$ on a shortest horizontal segment
$X_{\varepsilon}^{(n)}$, adjacent to one of the poles. If $n$ is large
enough, then the segment $X^{(n)}\subset S$ corresponding to the limit of 
$X_{\varepsilon_k}^{(n)}$ has no singularity on its interior. Since the
surgery corresponding to contracting $\varepsilon$ does not change the
the vertical foliation, the first return map of the vertical foliation of $S$ on
$X^{(n)}$ is precisely $\mathcal{R}^n(T)$.  

As in the special case, the surface $S$ does not have any vertical
saddle connection, so the generalized permutation corresponding to
$\mathcal{R}^n(T)$ is irreducible and the proposition is proven. 

\end{proof}

\subsection{Recurrence}
The following lemma is analogous to Proposition~9.1 in~\cite{Ve82}. 

\begin{lem}
\label{lem:sing}
Let $T$ be a linear involution on
$X=(0,L)$ with no connection and let $(x,\varepsilon)\in X\times\{0,1\}$ be a
singularity for $T$. Let $X^{(n)}\subset X$ be the subinterval
corresponding to the linear involution
$\mathcal{R}^n(T)$. There exists $n>0$ such that $X^{(n)}=(0,x)$. 
\end{lem}

\begin{proof}
 Since $T$ has no connection, there exists a first $n>0$ such that
 $x\notin X^{(n)}$. So $x\in X^{(n-1)}$, and $(x,\varepsilon)$ is
 still a singularity for $\mathcal{R}^{n-1}(T)$. We obtain
 $\mathcal{R}^{n}(T)$ from  $\mathcal{R}^{n-1}(T)$ by considering the
 first return map on the largest subinterval $X^{(n)}\subset
 X^{(n-1)}$ whose right endpoint corresponds to a singularity of
 $\mathcal{R}^{n-1}(T)$. So $X^{(n)}=(0,x)$. 
\end{proof}

Let $\pi_0$ be an irreducible generalized permutation, and let $C$ be
the set of generalized permutations that can be obtained by iterations
of the maps $\mathcal{R}_0$ and $\mathcal{R}_1$ (when possible).

We define $\mathcal{T}_C=\{(\pi,\zeta),\ \pi\in
C,\zeta \textrm{ is a suspension data for $\pi$} \}$. We have defined the
Rauzy-Veech map on the space $\mathcal{T}_C$. It
defines an almost everywhere invertible map: If $\sum_{i=1}^l Im(\zeta_{\pi(i)})
\not = 0$ then $(\pi,\zeta)$ has exactly one preimage for
$\mathcal R$.

We define the quotient $\mathcal{Q}_C$ of $\mathcal{T}_C$ by the
equivalence relation generated by $(\pi,\zeta)\sim
\mathcal{R}(\pi,\zeta)$. 

One will denote by $m$ the natural Lebesgue measure  on
$\mathcal{T}_C$ i.e. $m=d\gp d\zeta$, where $d\zeta$ is the natural
Lebesgue measure on the hyperplane $\sum_{i=1}^l
\zeta_{\gp(i)}=\sum_{j=l+1}^{2d} \zeta_{\gp(j)}$, and $d\gp$ is the
counting measure. The mapping $\mathcal{R}$ preserves $m$, so it
induces a measure, again denoted by $m$ on $\mathcal{Q}_C$.
 
The matrix $g_t$ acts on $\mathcal{T}_C$ by
$g_t(\pi,\zeta)=(\pi,(g_t(\zeta_\alpha))_\alpha)$, where $g_t$ acts on
$\zeta_\alpha\in \mathbb{C}=\mathbb{R}^2$ linearly. This action
preserves the measure $m$ on $\mathcal{T}_C$ and commutes with
$\mathcal{R}$, so it descends to a measure preserving flow on
$\mathcal{Q}_C $ called the Teichm\"uller flow.

If $(\pi,\zeta)$ is a suspension data, we denote by
$|Re(\zeta)|_{\pi}$ the length of the corresponding interval,
\emph{i.e.} $\sum_{i=1}^l Re(\zeta_{\pi(i)})$. The subset
$$
\left\{(\pi,\zeta)\in \mathcal{T}_C;\ 1 \leq |Re(\zeta)|_{\pi} \leq
1+\min\bigl(Re(\zeta_{\gp(l)}),Re(\zeta_{\gp(2d)})\bigr) \right\}
$$
is a fundamental domain of $\mathcal{T}_C$ for the relation $\sim$ and
the first return map of the  Teichm\"uller flow on
$$
\mathcal{S}=\{(\pi,\zeta);\ \gp\in C,\ |Re(\zeta)|_{\pi}=1\} / \sim
$$
is the renormalized Rauzy-Veech induction on suspensions.

\begin{prop}
The zippered rectangle construction  provides a finite covering $Z$
from $\mathcal{Q}_C$ to a subset of full measure in a connected
component of a stratum $\mathcal{Q}(k_1,\dots,k_n)$ of the 
moduli space of quadratic differentials. The degree of this cover is
$h!$ where $h$ is the dimension of the stratum. Moreover $h=2g+n-2$ 
($g$ is the genus of the surfaces).
\end{prop}

\begin{proof}
Let $S$ be a (generic) flat surface in $\mathcal{Q}(k_1,\dots,k_n)$
with no vertical and no horizontal saddle connection. Consider a
horizontal separatrix $l$ adjacent to a given singularity $P$. We call
\emph{admissible} a segment $X$ adjacent to $P$, such that the
vertical geodesic passing through the right endpoint of $X$  meets a
singularity before returning to $X$, in positive or negative
direction. Then Proposition~\ref{prop:irr:natural} implies that there
exists a corresponding suspension datum $\zeta$ such that
$S=Z(\pi,\zeta)$. Conversely, any $\zeta$ such that $S=Z(\pi,\zeta)$ is
obtained by this construction.

Now let $X_0,X_1$  be two admissible segments, and let
$\zeta_0,\zeta_1$  be the corresponding suspension data. One can
assume without loss of generality that $X_0\subset X_1$ and their left
endpoint is the singularity $P$.  Let $T_0,T_1$ be the linear involutions
corresponding to $X_0,X_1$. The right endpoint
of $X_0$ corresponds to a singularity of $T_1$. Hence there exists
$n\geq 0$ such that $\mathcal{R}^n(T_1)=T_0$, and therefore
$\mathcal{R}^n(\zeta_1)=\zeta_0$.

So we have proven that for each separatrix $l$ adjacent to a
singularity, there is only one preimage of $S$ by the mapping
$Z$. So $Z$ is a finite covering. The degree of $Z$ is obvious by
construction. If $2h$ is the number of possible choices of horizontal
separatrices then the degree of $Z$ is $h!$ (choices of labels 
and the choice of the intervals $X\times \{0,1\}$).

For each singularity, one has $k_i+2$ separatrices. Thus
$$
2h=\sum_{i=1}^n (k_i + 2) = 4g-4+2n = 2(2g+n-2).
$$
The proposition is proven.
\end{proof}

\begin{proof}[Proof of the second part of Theorem~$C$]
The subset  $\mathcal{Q}_C^1$ corresponding to surfaces of area~1 is a
finite ramified cover of a connected component of a stratum of
quadratic differentials, and the corresponding  Lebesgue measures are
proportional.
     
By Theorem 0.2 in~\cite{Ve3} the volume of the moduli space of
quadratic differentials is finite, and so, $\mathcal{Q}_C^1$ has
finite measure. Hence the Teichm\"uller geodesic flow on $\mathcal{Q}_C$ is recurrent for
the Lebesgue measure. Recall that the Rauzy-Veech renormalization for
suspensions $\mathcal R_r$ is the cross section of the Teichm\"uller
geodesic flow on $\mathcal{S}$; therefore the Rauzy-Veech
renormalization for suspension is recurrent. \medskip

We have $d\zeta=d\lambda d\tau$, and the Rauzy-Veech induction
commutes with the projection $(\gp,\zeta)\mapsto (\gp,\lambda)$. So,
for almost all parameters $\lambda$, the  sequence
$(\mathcal{R}_r^n(\gp,\lambda))_n$ is recurrent.
\end{proof}

\begin{rem}
Note that the proof of theorem~$C$ does not use the fact that a
linear involution satisfying the Keane's property is
minimal. We can use this theorem to give an alternative proof of the
minimality of such map.
let $T$ be a linear involution with the Keane's
property. From Theorem~$C$, there exists $n \geq 0$ such
that $\mathcal R^n(T)=(\pi,\lambda)$ is the cross section of the
vertical foliation on a flat surface with no vertical saddle
connection. Any infinite vertical geodesic on such a surface is dense
(see e.g.~\cite{Masur:Tabachnikov}). Thus $\mathcal R^n(T)$ is minimal
and so is $T$.
\end{rem}

\section{Rauzy classes}\label{rauzy:classes:giem}

As we have seen previously, the irreducible generalized permutations
are an attractor for the Rauzy-Veech induction. In this section, we
prove that there is no smaller attractor. We also prove Theorem~$D$.

We first define the Rauzy classes and then the extended Rauzy
classes. \medskip

Given a permutation $\pi$, we can define {\it at most} two other
permutations $\mathcal R_\varepsilon(\pi)$ with $\varepsilon=0,1$ when
$\mathcal{R}_\varepsilon$ is well defined. The relation  $\pi\sim
\mathcal{R}_\varepsilon(\pi)$ generates a partial order on the set
of generalized permutations; we represent it as a directed graph $G$, and
as for permutations, we will call Rauzy classes the connected components of
this graph.

In the case of interval exchanges, the periodicity of the maps $\mathcal{R}_0$ 
and $\mathcal{R}_1$ gives an easy proof of the fact that the above relation is 
an equivalence relation (proposition of section~\ref{rauzy:classes:iem}). 
Here the argument fails because these maps are not 
always defined, and it may happen that $\mathcal{R}_0(\gp)$ is well defined, but not 
$\mathcal{R}_0^2(\gp)$. However, the corresponding statement is still true.

\begin{prop}
\label{prop:classes}
The above partial order is an equivalence relation on the set of
{\it irreducible} generalized permutations.
\end{prop}

\begin{proof}
Let $\gp$ and $\gp^{\prime}$ be two generalized permutations. Assume
that there  is a sequence of maps $\mathcal{R}_0$ and $\mathcal{R}_1$
that sends $\gp$ to $\gp^{\prime}$.  
If $\gp^{\prime}=\mathcal{R}_\varepsilon(\gp^{\prime\prime})$, then
for any parameters $\lambda^{\prime}$, there exist parameters
$\lambda^{\prime\prime}$ such that
$\mathcal{R}(\gp^{\prime\prime},\lambda^{\prime\prime})=(\gp^{\prime},
\lambda^{\prime})$.  Iterating this argument, there exists
$(\gp,\lambda^{0})$ and $n_0$ such that $\mathcal{R}^{n_0}
(\gp,\lambda^{0}) =(\gp^{\prime},\lambda^{\prime})$. But for any
$\lambda$ in a sufficiently small neighborhood $U$ of $\lambda^{0}$,  
the generalized permutation corresponding to
$\mathcal{R}^{n_0}(\gp,\lambda^{0})$ is $\pi^{\prime}$.  

Recall that renormalized Rauzy-Veech induction map is recurrent
(Theorem~$C$) thus one can find $\lambda \in U$ such that the sequence
$(\mathcal{R}_r^n(\gp,\lambda))_n$ come back in a neighborhood of
$(\gp,\lambda)$ infinitely many time. Furthermore,
$\mathcal{R}_r^{n_0}(\gp,\lambda)=(\gp',\lambda^{(n_0)})$.  Thus
$(\mathcal{R}_r^n(\gp,\lambda))_n$ gives a sequence of generalized
permutations that reach $\gp'$ and then reach $\gp$. So, it gives a
combination of the maps $\mathcal{R}_0$ and $\mathcal{R}_1$ that sends
$\gp^{\prime}$ to $\gp$. This proves the proposition.
\end{proof}

\begin{defs}
Let $2d=l+m$. We define the symmetric permutation $s$ of $\{1,\dots,2d\}$ by
$s(i)=2d+1-i, \ \forall i=1,\dots,2d$. If $\gp$ is a generalized permutation of
type $(l,m)$ defined over an alphabet $\mathcal{A}$ of $d$ letters, we
define the generalized permutation $s\gp$ to be of type $(m,l)$ by
$$
(s\gp)(k):=\gp\circ s (k).
$$
We start from an irreducible generalized permutation $\gp$ and 
we construct the subset of irreducible generalized permutation that can be obtained from $\gp$ 
by some composition of the maps $\mathcal R_0,\  \mathcal{R}_1$, and $s$. 
The quotient of this set by 
the equivalence relation generated by $\pi \sim f\circ \pi$ for any bijective map $f$ 
from $\mathcal A$ onto $\mathcal A$ is called the {\it extended Rauzy class} of $\gp$.
\end{defs}

\begin{rem}
The quotient by the equivalence relation generated by  $\pi \sim
f\circ \pi$ means that we look at generalized permutations defined up
to renumbering. This is needed for technical reasons in the proof of
Theorem~$D$.
\end{rem}

\begin{rem}
In opposite to the case of interval exchange maps, the definition of
irreducibility we gave in section~\ref{combinatoric} is not invariant
by the map $s$: for instance, the generalized permutation $\pi=\left(
\begin{smallmatrix} 1 & 2 & 1 \\ 2 & 3 & 3 & 4 & 4 \end{smallmatrix}
\right)$ is irreducible while $s\pi=\left(
\begin{smallmatrix} 4&4&3&3&2\\ 1&2&1\end{smallmatrix}
\right)$ is reducible. 

So an extended Rauzy class is obtained after considering the set of
generalized permutations obtained from $\gp$ by the extended Rauzy
operations, and  intersecting this set by irreducible generalized
permutations. The results from  the previous section shows that our
definition of irreducibility is the good one  with respect to the
Rauzy-Veech induction, but we see that the convention of the
``left-end singularity'' is a real choice.
\end{rem}

\begin{rem}
Let $T$ be a linear involution defined on an interval
$X=(0,L)$. Recall that Rauzy-Veech induction applied on $T$ consists in
considering the first return map on $(0,L^{\prime})$, where~$L^{\prime}$ is the maximal element of $(0,L)$ that
corresponds to a singularity of $T$. In terms of generalized
permutation, this corresponds to the $\mathcal{R}_\varepsilon$ mapping.  

One can consider the first return map of $T$ on the interval 
$(L^{\prime\prime},L)$, where~$L^{\prime\prime}$ is
the minimal element of $(0,L)$ that corresponds to a singularity
of $T$. In terms of generalized permutations, this corresponds to the
the conjugaison of $s\circ \mathcal{R}_\varepsilon \circ s$ map. We
will call this the ``Rauzy-Veech induction of $T$ by cutting on the
left of $X$'', while the usual Rauzy-Veech induction will on the
opposite called the ``Rauzy-Veech induction of $T$ by cutting on the
right of $X$''.
\end{rem}

\begin{proof}[Proof of Theorem~$D$]
Let $\gp_1$ be an irreducible generalized permutation. The
corresponding set of suspension data is connected (even convex), so the set of
surfaces constructed from a suspension data, using the zippered
rectangle construction, belongs to a connected component of the moduli space of
quadratic differentials. 

It is also open and invariant by the action
of the Teichm\"uller geodesic flow, hence it is a subset of full measure
by ergodicity. 

Let $\gp_2$ be a generalized permutation that corresponds to the same
connected component of the 
moduli space. Then there exists a surface $S$ and two segments $X_1$
and $X_2$, each one being adjacent to a singularity $x_1$ and $x_2$,
such that for each $i$, the linear involution $T_i$ given
by the first return maps on $X_i$ has combinatorial data $\gp_i$. We
can assume that $S$ has no vertical saddle connection. 

We recall that each $X_i$ has an orientation so that the
corresponding singularity $x_i$ is in its left end. Consider the
vertical separatrix $l$ starting  from $x_2$, in the positive
direction and let $y_1$ be its first intersection point  with $X_1\cup
\{x_1\}$. 

Applying the usual Rauzy-Veech
induction for $T_2$, the map  $\mathcal{R}^n(T_2)$ is a first return
map of the vertical flow on a subinterval $X_2^{(n)}\subset X_2$,
adjacent to $x_2$. If $n$ is large enough, then  $\mathcal{R}^n(T_2)$
is isomorphic to the first return map on the subinterval  
$(y_1,y_2)\subset X_1$,  of the same length as $X_2^{(n)}$. We assume
first that $y_1<y_2$, hence this first return map is  consistent with
the positive direction on $X_1$.   

We now have to apply Rauzy-Veech inductions (on the right and on the
left) on $T_1$ until we get a first return map on $(y_1,y_2)$ with
corresponding generalized permutation $\gp_3$. Since $\gp_3$ is by
construction, up to renumbering the alphabet,  in the same Rauzy
class as $\gp_2$, we will therefore find some composition of the maps
$\mathcal{R}_\varepsilon$, $s\circ \mathcal{R}_\varepsilon\circ s$
that send $\gp_1$ to $\gp_2$.

Note that $y_2$ might not correspond a priory to some singularities
of $T_1$, so naive Rauzy-Veech induction on $X_1$ might  miss the
interval $(y_1,y_2)$. But $(y_1,0)$ or $(y_1,1)$  is a singularity, so we can cut the
interval on the left until $y_1$ is the left end, this will eventually
occurs because of Lemma~\ref{lem:sing}. Then after cutting on the
right $y_2$ will become the right end of the corresponding interval.
\smallskip

If $y_2< y_1$, then similarly, by cutting on the right and then on the
left, we get two linear involutions corresponding to
first returns maps that only differ by a different choice of
orientation. Hence we have found  some composition of the maps
$\mathcal{R}_\varepsilon$, $s\circ \mathcal{R}_\varepsilon\circ s$
that send $\gp_1$ to some $\gp_3$, such that $s\gp_3$ is in the same
Rauzy class as $\gp_2$.  

Hence we have proved that if two irreducible generalized permutations
correspond to the same connected component, then they are in the same
extended Rauzy class. To prove the converse, we must consider a
slightly more general kind of suspensions that do not necessary
corresponds to a singularity on the left. The corresponding
``extended'' suspension data satisfy
\begin{enumerate}
\item $\forall \alpha \in \mathcal{A} \quad Re(\zeta_{\alpha})>0$.  
\item $\forall 1\leq i \leq l-1 \quad t+Im(\sum_{j\leq i} \zeta_{\gp(j)})>0$ 
\item $\forall 1\leq i \leq m-1 \quad t+Im(\sum_{1\leq j\leq i} \zeta_{\gp(l+j)})<0$
\item $\sum_{1\leq i\leq l} \zeta_{\gp(i)} = \sum_{1\leq j\leq m}\zeta_{\gp(l+j)}$
\end{enumerate}
for some $t\in \mathbb{R}$ (the case $t=0$ corresponds to suspension
data as seen previously). \medskip 

Then we can extend the zippered rectangle construction to these
extended suspension data. As in the usual case, the space of extended
suspension data corresponding to a generalized permutation is convex,
so the set of surfaces corresponding to a given generalized
permutation belong to a connected component of stratum. Then it is
easy to see that if $\pi'$ is obtained from $\pi$ by the map
$\mathcal{R}_0$, $\mathcal{R}_1$ or $s$, then the corresponding
connected component is the same.
\end{proof}

Historically, extended Rauzy classes have been used to prove the non
connectedness of some stratum of Abelian differentials (see for
instance \cite{Ve3}). For this case, some topological invariants were
found by Kontsevich and Zorich \cite{Kontsevich:Zorich}
(hyperellipticity and spin structure). For the case of quadratic
differentials, all non-connected components (except four special
cases) are distinguished by hyperellipticity \cite{La1}. For the four
``exceptional ones'', the only known proof up to now is an explicit
computation of the corresponding extended Rauzy classes. Theorem~$D$,
which is now formally proven complete the proof of the following

\begin{NNths}[Zorich]
The strata $\mathcal Q(-1,9)$, $\mathcal Q(-1,3,6)$, $\mathcal Q(-1,3,3,3)$ and  
$\mathcal Q(12)$ are non connected.
\end{NNths}

\begin{proof}
The generalized permutations 
$\left( \begin{smallmatrix} 1 & 1 & 2 & 3 & 2 & 3 & 4\\
5 & 4 & 5 & 6 & 7 & 6 & 7 \end{smallmatrix} \right)$ and 
$ \left( \begin{smallmatrix} 1 & 1 & 2 & 3 & 4 & 5 & 6\\
3 & 2 & 7 & 5 & 7 & 6 & 4 \end{smallmatrix} \right)$ are
irreducible. The corresponding suspension surfaces belong to 
the stratum $\QQQ(-1,9)$. According to Zorich's computation, these two
permutations do not belong to the same extended Rauzy classes (see Table~\ref{table:rauzy} in the Appendix). 
Hence the stratum $\QQQ(-1,9)$ is not connected. In fact this stratum has
precisely two connected components corresponding to the two extended
Rauzy classes.

\noindent We have similar conclusions for other strata with the following
generalized permutations. For the stratum $\QQQ(-1,3,6)$ one can consider 
the generalized permutations
$$
\left( \begin{smallmatrix} 1 & 1 & 2 & 3 & 2 & 3 & 4 & 5 \\
4 & 6 & 5 & 6 & 7 & 8 & 7 & 8 \end{smallmatrix} \right) \qquad \textrm{and} 
\qquad
\left( \begin{smallmatrix} 1 & 2 & 3 & 4 & 5 & 6 & 2 & 3 \\
7 & 1 & 7 & 6 & 5 & 4 & 8 & 8 \end{smallmatrix} \right).
$$
For the stratum $\QQQ(-1,3,3,3)$ one can consider 
the generalized permutations
$$
\left( \begin{smallmatrix} 1 & 1 & 2 & 3 & 4 & 5 & 6 & 7 & 6 \\
7 & 8 & 5 & 8 & 2 & 4 & 9 & 3 & 9 \end{smallmatrix} \right)\qquad \textrm{and} 
\qquad
\left( \begin{smallmatrix} 1 & 1 & 2 & 3 & 2 & 3 & 4 & 5 & 6 \\
4 & 7 & 8 & 9 & 7 & 8 & 6 & 5 & 9 \end{smallmatrix} \right).
$$
For the stratum $\QQQ(12)$ one can consider the generalized permutations
$$
\left( \begin{smallmatrix} 1 & 2 & 1 & 2 & 3 & 4 & 5 & 3 \\
6 & 7 & 6 & 7 & 5 & 8 & 4 & 8 \end{smallmatrix} \right)\qquad \textrm{and} 
\qquad
\left( \begin{smallmatrix} 1 & 2 & 3 & 4 & 5 & 6 & 7 & 6 \\
8 & 7 & 5 & 8 & 4 & 3 & 2 & 1 \end{smallmatrix} \right).
$$
The theorem is proven.
\end{proof}

\newpage
\appendix

\section{Computation of the Rauzy classes}
Here we give explicit examples of  \emph{reduced} Rauzy classes
(\emph{i.e.} up to the equivalence $\pi \sim f\circ \pi$, for any
permutation $f$ of $\mathcal{A}$). \medskip

It is easy to see that there is only one Rauzy class filled by
(irreducible) generalized permutations defined over $3$ letters. In that
case the Rauzy class contains $4$ generalized permutations and a
permutation is irreducible if and only if it is dynamically
irreducible. Thus there is no interesting phenomenon in this
``simple'' case. \medskip

\begin{figure}[htbp]
   \begin{center}
 \psfrag{1}[][]{\tiny $\left(\begin{smallmatrix} 1& 1& 2& 2\\ 3& 4& 3& 4\end{smallmatrix}\right)$}
\psfrag{2}[][]{\tiny $\left(\begin{smallmatrix} 1& 1& 2\\ 3& 2& 4& 3& 4\end{smallmatrix}\right)$}
\psfrag{3}[][]{\tiny $\left(\begin{smallmatrix} 1& 1\\ 2& 3& 3& 4& 2& 4\end{smallmatrix}\right)$}
\psfrag{4}[][]{\tiny $\left(\begin{smallmatrix} 1& 2& 2\\ 3& 4& 4& 1& 3\end{smallmatrix}\right)$}
\psfrag{5}[][]{\tiny $\left(\begin{smallmatrix} 1& 2& 3& 3\\ 2& 4& 4& 1\end{smallmatrix}\right)$}
\psfrag{6}[][]{\tiny $\left(\begin{smallmatrix} 1& 2& 3& 2\\ 3& 4& 4& 1\end{smallmatrix}\right)$}
\psfrag{7}[][]{\tiny $\left(\begin{smallmatrix} 1& 2& 2& 3\\ 3& 4& 4& 1\end{smallmatrix}\right)$}
\psfrag{8}[][]{\tiny $\left(\begin{smallmatrix} 1& 1& 2& 3& 3\\ 4& 2& 4\end{smallmatrix}\right)$}
\psfrag{9}[][]{\tiny $\left(\begin{smallmatrix} 1& 1& 2& 3\\ 3& 4& 2& 4\end{smallmatrix}\right)$}
\psfrag{10}[][]{\tiny $\left(\begin{smallmatrix} 1& 1& 2\\ 3& 3& 4& 2& 4\end{smallmatrix}\right)$}
\psfrag{11}[][]{\tiny $\left(\begin{smallmatrix} 1& 1\\ 2& 2& 3& 4& 3& 4\end{smallmatrix}\right)$}
\psfrag{12}[][]{\tiny $\left(\begin{smallmatrix} 1& 2& 2\\ 3& 3& 4& 1& 4\end{smallmatrix}\right)$}
\psfrag{13}[][]{\tiny $\left(\begin{smallmatrix} 1& 2& 3& 3\\ 4& 4& 2& 1\end{smallmatrix}\right)$}
\psfrag{14}[][]{\tiny $\left(\begin{smallmatrix} 1& 2& 1& 3& 3\\ 4& 4& 2\end{smallmatrix}\right)$}
\psfrag{15}[][]{\tiny $\left(\begin{smallmatrix} 1& 2& 1& 2& 3& 3\\ 4& 4\end{smallmatrix}\right)$}
\psfrag{16}[][]{\tiny $\left(\begin{smallmatrix} 1& 2& 1& 2& 3\\ 3& 4& 4\end{smallmatrix}\right)$}
\psfrag{17}[][]{\tiny $\left(\begin{smallmatrix} 1& 1& 2& 3\\ 3& 4& 4& 2\end{smallmatrix}\right)$}
\psfrag{18}[][]{\tiny $\left(\begin{smallmatrix} 1& 1& 2\\ 3& 2& 4& 4& 3\end{smallmatrix}\right)$}
\psfrag{19}[][]{\tiny $\left(\begin{smallmatrix} 1& 1\\ 2& 3& 2& 4& 4& 3\end{smallmatrix}\right)$}
\psfrag{20}[][]{\tiny $\left(\begin{smallmatrix} 1& 2& 2\\ 3& 1& 3& 4& 4\end{smallmatrix}\right)$}
\psfrag{21}[][]{\tiny $\left(\begin{smallmatrix} 1& 2& 3& 3\\ 4& 1& 4& 2\end{smallmatrix}\right)$}
\psfrag{22}[][]{\tiny $\left(\begin{smallmatrix} 1& 2& 2& 3& 3\\ 4& 1& 4\end{smallmatrix}\right)$}
\psfrag{1p}[][]{\tiny $\left(\begin{smallmatrix} 1& 2& 1& 2\\ 3& 3& 4& 4\end{smallmatrix}\right)$}
\psfrag{2p}[][]{\tiny $\left(\begin{smallmatrix} 1& 2& 3& 1& 3\\ 4& 4& 2\end{smallmatrix}\right)$}
\psfrag{3p}[][]{\tiny $\left(\begin{smallmatrix} 1& 2& 2& 3& 1& 3\\ 4& 4\end{smallmatrix}\right)$}
\psfrag{4p}[][]{\tiny $\left(\begin{smallmatrix} 1& 2& 2& 3& 1\\ 3& 4& 4\end{smallmatrix}\right)$}
\psfrag{5p}[][]{\tiny $\left(\begin{smallmatrix} 1& 2& 2& 3\\ 3& 1& 4& 4\end{smallmatrix}\right)$}
\psfrag{6p}[][]{\tiny $\left(\begin{smallmatrix} 1& 2& 2& 3\\ 3& 4& 1& 4\end{smallmatrix}\right)$}
\psfrag{8p}[][]{\tiny $\left(\begin{smallmatrix} 1& 2& 1\\ 3& 3& 2& 4& 4\end{smallmatrix}\right)$}
\psfrag{9p}[][]{\tiny $\left(\begin{smallmatrix} 1& 2& 3& 2\\ 4& 4& 3& 1\end{smallmatrix}\right)$}
\psfrag{10p}[][]{\tiny $\left(\begin{smallmatrix} 1& 1& 2& 3& 2\\ 4& 4& 3\end{smallmatrix}\right)$}
\psfrag{11p}[][]{\tiny $\left(\begin{smallmatrix} 1& 1& 2& 3& 2& 3\\ 4& 4\end{smallmatrix}\right)$}
\psfrag{12p}[][]{\tiny $\left(\begin{smallmatrix} 1& 1& 2& 3& 2\\ 3& 4& 4\end{smallmatrix}\right)$}
\psfrag{13p}[][]{\tiny $\left(\begin{smallmatrix} 1& 1& 2& 3\\ 3& 2& 4& 4\end{smallmatrix}\right)$}
\psfrag{14p}[][]{\tiny $\left(\begin{smallmatrix} 1& 1& 2\\ 3& 2& 3& 4& 4\end{smallmatrix}\right)$}
\psfrag{15p}[][]{\tiny $\left(\begin{smallmatrix} 1& 1\\ 2& 3& 2& 3& 4& 4\end{smallmatrix}\right)$}
\psfrag{16p}[][]{\tiny $\left(\begin{smallmatrix} 1& 2& 2\\ 3& 4& 3& 4& 1\end{smallmatrix}\right)$}
\psfrag{17p}[][]{\tiny $\left(\begin{smallmatrix} 1& 2& 2& 3\\ 4& 4& 3& 1\end{smallmatrix}\right)$}
\psfrag{18p}[][]{\tiny $\left(\begin{smallmatrix} 1& 2& 3& 3& 1\\ 4& 4& 2\end{smallmatrix}\right)$}
\psfrag{19p}[][]{\tiny $\left(\begin{smallmatrix} 1& 2& 1& 3& 3& 2\\ 4& 4\end{smallmatrix}\right)$}
\psfrag{20p}[][]{\tiny $\left(\begin{smallmatrix} 1& 2& 1& 3& 3\\ 2& 4& 4\end{smallmatrix}\right)$}
\psfrag{21p}[][]{\tiny $\left(\begin{smallmatrix} 1& 2& 1& 3\\ 2& 3& 4& 4\end{smallmatrix}\right)$}
\psfrag{22p}[][]{\tiny $\left(\begin{smallmatrix} 1& 2& 1\\ 2& 3& 3& 4& 4\end{smallmatrix}\right)$}

\psfrag{0}[][]{\tiny $\left(\begin{smallmatrix} 1& 1& 2& 2& 3\\ 4& 3& 4\end{smallmatrix}\right)$}
\psfrag{0p}[][]{\tiny $\left(\begin{smallmatrix} 1& 2& 1\\ 3& 3& 4& 4& 2\end{smallmatrix}\right)$}
\psfrag{x}[][]{\tiny 0}
\psfrag{y}[][]{\tiny 1} 
     \includegraphics[scale=0.80]{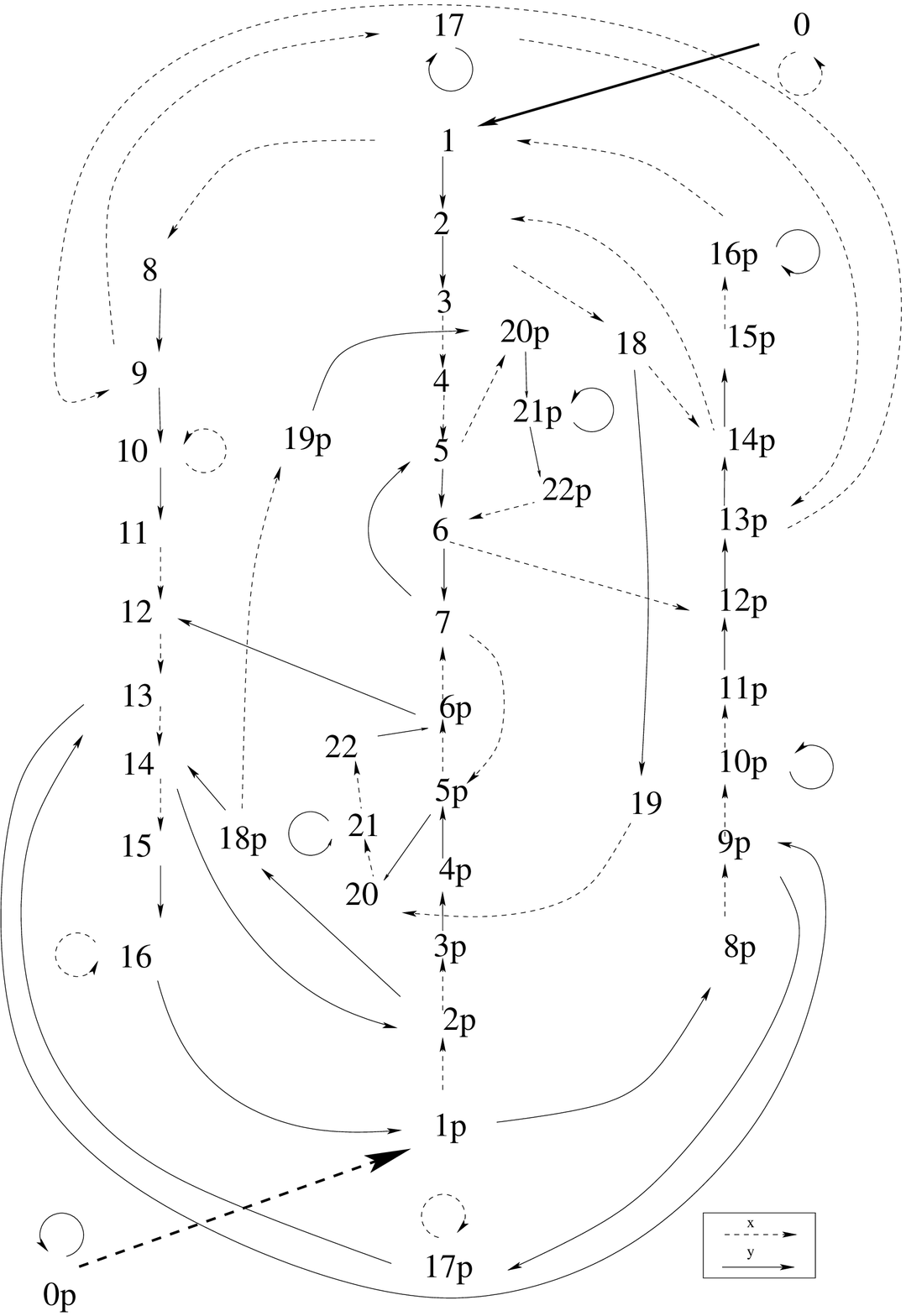}
     \caption{A (reduced) Rauzy class in $\QQQ(2,-1,-1)$.}
          \label{rauzy:class:2ii}
   \end{center}

\end{figure}

\begin{figure}[htbp]
\psfrag{1}{\tiny $\left( \begin{smallmatrix} 1 & 2 & 2 & 3 & 3 & 4 & 1\\
5 & 6 & 7 & 7 & 5 & 6 & 4 \end{smallmatrix} \right)$} 

\psfrag{2}{\tiny $\left( \begin{smallmatrix} 1 & 2 & 3 & 3 & 4 & 4 & 1 & 2\\
5 & 6 & 7 & 7 & 5 & 6 \end{smallmatrix} \right)$} 

\psfrag{3}{\tiny $\left( \begin{smallmatrix} 1 & 2 & 3 & 3 & 4 & 4 & 1\\
5 & 2 & 6 & 7 & 7 & 5 & 6 \end{smallmatrix} \right)$} 

\psfrag{4}{\tiny $\left( \begin{smallmatrix} 1 & 2 & 3 & 3 & 4 & 4\\
5 & 2 & 1 & 6 & 7 & 7 & 5 & 6 \end{smallmatrix} \right)$} 

\psfrag{5}{\tiny $ \left( \begin{smallmatrix} 1 & 2 & 3 & 4 & 4 & 5&5&2\\
6 & 3 & 1 & 7 & 7 & 6 \end{smallmatrix} \right)$}

\psfrag{6}{\tiny $ \left( \begin{smallmatrix} 1 & 2 & 3 & 4 & 5&5&6&6&3\\
2&4 & 1 & 7 & 7 \end{smallmatrix} \right)$} 

\psfrag{7}{\tiny $ \left( \begin{smallmatrix} 1 & 2 & 3 & 4&| & 5&5&6&6\\
2&4 & 1 &3&|& 7 & 7 \end{smallmatrix} \right)$} 

\psfrag{8}{\tiny $ \left( \begin{smallmatrix} 1 & 2 & 3 & 4 & 4 & 5 & 5&1&3\\
6 & 2&7 & 7 & 6 \end{smallmatrix} \right)$} 

\psfrag{9}{\tiny $ \left( \begin{smallmatrix} 1 & 2 & 3 & 4 & 5 & 5 & 1&3\\
3 & 2&7 & 7 \end{smallmatrix} \right)$}

\psfrag{10}{\tiny $ \left( \begin{smallmatrix} 1 & 2 & 3 & 4 & 5 & 5 & 1\\
3 & 2&3&7 & 7 \end{smallmatrix} \right)$} 

\psfrag{11}{\tiny $ \left( \begin{smallmatrix} 1 & 2 & 3 & 4 & 4 & 5 & 5&1\\
3&6 & 2&7 & 7 & 6 \end{smallmatrix} \right)$}

\psfrag{12}{\tiny $ \left( \begin{smallmatrix} 1 & 2 & 3 & 4 & 5 & 5&6&6&2\\
4&1 & 3&7 & 7  \end{smallmatrix} \right)$} 

\psfrag{13}{\tiny $ \left( \begin{smallmatrix} 1 & 2 & 3 & 4 & 4 & 5&5\\
2&6 & 3 & 1 & 7 & 7 & 6 \end{smallmatrix} \right)$} 

\psfrag{14}{\tiny $ \left( \begin{smallmatrix} 1 & 2 & 3 & 4 & 4 & 5 & 5\\
3&1&6 & 2&7 & 7 & 6 \end{smallmatrix} \right)$} 

\psfrag{15}{\tiny $ \left( \begin{smallmatrix} 1 & 2 & 3 & 4 &|& 5 & 5\\
3 & 2&3&1&|&7 & 7 \end{smallmatrix} \right)$} 

\psfrag{16}{\tiny $ \left( \begin{smallmatrix} 1 & 2 & 3 & 4 &|& 5 & 5&6&6\\
4&1 & 3&2&|&7 & 7  \end{smallmatrix} \right)$}

\psfrag{irred}{{\tiny $28884$ irreducible generalized permutations}}
\psfrag{red}{{\tiny $12$ reducible generalized permutations}}
\psfrag{x}[][]{\tiny 0}
\psfrag{y}[][]{\tiny 1} 

\begin{center}
\includegraphics[width=13cm]{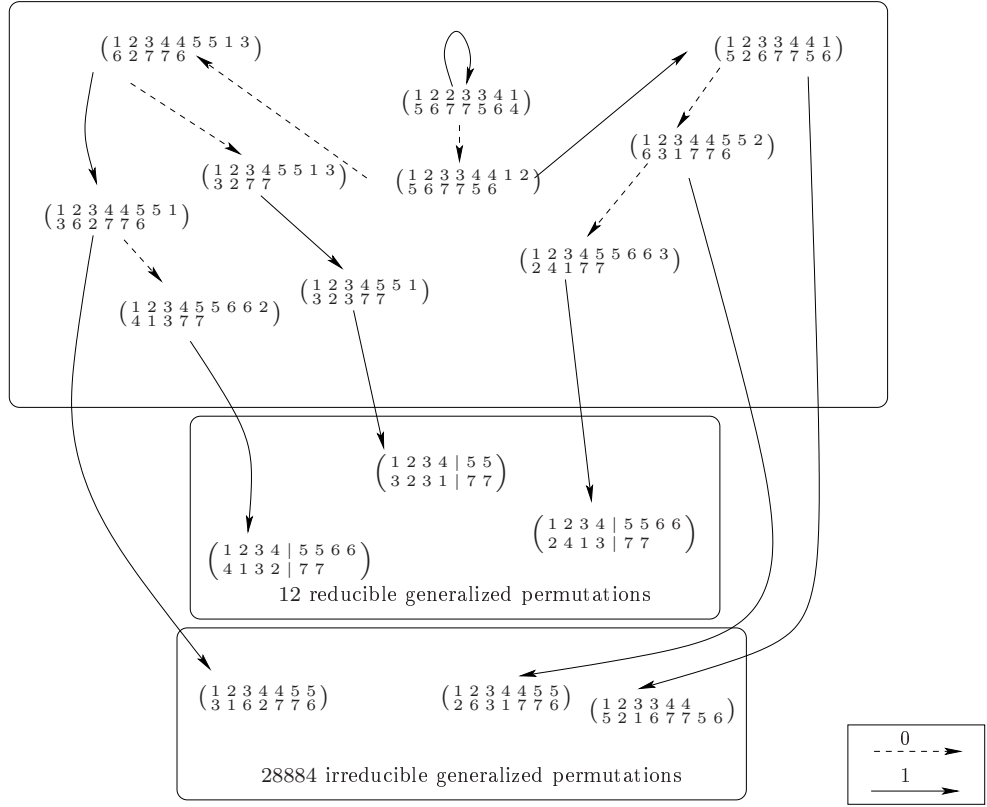}
\end{center}
\caption{
\label{fig:rauzy:2}
An example of a Rauzy class. The corresponding stratum is $\mathcal Q(-1,-1,-1,7)$. 
There are $28906$ permutations in the whole ``class'' and $28884$ permutations in the 
``good'' Rauzy class. The $28906-28884=22$ remaining permutations belong to the reducible part 
($12$ permutations) and the ``unstable'' part ($10$ permutations). Note that 
there is no smaller attractor set: the three irreducible permutations belong to the same Rauzy class.
Let us also note that the extended Rauzy class has $38456$ elements.
}
\label{classii60}
\end{figure}

If we consider a slightly more complicated case, for instance
permutations defined over $4$ letters we get some interesting phenomenon.
Figure \ref{rauzy:class:2ii} illustrates such a Rauzy class. It
corresponds to the stratum $\mathcal{Q}(2,-1,-1)$. The generalized
permutations $\left(\begin{smallmatrix} 1& 1& 2& 2& 3\\ 4& 3&
4\end{smallmatrix}\right)$ and $\left(\begin{smallmatrix} 1& 2& 1\\
3& 3& 4& 4& 2\end{smallmatrix}\right)$ are not formally in the
Rauzy class since they are reducible, but we can see there is
concretely the ``attraction'' phenomenon. As we can see the (reduced)
Rauzy classes for generalized permutations are in 
general much more complicated than the one for usual
permutation since the vertex are either of valence one or of valence
two. 
In Figure~\ref{fig:rauzy:2} we present a more complicated case with
an ``unstable'' set of permutations. \medskip

We end this section with an explicit calculation of the cardinality of the Rauzy classes
of the four exceptional strata (performed with Anton Zorich's
software~\cite{Zorich:experiments}). \medskip

\begin{table}[h]
\begin{center}
\begin{tabular}{|c|c|c|}
\hline
\textrm{connected} & \textrm{representatives} & 
\textrm{cardinality of} \\
\textrm{components} &\textrm{elements} & \textrm{extended Rauzy classes}  \\
\hline
$\QQQ(-1,9)^{adj}$ & $\left( \begin{smallmatrix} 1 & 1 & 2 & 3 & 2 & 3 & 4\\
5 & 4 & 5 & 6 & 7 & 6&  7 \end{smallmatrix} \right)$ & 95944 \\
$\QQQ(-1,9)^{irr}$ & $\left( \begin{smallmatrix} 1 & 1 & 2 & 3 & 4 & 5 & 6\\
3 & 2 & 7 & 5 & 7 & 6 & 4 \end{smallmatrix} \right)$ & 12366 \\
&& \\
$\QQQ(-1,3,6)^{adj}$ & $\left( \begin{smallmatrix} 1 & 1 & 2 & 3 & 2 & 3 & 4 & 5 \\
4 & 6 & 5 & 6 & 7 & 8 & 7 & 8 \end{smallmatrix} \right)$ & 531674 \\
$\QQQ(-1,3,6)^{irr}$ & $\left( \begin{smallmatrix} 1 & 2 & 3 & 4 & 5 & 6 & 2 & 3 \\
7 & 1 & 7 & 6 & 5 & 4 & 8 & 8 \end{smallmatrix} \right)$ & 72172 \\
&& \\
$\QQQ(-1,3,3,3)^{adj}$ & $\left( \begin{smallmatrix} 1 & 1 & 2 & 3 & 4 & 5 & 6 & 7 & 6 \\
7 & 8 & 5 & 8 & 2 & 4 & 9 & 3 & 9 \end{smallmatrix} \right)$ & 612838 \\
$\QQQ(-1,3,3,3)^{irr}$ & $\left( \begin{smallmatrix} 1 & 1 & 2 & 3 & 2 & 3 & 4 & 5 & 6 \\
4 & 7 & 8 & 9 & 7 & 8 & 6 & 5 & 9 \end{smallmatrix} \right)$ & 88374  \\
&& \\
$\QQQ(12)^{adj}$ & $\left( \begin{smallmatrix} 1 & 2 & 1 & 2 & 3 & 4 & 5 & 3 \\
6 & 7 & 6 & 7 & 5 & 8 & 4 & 8 \end{smallmatrix} \right)$ & 881599  \\
$\QQQ(12)^{irr}$ & $\left( \begin{smallmatrix} 1 & 2 & 3 & 4 & 5 & 6 & 7 & 6 \\
8 & 7 & 5 & 8 & 4 & 3 & 2 & 1 \end{smallmatrix} \right)$ & 146049 \\
\hline 
\end{tabular} 
\end{center}
\caption{
\label{table:rauzy}
Representatives elements for the special strata.}
\end{table}

\newpage

\section{An other definition of the extended Rauzy class}

In section~\ref{rauzy:classes:giem}, we have defined an extended Rauzy
class by considering the set of generalized permutations obtained from
an irreducible permutation $\gp$ by the extended Rauzy
operations. This set is not in general a subset of the irreducible
generalized permutations, therefore we must intersect it with the set
of irreducible generalized permutations to get an extended Rauzy
class.  

One could also define an extended Rauzy class in the following way: it
is a minimal subset of the irreducible generalized permutations stable
by the operations $\mathcal{R}_0$, $\mathcal{R}_1$, and $s$. It is
equivalent to say  that we forbid the operation $s$ for $\pi'$ such
that $s\pi'$ is reducible.  
For the purpose of this section, let us call this new class a~\emph{weakly
  extended Rauzy class}. A priory, an extended Rauzy class is a union
of weakly extended Rauzy classes. We will prove:

\begin{prop}
The extended Rauzy classes and the weakly extended Rauzy classes
coincide.
\end{prop}

\begin{proof}
All we have to prove is that if two irreducible generalized permutations
$\pi_1$ and $\pi_2$ correspond to the same connected component of a
stratum of quadratic differentials, then we can join them (up to
relabelling) by a combination of the maps $\mathcal{R}_0$, 
$\mathcal{R}_1$, and $s$, such that \emph{all} the corresponding
intermediary generalized permutations are irreducible. Recall that if
$\pi$ is irreducible, then so are $\mathcal{R}_0(\pi)$ and
$\mathcal{R}_1(\pi)$ (when defined). 

The idea is now to modify the proof of Theorem~$D$, by using the three
following elementary remarks. Let $\zeta$ be a suspension datum over
an irreducible generalized permutation $\pi$ (of type $(l,m)$).

\begin{enumerate}

\item In Remark~\ref{rem:strongly:irreducible} we gave a condition in
order to have $Im(\sum_{i=1}^l \zeta_{\gp(i)})=0$. Equivalently if a
decomposition of $\pi$ holds then there is no empty corner. It is
obvious to check that, under this condition, $s\gp$ is irreducible.

\item Let us assume that the two lines joining the end points of $L_0$
and the end points of $L_1$ do not have any other intersection point
with $L_0$ and $L_1$. Then applying to $\zeta_i$ the matrix $\left(
\begin{smallmatrix} 1 & 0\\ t & 1\end{smallmatrix}\right)$ for a
suitable $t$, we get a new suspension data $\zeta'$ over $T$ with 
$Im(\sum_{i=1}^l \zeta'_{\gp(i)})=0$. Hence $s\gp$ is irreducible.

\item Let $k\leq l$ minimize the value $Im(\sum_{i\leq
k}\zeta_{\pi(i)})$. Lemma~\ref{lem:sing} implies that there exists 
$n>0$ such that $\mathcal{R}^n(T)$ is the first return map of $T$ to 
the subinterval $\left(0,Re(\sum_{i=1}^k \lambda_i)\right)$. Let us 
consider $(\pi^{(n)},\zeta^{(n)})=\mathcal{R}^n(\pi,\zeta)$. By
construction $\zeta^{(n)}$ satisfies the previous condition, hence
$s\pi^{(n)}$ is irreducible.

\end{enumerate}

Let us now prove the proposition. Let $\pi_1$ and $\pi_2$ be two
generalized permutations in the same extended Rauzy class. The proof
of Theorem~$D$ asserts that there exists a surface $S$ and two
segments $X_1$ and $X_2$, each one being adjacent to a singularity
$x_1$ and $x_2$, such that for each $i$, the linear involution 
$T_i$ given by the first return maps on $X_i$ has
combinatorial datum $\gp_i$. We can assume that $S$ has no vertical
saddle connection.  

The previous remark implies that, up to replacing $T_1$ by some
$\mathcal{R}^{n_0}(T_1)$ for some well chosen $n_0$, one can and do
assume that $s\pi_1$ is irreducible. Let $(\pi_1,\zeta)$ be the
suspension over $T_1$ that corresponds to the 
surface $S$, then up to applying to $\zeta$ the matrix $\left(
\begin{smallmatrix} 1 & 0\\ t & 1\end{smallmatrix}\right)$ for a
suitable $t$ (which does not change the vertical foliation), we can
assume that $Im(\sum_{i=1}^l \zeta_{\gp(i)})=0$.

For $n$ large enough, $\mathcal{R}^n(T_2)$ is isomorphic to the first
return map on a subinterval $(y_1,y_2)$ of $X_1$, with $(y_1,0)$ or
$(y_1,1)$ a singularity of $T_1$. 
Let $k\leq l$ that minimizes the value $Im(\sum_{i\leq k}
\zeta_{\pi(i)})$ and let $x\in X_1$ be the corresponding point. If
$y_1<x$ then we also have $y_2<x$ (since $y_2$ can be chosen
arbitrarily close to $y_1$). We then apply the Rauzy-Veech induction to $T_1$
until we get a first return map on $(x_1,x)$. If $y_1>x$ then we also 
have $y_2>x$.
By definition $\zeta$ is a suspension data over $(\lambda,s\pi_1)$
(i.e. we are ``rotating by $180^\circ$'' the polygon and the
linear involution $T_1$). We apply the Rauzy-Veech induction on
$(\lambda,s\pi_1)$ until we get a  a first return map on $(x_1',x)$
that contains $y_1,y_2$. The result is a linear involution
$T_1'=(\lambda',\pi')$ such that $s\pi'$ is irreducible, and a
suspension $\zeta'$ over $T_1'$. As before we can assume that
$Im(\sum_{i=1}^l \zeta'_{\gp'(i)})=0$ and then $(\zeta',s\pi_1')$
is a suspension over $(\lambda',s\pi_1')$ that corresponds to a
first return map of $T_1$ on the subinterval $X_1'=(x,x_1')$. 
Moreover the sequence of generalized permutations joining $\pi_1$
to $s\pi_1'$ corresponding to our description consists entirely of
irreducible elements. 

Iterating this argument, there will be a step where the point $x''$
minimizing the value $Im(\sum_{i\leq k} \zeta''_{\pi''(i)})$ is
precisely $y_1$ (because the surface admits a finite number of
vertical separatrices starting from the singularities). The same 
argument produces a sequence of irreducible generalized
permutations joining $\pi''$ to $\pi_2$. 
 
This proves the equivalence of the two definitions of the extended
Rauzy classes.
\end{proof}


\end{document}